\documentclass[amscd, leqno]{amsart}
\usepackage{amscd}
\usepackage{amsmath}
\usepackage{amssymb}
\usepackage{amsfonts}
\usepackage{latexsym}
\usepackage[all]{xy}

\theoremstyle{plain}
\newtheorem{theorem}{Theorem}[section]
\newtheorem{proposition}[theorem]{Proposition}
\newtheorem{lemma}[theorem]{Lemma}

\newtheorem{corollary}[theorem]{Corollary}
\theoremstyle{definition}

\theoremstyle{remark}

\newenvironment{pf}{\begin{proof}}{\end{proof}}

\makeatletter

\@addtoreset{equation}{section}
\makeatother

\newcommand{\ZZ}{{\bf{Z}}}

\newcommand{\LL}{{\mathbb{L}_{K3}}}

\newcommand{\LAM}{{\bf\Lambda}}

\newcommand{\EE}{{\mathbb{E}}}

\newcommand{\UU}{{\mathbb{U}}}

\newcommand{\OO}{{\mathcal{O}}}
\renewcommand{\setminus}{\smallsetminus}

\renewcommand{\l}{\lambda}

\newcommand{\Ccal}{{\mathcal C}}
\newcommand{\Dcal}{{\mathcal D}}

\newcommand{\Mcal}{{\mathcal M}}

\begin{document}

\title[]
{Resultants and the Borcherds $\Phi$-function}
\author{Shu Kawaguchi}
\address{
Department of Mathematics, 
Faculty of Science,
Kyoto University, 
Kyoto 606-8502, 
Japan}
\email{kawaguch@math.kyoto-u.ac.jp}
\author{Shigeru Mukai}
\address{
Research Institute for Mathematical Sciences, 
Kyoto University, 
Kyoto 606-8502, 
Japan}
\email{mukai@kurims.kyoto-u.ac.jp}
\author{Ken-Ichi Yoshikawa}
\address{
Department of Mathematics, 
Faculty of Science,
Kyoto University, 
Kyoto 606-8502, 
Japan}
\email{yosikawa@math.kyoto-u.ac.jp}

\begin{abstract}
The Borcherds $\Phi$-function is the automorphic form on the moduli space of Enriques surfaces characterizing the discriminant locus. 
In this paper, we give an algebro-geometric construction of the Borcherds $\Phi$-function.
\end{abstract}

\maketitle

\tableofcontents


\section
{Introduction}
\label{sec:Introduction}
\par
In \cite{Borcherds96}, Borcherds discovered a beautiful automorphic form on the moduli space of Enriques surfaces,
which is closely analogous to the Dedekind $\eta$-function.
In this paper, this automorphic form is called the {\em Borcherds $\Phi$-function}. 
Up to a constant, the Borcherds $\Phi$-function is characterized as the automorphic form of weight $4$ 
on the period domain for Enriques surfaces with respect to the automorphism group of the Enriques lattice, 
whose zero divisor is exactly the discriminant locus. 
\par
Besides the original constructions \cite{Borcherds96}, \cite{Borcherds98}, several distinct understandings of the Borcherds $\Phi$-function 
are known \cite{Scheithauer00}, \cite{JorgensonKramer01}, \cite{Yoshikawa04}. 
However, even though its direct connection with the moduli space of Enriques surfaces,
no explicit algebro-geometric construction of the Borcherds $\Phi$-function has been known until now. 
The purpose of the present paper is to give such a construction modeled after the following classical result
for the Dedekind $\eta$-function: For an elliptic curve $y^{2}=4x^{3}-g_{2}x-g_{3}$ equipped with a symplectic basis $\{\alpha,\beta\}$ 
of its first integral homology group, the value of the Dedekind $\eta$-function evaluated at its period 
$\tau_{\alpha,\beta}=(\int_{\beta}dx/y)/(\int_{\alpha}dx/y)$, $\Im\tau_{\alpha,\beta}>0$ is given by the formula
\begin{equation}
\label{eqn:Dedekind:eta}
\eta\left(\tau_{\alpha,\beta}\right)^{24}
=
\left(g_{2}^{3}-27g_{3}^{2}\right)\left(\frac{1}{2\pi}\int_{\alpha}\frac{dx}{y}\right)^{12}.
\end{equation}
\par
Since the period of an Enriques surface, i.e., the period of its universal covering $K3$ surface, does depend on the choice of a marking, 
the value of the Borcherds $\Phi$-function makes sense only for {\em marked } Enriques surfaces. 
However, the Petersson norm of the Borcherds $\Phi$-function is independent of the choice of a marking and hence is well defined 
even for non-marked Enriques surfaces. 
For an Enriques surface $Y$, the Petersson norm of the Borcherds $\Phi$-function evaluated at its period point is denoted by $\|\Phi(Y)\|$. 
If $Y$ has rational double points (RDP's for short), then $\|\Phi(Y)\|$ is defined as the value $\|\Phi(\widetilde{Y})\|$, 
where $\widetilde{Y}\to Y$ is the minimal resolution.
\par
Firstly, we give an algebro-geometric expression of $\|\Phi(Y)\|$.
For this, we need explicit equations defining Enriques surfaces.
Let $f_{1},g_{1},h_{1}\in{\bf C}[x_{1},x_{2},x_{3}]$ and $f_{2},g_{2},h_{2}\in{\bf C}[x_{4},x_{5},x_{6}]$ 
be quadratic forms and define $f,g,h\in{\bf C}[x_{1},\ldots,x_{6}]$ as
$$
f=f_{1}+f_{2}, 
\qquad\quad
g=g_{1}+g_{2}, 
\qquad\quad
h=h_{1}+h_{2}.
$$
We set
$$
X_{(f,g,h)}
=
\{[x]\in{\bf P}^{5};\,f(x)=g(x)=h(x)=0\},
$$
which is preserved by the involution on ${\bf P}^{5}$
$$
\iota(x_{1},x_{2},x_{3},x_{4},x_{5},x_{6})=(x_{1},x_{2},x_{3},-x_{4},-x_{5},-x_{6}).
$$
If $(f,g,h)$ is not special, then $X_{(f,g,h)}$ is a $K3$ surface on which acts $\iota$ freely. Hence
$$
Y_{(f,g,h)}
=
X_{(f,g,h)}/\iota
$$
is an Enriques surface. In fact, every Enriques surface is expressed as the minimal resolution of some $Y_{(f,g,h)}$
by Verra \cite{Verra83} and Cossec \cite{Cossec85}.
\par
We define a canonical differential $\omega$ on $X_{(f,g,h)}$ as the residue of $(f,g,h)$ (cf. \cite{GriffithsSchmid75}).
Namely, let $F,G,H\in{\bf C}[x_{1},\ldots,x_{5}]$ be the inhomogeneous equations of $f,g,h$ obtained by setting $x_{6}=1$, respectively.
The meromorphic canonical form on ${\bf C}^{5}$
$$
\frac{dx_{1}\wedge\cdots\wedge dx_{5}}{F(x)G(x)H(x)}
$$
extends to a meromorphic canonical form on ${\bf P}^{5}$ with logarithmic poles along the divisor ${\rm div}(fgh)\subset{\bf P}^{5}$.
We set
$$
\omega=\Upsilon|_{X_{(f,g,h)}},
$$
where $\Upsilon$ is a locally defined holomorphic $2$-form on ${\bf P}^{5}$ such that
$$
\frac{dF}{F}\wedge\frac{dG}{G}\wedge\frac{dH}{H}\wedge\Upsilon
=
\frac{dx_{1}\wedge\cdots\wedge dx_{5}}{FGH}.
$$
Then $\omega=\Upsilon|_{X_{(f,g,h)}}$ is independent of the choice of $\Upsilon$ as above.
\par
Our first theorem is stated as follows (cf. Theorem~\ref{thm:valuePhi:general}).

\begin{theorem}
\label{thm:MainTheorem1}
If $Y_{(f,g,h)}$ is an Enriques surface with possibly RDP's, then
\begin{equation}
\label{eqn:MainTheorem1}
\left\|\Phi\left(Y_{(f,g,h)}\right)\right\|^{2}
=
\left|
R(f_{1},g_{1},h_{1})R(f_{2},g_{2},h_{2})
\right|
\left(
\frac{2}{\pi^{4}}
\int_{X_{(f,g,h)}}\omega\wedge\overline{\omega}
\right)^{4},
\end{equation}
where $R(f_{i},g_{i},h_{i})$ is the resultant of $f_{i},g_{i},h_{i}$.
\end{theorem}

We remark that a weaker version of \eqref{eqn:MainTheorem1} was obtained by Maillot-R\"ossler \cite{MaillotRoessler08}.
As a consequence of Theorem~\ref{thm:MainTheorem1}, we get an extension of \eqref{eqn:Dedekind:eta} to 
the Borcherds $\Phi$-function (cf. Theorem~\ref{theorem:Thomae:Phi}).
Write $\langle\cdot,\cdot\rangle$ for the cup-product on $H^{2}(X_{(f,g,h)},{\bf Z})$.

\begin{theorem}
\label{thm:MainTheorem2}
Assume that $Y_{(f,g,h)}$ is smooth.
Let ${\bf v}\in H^{2}(X_{(f,g,h)},{\bf Z})$ be an anti-$\iota$-invariant primitive isotropic vector of level $\ell\in\{1,2\}$ 
and let ${\bf v}^{\lor}\in H_{2}(X_{(f,g,h)},{\bf Z})$ be its Poincar\'e dual.
Let ${\bf v}'\in H^{2}(X_{(f,g,h)},{\bf Z})$ be another anti-$\iota$-invariant primitive isotropic vector of level $\ell$ 
with $\langle{\bf v},{\bf v}'\rangle=\ell$. Then, by choosing a suitable marking $\alpha_{{\bf v},{\bf v}'}$ of $X_{(f,g,h)}$, the vector
$$
z_{{\bf v},{\bf v}'}
=
\frac{
\omega
-
\langle\omega,{\bf v}'/\ell\rangle{\bf v}
-
\langle\omega,{\bf v}\rangle({\bf v}'/\ell)
}
{\langle\omega,{\bf v}\rangle}
\in({\bf Z}{\bf v}+{\bf Z}{\bf v}')^{\perp}\otimes{\bf C}
$$ 
is the period of $Y_{(f,g,h)}$ via $\alpha_{{\bf v},{\bf v}'}$ (cf. Sect.~ \ref{sect:Thomae:Phi}) and the following equality holds
\begin{equation}
\label{eqn:MainTheorem2}
\Phi_{\ell}\left(\alpha_{{\bf v},{\bf v}'}(z_{{\bf v},{\bf v}'})\right)^{2} 
=
R(f_{1},g_{1},h_{1})R(f_{2},g_{2},h_{2})
\left(
\frac{2}{\pi^{2}}
\int_{{\bf v}^{\lor}}\omega
\right)^{8}.
\end{equation}
Here $\Phi_{\ell}(z)$ is the Borcherds $\Phi$-function with respect to the level $\ell$ cusp (cf. Sect.~\ref{subsec:summary:Borcherds:Phi}).
\end{theorem}
We point out a similarity of \eqref{eqn:MainTheorem2} to the Thomae type formula of Matsumoto-Terasoma \cite{MatsumotoTerasoma12}
for certain $K3$ surfaces.
\par
Let $M_{3,6}^{o}({\bf C})$ be the set of complex $3\times 6$-matrices without vanishing $3\times3$-minors.
Considering the loci of the moduli space of Enriques surfaces parametrized by $M_{3,6}^{o}({\bf C})$, 
we get an explicit relation between the Borcherds $\Phi$-function and certain theta functions.
For $N=(n_{ij})=({\bf n}_{1},\ldots,{\bf n}_{6})\in M_{3,6}^{o}({\bf C})$, we define
$$
X_{N}
=
\left\{
[x]\in{\bf P}^{5};\,
{\bf n}_{1}x_{1}^{2}+{\bf n}_{2}x_{2}^{2}+{\bf n}_{3}x_{3}^{2}+{\bf n}_{4}x_{4}^{2}+{\bf n}_{5}x_{5}^{2}+{\bf n}_{6}x_{6}^{2}={\bf 0}
\right\}.
$$
Then the involution $\iota_{\binom{pqr}{stu}}(x_{p},x_{q},x_{r},x_{s},x_{t},x_{u})=(x_{p},x_{q},x_{r},-x_{s},-x_{t},-x_{u})$ acts on $X_{N}$.
Here $\binom{pqr}{stu}$ denotes the partition $\{p,q,r\}\amalg\{s,t,u\}=\{1,\ldots,6\}$.
For all $N\in M_{3,6}^{o}({\bf C})$ and for all partitions $\binom{pqr}{stu}$,  
$X_{N}/\iota_{\binom{pqr}{stu}}$ is an Enriques surface.
\par
There is another $K3$ surface $Z_{N}$ associated to $N\in M_{3,6}^{o}({\bf C})$
$$
Z_{N}
=
\{
((w_{1}:w_{2}:w_{3}),y)\in{\mathcal O}_{{\bf P}^{2}}(3);\,
y^{2}=\prod_{k=1}^{6}(n_{1k}w_{1}+n_{2k}w_{2}+n_{3k}w_{3})
\},
$$
which is identified with its minimal resolution $\widetilde{Z}_{N}$.
After a suitable choice of a system of transcendental cycles of $\widetilde{Z}_{N}$,
the period of $\widetilde{Z}_{N}$, denoted by $\varOmega_{N}$, lies in
$$
{\Bbb D}=\{T\in M_{2}({\bf C});\,(T-{}^{t}\overline{T})/2\sqrt{-1}>0\}\subset M_{2}({\bf C}),
$$ 
where $M_{2}({\bf C})$ is the set of complex $2\times 2$-matrices 
(cf. \cite{Yoshida97}, \cite{MatsumotoTerasoma12} and Sect.\,\ref{subsect:Phi:Freitag:theta}).
\par
For $\varOmega\in{\Bbb D}$ and $\binom{pqr}{stu}$, the theta function
$\Theta_{\binom{pqr}{stu}}(\varOmega)$ is defined as a certain Fourier series on ${\Bbb D}$ 
(cf. Sect.~\ref{subsect:Phi:Freitag:theta}), whose Petersson norm is denoted by $\|\Theta_{\binom{pqr}{stu}}(\varOmega)\|$.

\begin{theorem}
\label{thm:MainTheorem3}
For all $N\in M^{o}(3,6)$, the following equality holds
\begin{equation}
\label{eqn:MainTheorem3}
\left\|
\Phi\left(X_{N}/\iota_{\binom{pqr}{stu}}\right)
\right\|
=
\left\|
\Theta_{\binom{pqr}{stu}}(\varOmega_{N})
\right\|^{4}.
\end{equation}
\end{theorem}

Since the Borcherds $\Phi$-function is expressed as an additive Borcherds lift by Freitag-Salvati-Manni \cite{FreitagSalvatiManni07}, 
it may not be surprising that the restriction of the Borcherds $\Phi$-function to a subdomain isomorphic to ${\Bbb D}$ is expressed as a theta series.
As a by-product of Theorem~\ref{thm:MainTheorem3}, we shall show that $\theta_{a,b}(T)^{8}$ is expressed as an infinite product of Borcherds type,
where $\theta_{a,b}(T)$ is an arbitrary even theta constant of genus $2$ (cf. Corollary~\ref{cor:infinite:product:theta:genus:2}).
We remark that the product of {\em all} even theta constants of genus $2$ is expressed as a Borcherds product 
by Gritsenko-Nikulin \cite{GritsenkoNikulin97}.
\par
To prove Theorem~\ref{thm:MainTheorem1}, we compare the $\partial\bar{\partial}\log(\cdot)$ of the both sides of \eqref{eqn:MainTheorem1} 
as currents on the Grassmann variety parametrizing $3$-dimensional subspaces of the $\iota$-invariant quadratic forms 
in the variables $x_{1},\ldots,x_{6}$. Outside the locus of vanishing resultants, the comparison of the curvature is easy.
We give a criterion that the period map for a one-parameter degenerating family of Enriques surfaces intersects
the discriminant locus transversally, and we use it to determine the singularity of $\|\Phi\|$ near the locus of vanishing resultants. 
In this way, we prove Theorem~\ref{thm:MainTheorem1}, up to a constant. 
To determine the constant, we evaluate the both sides of \eqref{eqn:MainTheorem1} 
for those Enriques surfaces studied by Mukai \cite{Mukai11} and Ohashi \cite{Ohashi07}. 
\par
This paper is organized as follows.
In Section $2$, we recall Enriques surfaces and Borcherds $\Phi$-function.
In Section $3$, we prove Theorem~\ref{thm:MainTheorem1}.
In Section $4$, we prove Theorem~\ref{thm:MainTheorem2}.
In Section $5$, we prove Theorem~\ref{thm:MainTheorem3}. 
In Section $6$, we study infinite product expansion of theta constants of genus $2$.
In Section $7$, we study the Borcherds $\Phi$-function for those Enriques surfaces studied by Mukai and Ohashi.
\par
{\bf Acknowledgements }
The first author is partially supported by JSPS Grants-in-Aid for young scientists (B) 24740015
and was supported by JSPS-MAEE bilateral program Sakura.
The second author is partially supported by JSPS Grants-in-Aid (S) 25220701, (A) 22244003. 
The third author is partially supported by JSPS Grants-in-Aid (B) 23340017, (A) 22244003, (S) 22224001, (S) 25220701.

\section
{Enriques surfaces and the Borcherds $\Phi$-function}
\label{sec:Enriques:Phi}
\par 
A free ${\bf Z}$-module of finite rank endowed with a non-degenerate, 
integral, symmetric bilinear form is called a lattice. 
For a lattice $L=({\bf Z}^{r},\langle\cdot,\cdot\rangle)$ and $k\in{\bf Q}^{\times}$, 
we define $L(k):=({\bf Z}^{r},k\langle\cdot,\cdot\rangle)$.
The group of isometries of $L$ is denoted by $O(L)$.
The set of roots of $L$ is defined by $\Delta_{L}:=\{l\in L;\,\langle l,l\rangle=-2\}$.
\par
Let ${\Bbb U}=({\bf Z}^{2},\binom{0\,1}{1\,0})$ and let ${\Bbb E}_{8}$ be the even unimodular
{\em negative-definite} lattice of rank $8$. 
The $K3$-lattice ${\Bbb L}_{K3}$ and the {\em Enriques lattice} $\LAM$ are defined as
$$
{\Bbb L}_{K3}
:=
{\Bbb U}\oplus{\Bbb U}\oplus{\Bbb U}\oplus
{\Bbb E}_{8}\oplus{\Bbb E}_{8},
\qquad
\LAM := \UU(2) \oplus \UU \oplus \EE_8(2).
$$
We fix a primitive embedding $\LAM \subset \LL$.
Then $\LAM^{\perp_{\LL}} \cong \UU(2)\oplus\EE_8(2)$.

\subsection
{$K3$ and Enriques surfaces}
\par
A {\em $K3$ surface} $X$ is a smooth compact complex surface with $H^1(X, \OO_X) = 0$ and 
the trivial canonical line bundle $K_{X}\cong{\mathcal O}_{X}$. By \cite[p.241]{BPV84}, 
$H^2(X, \ZZ)$ endowed with the cup-product is isometric to $\LL$. 
\par
An {\em Enriques surface} $Y$ is a compact connected complex surface with 
$H^1(Y, \OO_Y) = H^2(Y, \OO_Y) =0$, $K_{Y}\not\cong{\mathcal O}_{X}$ and 
$K_{Y}^{\otimes2} \cong{\mathcal O}_{X}$. 
Let $X$ be the universal covering of $Y$ and 
let $\theta\colon X\to X$ be the non-trivial covering transformation of $X$ over $Y$. 
Then $X$ is a $K3$ surface and $\theta$ is an anti-symplectic holomorphic involution without fixed points. 
By e.g. \cite{BPV84}, there is an isometry of lattices $\alpha\colon H^{2}(X,{\bf Z})\cong{\Bbb L}_{K3}$ such that 
\begin{equation}
\label{eqn:alpha:+:-}
\alpha(H^{2}(X,{\bf Z})_{+})
=
\LAM^{\perp_{\LL}},
\qquad
\alpha(H^{2}(X,{\bf Z})_{-})
=
\LAM,
\end{equation}
where $H^{2}(X,{\bf Z})_{\pm}:=\{l\in H^{2}(X,{\bf Z});\,\theta^{*}(l)=\pm l\}$.
The pair $(Y,\alpha)$ is called a marked Enriques surface
if the isometry $\alpha$ satisfies \eqref{eqn:alpha:+:-}.
\par
We set
$$
\Omega_{\LAM}
:=
\{[\eta]\in{\bf P}(\LAM\otimes{\bf C});\,
\langle\eta,\eta\rangle=0,\,\langle\eta,\bar{\eta}\rangle>0\}.
$$
Let $\Omega_{\LAM}^{+}$ be one of the connected components of $\Omega_{\LAM}$, 
which we give in Sect.\,\ref{subsec:summary:Borcherds:Phi}.
Then $\Omega_{\LAM}^{+}$ is a bounded symmetric domain of type IV of dimension $10$. 
Set 
$$
O^{+}(\LAM):=\{g \in O(\LAM); g(\Omega_{\LAM}^{+}) = \Omega_{\LAM}^{+}\}.
$$
The projective action of $O^{+}(\LAM)$ on $\Omega_{\LAM}^{+}$ is proper and discontinuous. We define
$$
{\mathcal M}:=\Omega_{\LAM}^{+}/O^{+}(\LAM). 
$$
The period of a marked Enriques surface $(Y,\alpha)$ is defined as
$$
\varpi(Y,\alpha):=[\alpha(H^{0}(X,\Omega^{2}_{X}))]\in\Omega_{\LAM}^{+}
$$
and the period of an Enriques surface $Y$ is defined as the $O^{+}(\LAM)$-orbit of $\varpi(Y,\alpha)$:
$$
\overline{\varpi}(Y):=[\varpi(Y,\alpha)]\in{\mathcal M}.
$$
For $d \in\Delta_{\LAM}$, set $H_{d}:= \{\eta \in\Omega_{\LAM}^{+}; \langle\eta, d\rangle = 0\}$. 
The $O^{+}(\LAM)$-invariant reduced divisor
$$
{\mathcal D}
:=
\sum_{d\in\Delta_{\LAM}/\pm1}H_{d}
\subset
\Omega_{\LAM}^+
$$
is called the discriminant locus. We set $\overline{\mathcal D}:={\mathcal D}/O^{+}(\LAM)$.
Then $\overline{\varpi}(Y)\not\in\overline{\mathcal D}$. 
By \cite{BPV84}, 
$$
{\mathcal M}^{o}
:=
(\Omega_{\LAM}^{+}\setminus{\mathcal D})/O^{+}(\LAM)
={\mathcal M}\setminus\overline{\mathcal D}
$$ 
is the coarse moduli space of Enriques surfaces via the period map.

\subsection
{The Borcherds $\Phi$-function and its basic properties}
\label{subsec:summary:Borcherds:Phi}
\par

\subsubsection{The Borcherds $\Phi$-function}
\label{subsubsec:Borcherds:Phi:Enriques:lattice}
\par
For a subset $S\subset{\bf P}(\LAM\otimes{\bf C})$, the cone over $S$ is denoted by $C(S):=\{\eta\in(\LAM\otimes{\bf C})\setminus\{0\};\,[\eta]\in S\}$.
Up to a constant, the Borcherds $\Phi$-function is defined as the holomorphic function $\Phi(Z)$ on $C(\Omega_{\LAM}^{+})$ 
with the following properties (C1), (C2), (C3):
\begin{itemize}
\item[(C1)]
$\Phi(\lambda Z)=\lambda^{-4}\,\Phi(Z)$ for all $\lambda\in{\bf C}^{*}:={\bf C}\setminus\{0\}$.
\item[(C2)]
$\Phi(g(Z))=\chi(g)\,\Phi(Z)$ for all $g\in O^{+}(\LAM)$, where $\chi\in{\rm Hom}(O^{+}(\LAM),{\bf C}^{*})$.
\item[(C3)]
The zero divisor of $\Phi$ is the cone $C({\mathcal D})$.
\end{itemize}
In \cite{Borcherds96}, $\Phi(Z)$ (precisely speaking $\Phi_{2}(w)$ below) was constructed.
In \cite[Example 13.7]{Borcherds98}, $\Phi(Z)$ (precisely speaking $\Phi_{1}(z)$, $\Phi_{2}(w)$ below) was constructed as the Borcherds lift 
of certain vector valued modular form for $Mp_{2}({\bf Z})$. 
Since the lattice used in \cite[Example 13.7]{Borcherds98} is distinct from $\LAM$,
we also refer to \cite[Example 8.9]{Yoshikawa13}, where $\LAM$ is used in the construction of $\Phi$ as the Borcherds lift.
\par
The Petersson norm of $\Phi(Z)$ is the $C^{\infty}$ function on $C(\Omega_{\LAM}^{+})$ defined by
$$
\|\Phi(Z)\|^{2}:=2^{-4}\langle Z,\overline{Z}\rangle_{\LAM}^{4}|\Phi(Z)|^{2}.
$$
Since $\chi^{2}$ is trivial by Lemma~\ref{cor:character:Phi:square} below, $\|\Phi(Z)\|^{2}$ is invariant under the actions of the groups
$O^{+}(\LAM)$ and ${\bf C}^{*}$. We regard $\|\Phi\|^{2}$ as a $C^{\infty}$ function on ${\mathcal M}$ in what follows.
For an Enriques surface $Y$, we define
$$
\|\Phi(Y)\|:=\|\Phi(\overline{\varpi}(Y))\|.
$$
\par
Let $\kappa$ be the K\"ahler form of the Bergmann metric on $\Mcal$ and 
let $\delta_{\overline{\Dcal}}$ be the Dirac $\delta$-current associated to the divisor $\overline{\Dcal}$.
By the Poincar\'e-Lelong formula, we get 
\begin{equation}
\label{eqn:Poincare:Lelong:Phi:moduli}
-dd^{c}\log\|\Phi\|^{2}=4\kappa-\frac{1}{2}\delta_{\overline{\Dcal}}
\end{equation}
as currents on $\Mcal$.
Here the coefficient $1/2$ enters in the formula, because the projection $\Omega_{\LAM}^{+}\to\Mcal$ is doubly
ramified along $\Dcal$.

\begin{lemma}
\label{cor:character:Phi:square}
The character of $\Phi^{2}$ is trivial.
\end{lemma}

\begin{pf}
Since $O^{+}(\LAM^{\lor}(2))\cong O^{+}(\LAM^{\lor})\cong O^{+}(\LAM)$, the result follows from the fact that 
$O^{+}(\LAM^{\lor}(2))$ is generated by reflections \cite[proof of Prop.\,5.6]{GritsenkoHulekSankaran09}.
\end{pf}

\subsubsection{The tube domain realization of $\Omega_{\LAM}$ and the automorphic factor}
\label{subsubsec:tubedomain:automorphic:factor}
\par
We define the positive cone of a Lorentzian lattice $L$ as ${\mathcal C}_{L}=\{x\in L\otimes{\bf R};\,x^{2}>0\}$.
Then ${\mathcal C}_{L}$ consists of two connected components ${\mathcal C}_{L}^{+}$ and ${\mathcal C}_{L}^{+}$
such that ${\mathcal C}_{L}^{-}=-{\mathcal C}_{L}^{+}$.
\par
For $\ell=\{1,2\}$, set 
$$
{\Bbb M}_{\ell}:={\Bbb U}(2/\ell)\oplus{\Bbb E}_{8}(2).
$$ 
Let $\{{\frak e}_{\ell},{\frak f}_{\ell}\}$ be a basis of ${\Bbb U}(\ell)$ such that 
$\langle{\frak e}_{\ell},{\frak e}_{\ell}\rangle=\langle{\frak f}_{\ell},{\frak f}_{\ell}\rangle=0$
and $\langle{\frak e}_{\ell},{\frak f}_{\ell}\rangle=\ell$. Regarding ${\Bbb U}$ and ${\Bbb U}(2)$ as direct summands of $\LAM$, 
we get ${\frak e}_{\ell}^{\perp}/{\bf Z}{\frak e}_{\ell}\cong{\Bbb M}_{\ell}$.
\par
We define the isomorphism 
$\iota_{\ell}\colon{\Bbb M}_{\ell}\otimes{\bf R} +i\,{\mathcal C}_{{\Bbb M}_{\ell}}\ni u\to[\iota_{\ell}(u)]\in\Omega_{\LAM}$ 
by
\begin{equation}
\label{eqn:tube:domain:level1}
\iota_{\ell}(u)
:=
-(u^{2}/2){\frak e}_{\ell}+({\frak f}_{\ell}/\ell)+(-1)^{2/\ell}u.
\end{equation}
Let ${\mathcal C}_{{\Bbb M}_{1}}^{+}$ be the component of ${\mathcal C}_{{\Bbb M}_{1}}$ whose closure contains ${\frak e}_{2}$, ${\frak f}_{2}$.
Then $\Omega_{\LAM}^{+}$ is defined as the component corresponding to ${\Bbb M}_{1}\otimes{\bf R}+i\,{\mathcal C}_{{\Bbb M}_{1}}^{+}$,
and ${\mathcal C}_{{\Bbb M}_{2}}^{+}$ is defined as the component such that
${\Bbb M}_{2}\otimes{\bf R}+i\,{\mathcal C}_{{\Bbb M}_{2}}^{+}\cong\Omega_{\LAM}^{+}$ via \eqref{eqn:tube:domain:level1}.
Replacing $\{{\frak e}_{1},{\frak f}_{1}\}$ by $\{-{\frak e}_{1},-{\frak f}_{1}\}$ if necessary, we may and will assume that
${\frak e}_{1},{\frak f}_{1}\in\overline{\mathcal C}_{{\Bbb M}_{2}}^{+}$.
\par
Through the isomorphism \eqref{eqn:tube:domain:level1}, $O^{+}(\LAM)$ acts on ${\Bbb M}_{\ell}\otimes{\bf R}+i\,{\mathcal C}_{{\Bbb M}_{\ell}}^{+}$. 
For $u\in{\Bbb M}_{\ell}\otimes{\bf R}+i\,{\mathcal C}_{{\Bbb M}_{\ell}}^{+}$ and $g\in O^{+}(\LAM)$, 
we define $g\cdot u\in{\Bbb M}_{\ell}\otimes{\bf R}+i\,{\mathcal C}_{{\Bbb M}_{\ell}}^{+}$ by the formula
\begin{equation}
\label{eqn:group:action:domain:type:IV:level1}
\iota_{\ell}(g\cdot u)
=
-\frac{(g\cdot u)^2}{2}{\frak e}_{\ell}
+
\frac{{\frak f}_{\ell}}{\ell}
+
(-1)^{2/\ell}g\cdot u
:=
\frac{g(\iota_{\ell}(u))}{\langle g(\iota_{\ell}(u)),{\frak e}_{\ell}\rangle_{\LAM}}.
\end{equation}
\par
The automorphic factor for the $O^{+}(\LAM)$-action on ${\Bbb M}_{\ell}\otimes{\bf R}+i\,{\mathcal C}_{{\Bbb M}_{\ell}}^{+}$ is defined by
\begin{equation}
\label{eqn:automorphic:factor:domain:type:IV:level1}
j_{\ell}(g,u)
:=
\langle g(\iota_{\ell}(u)),{\frak e}_{\ell}\rangle_{\LAM}
=
\langle g(-(u^{2}/2){\frak e}_{\ell}+({\frak f}_{\ell}/\ell)+(-1)^{2/\ell}u),{\frak e}_{\ell}\rangle_{\LAM}.
\end{equation}
Then $j_{\ell}(gg',u)=j_{\ell}(g,g'\cdot u)j_{\ell}(g',u)$ for all $g,g'\in O^{+}(\LAM)$.
Since 
\begin{equation}
\label{eqn:relation:bergman:tube:bounded}
\begin{aligned}
2\langle\Im u,\Im u\rangle_{{\Bbb M}_{\ell}}
&=
\langle
-(u^{2}/2){\frak e}_{\ell}+({\frak f}_{\ell}/\ell)+(-1)^{2/\ell}u,
\overline{-(u^{2}/2){\frak e}_{\ell}+({\frak f}_{\ell}/\ell)+(-1)^{2/\ell}u}
\rangle_{\LAM}
\\
&=
\langle\iota_{\ell}(u),\overline{\iota_{\ell}(u)}\rangle_{\LAM},
\end{aligned}
\end{equation}
we get the automorphic property of $\langle\Im u,\Im u\rangle_{{\Bbb M}_{\ell}}$
by \eqref{eqn:group:action:domain:type:IV:level1}, \eqref{eqn:automorphic:factor:domain:type:IV:level1}
$$
\langle\Im(g\cdot u),\Im(g\cdot u)\rangle_{{\Bbb M}_{\ell}}=|j_{\ell}(g,u)|^{-2}\langle\Im u,\Im u\rangle_{{\Bbb M}_{\ell}}
\qquad
(g\in O^{+}(\LAM)).
$$

\subsubsection{The Borcherds $\Phi$-function with respect to the level $\ell$ cusp}
\label{subsubsec:Petersson:norm:automorphic:property}
\par 
We define the Borcherds $\Phi$-function with respect to the level $\ell$ cusp as the pullback of $\Phi(Z)$ to
${\Bbb M}_{\ell}\otimes{\bf R}+i\,{\mathcal C}_{{\Bbb M}_{\ell}}^{+}$ via the embedding $\iota_{\ell}$, i.e.,
$$
\Phi_{\ell}(u)
:=
\Phi(\iota_{\ell}(u))
=
\Phi(-(u^{2}/2){\frak e}_{\ell}+({\frak f}_{\ell}/\ell)+(-1)^{2/\ell}u).
$$
By (C1), $\Phi(Z)$ can be recovered from $\Phi_{\ell}(u)$. In this sense, $\Phi$ and $\Phi_{\ell}$ are equivalent.
\par
By \eqref{eqn:group:action:domain:type:IV:level1}, \eqref{eqn:automorphic:factor:domain:type:IV:level1}, (C1), (C2), 
$\Phi_{\ell}(u)$ satisfies the following functional equation 
\begin{equation}
\label{eqn:automorphic:property:Phi:level:N}
\Phi_{\ell}(g\cdot u)
=
\chi(g)\,j_{\ell}(g,u)^{4}\,\Phi_{\ell}(u)
\end{equation}
on ${\Bbb M}_{\ell}\otimes{\bf R}+i\,{\mathcal C}_{{\Bbb M}_{\ell}}^{+}$ for all $g\in O^{+}({\LAM})$.
\par
By \eqref{eqn:relation:bergman:tube:bounded}, the Petersson norm of $\Phi$ is expressed as follows on 
${\Bbb M}_{\ell}\otimes{\bf R}+i\,{\mathcal C}_{{\Bbb M}_{\ell}}^{+}$
$$
\|\Phi_{\ell}(u)\|
:=
\|\Phi(\iota_{\ell}(u))\|^{2}
=
\langle{\Im u},{\Im u}\rangle_{{\Bbb M}_{\ell}}^{4}|\Phi_{\ell}(u)|^{2}.
$$
\par
The relation between $\Phi_{1}$ and $\Phi_{2}$ is as follows.
For $z\in{\Bbb M}_{1}\otimes{\bf R}+i\,{\mathcal C}_{{\Bbb M}_{1}}^{+}$, 
the point $w\in{\Bbb M}_{2}\otimes{\bf R}+i\,{\mathcal C}_{{\Bbb M}_{2}}^{+}$ corresponding to $z$ is given by
$\iota_{1}(z)/\langle z,{\frak e}_{2}\rangle_{{\Bbb M}_{1}}=\iota_{2}(w)$. Hence
\begin{equation}
\label{eqn:transformation}
w
=
-\frac{1}{\langle z,{\frak e}_{2}\rangle_{{\Bbb M}_{1}}}
\left\{
-\frac{z^{2}}{2}{\frak e}_{1}
+
{\frak f}_{1}
+
\left(
z-\langle z,{\frak f}_{2}\rangle_{{\Bbb M}_{1}}\frac{{\frak e}_{2}}{2}-\langle z,{\frak e}_{2}\rangle_{{\Bbb M}_{1}}\frac{{\frak f}_{2}}{2}
\right)
\right\}.
\end{equation}
By the inequality $\langle{\Im w},{\frak e}_{1}\rangle_{{\Bbb M}_{2}}=\Im\left(-1/\langle z,{\frak e}_{2}\rangle_{{\Bbb M}_{1}}\right)>0$, 
we get $\Im w\in{\mathcal C}_{{\Bbb M}_{2}}^{+}$.
Since $\iota_{1}(z)/\langle z,{\frak e}_{2}\rangle_{{\Bbb M}_{1}}=\iota_{2}(w)$, we deduce from (C1) that
\begin{equation}
\label{eqn:relation:Phi1:Phi2:sec2.2}
\Phi_{2}(w)
=
\langle z,{\frak e}_{2}\rangle_{{\Bbb M}_{1}}^{4}\Phi_{1}(z).
\end{equation}

\subsubsection{The Borcherds $\Phi$-function with respect to the level $1$ cusp}
\label{subsubsec:Borcherds:Phi:level:1:cusp}
\par
By \cite[Example 13.7]{Borcherds98}, \cite[(8.5), Example 8.9]{Yoshikawa13},
$\Phi_{1}(z)$ is expressed as the following infinite product 
\begin{equation}
\label{eqn:Phi:expansion:level1}
\Phi_{1}(z)
=
\prod_{\lambda\in{\Bbb M}_{1}\cap\overline{\mathcal C}_{{\Bbb M}_{1}}^{+}\setminus\{0\}}
\left(
\frac{1-e^{\pi i\langle\lambda,z\rangle_{{\Bbb M}_{1}}}}
{1+e^{\pi i\langle\lambda,z\rangle_{{\Bbb M}_{1}}}}
\right)^{c(\lambda^{2}/2)}
\end{equation}
when $(\Im z)^{2}\gg0$. Here the series $\{c(n)\}\subset{\bf Z}$ is defined by the generating function:
$$
\sum_{n\in{\bf Z}}c(n)\,e^{2\pi in\tau}=\eta(\tau)^{-8}\eta(2\tau)^{8}\eta(4\tau)^{-8},
\qquad
\eta(\tau):=e^{2\pi i\tau/24}\prod_{n=1}^{\infty}(1-e^{2\pi in\tau}).
$$
By this explicit expression, $\Phi(Z)$ is defined without an ambiguity of constant now.

\subsubsection{The Borcherds $\Phi$-function with respect to the level $2$ cusp}
\label{subsubsec:Borcherds:Phi:level:2:cusp}
\par
Let ${\mathcal W}\subset{\mathcal C}_{{\Bbb M}_{2}}^{+}$ be a Weyl chamber with ${\frak e}_{1}\in\overline{\mathcal W}$
and set $\Pi^{+}:=\{\lambda\in{\Bbb M}_{2};\,\langle\lambda,{\mathcal W}\rangle_{{\Bbb M}_{2}}>0,\,\lambda^{2}\geq-2\}$.

\begin{theorem}
\label{theorem:Phi:expansion:level2}
If $(\Im w)^{2}\gg0$, $\Phi_{2}(w)$ is expressed as the following infinite product:
\begin{equation}
\label{eqn:Phi:expansion:level2}
\Phi_{2}(w) 
=
2^{8}e^{2\pi i\langle{\frak e}_{1}, w\rangle_{{\Bbb M}_{2}}}
\prod_{\lambda \in\Pi^{+}} 
\left(
1 - e^{2\pi i \langle \l, w\rangle_{{\Bbb M}_{2}}}
\right)^{(-1)^{\langle \lambda,{\frak e}_{1}-{\frak f}_{1}\rangle_{{\Bbb M}_{2}}}c(\lambda^{2}/2)}.
\end{equation}
\end{theorem}

\begin{pf}
By \cite{Borcherds96}, \cite[Th.\,13.3 (5)]{Borcherds98}, \cite[(8.5), Example 8.9]{Yoshikawa13}, there is a constant $C$ with $|C|=2^{8}$ such that
\begin{equation}
\label{eqn:Phi:expansion:level2:constant}
\Phi_{2}(w) 
=
C\,e^{2\pi i\langle{\frak e}_{1}, w\rangle_{{\Bbb M}_{2}}}
\prod_{\lambda \in\Pi^{+}} 
\left(
1 - e^{2\pi i \langle \l, w\rangle_{{\Bbb M}_{2}}}
\right)^{(-1)^{\langle \lambda,{\frak e}_{1}-{\frak f}_{1}\rangle_{{\Bbb M}_{2}}}c(\lambda^{2}/2)}.
\end{equation}
\par
For $z\in{\Bbb M}_{1}\otimes{\bf R}+i\,{\mathcal C}_{{\Bbb M}_{1}}^{+}$, write
$z=\langle z,{\frak f}_{2}\rangle_{{\Bbb M}_{1}}({\frak e}_{2}/2)+\langle z,{\frak e}_{2}\rangle_{{\Bbb M}_{1}}({\frak f}_{2}/2)+z_{{\Bbb E}_{8}(2)}$
with $z_{{\Bbb E}_{8}(2)}\in{\Bbb E}_{8}(2)\otimes{\bf C}$ and $\tau:=\langle{\frak e}_{2},z\rangle_{{\Bbb M}_{1}}$.
Consider the limit $\Im\langle{\frak f}_{2},z\rangle_{{\Bbb M}_{1}}\to+\infty$ with $\langle z,{\frak e}_{2}\rangle_{{\Bbb M}_{1}}$ 
and $z_{{\Bbb E}_{8}(2)}$ bounded. 
Since
$\Im\langle\lambda,z\rangle_{{\Bbb M}_{1}}\to+\infty$ for all $\lambda\in\overline{\mathcal C}_{{\Bbb M}_{1}}^{+}\setminus\{0\}$
with $\langle\lambda,{\frak e}_{2}\rangle_{{\Bbb M}_{1}}\not=0$, we get $\exp(2\pi i\langle\lambda,z\rangle_{{\Bbb M}_{1}})\to0$
for all $\lambda\in\overline{\mathcal C}_{{\Bbb M}_{1}}^{+}\setminus{\bf R}_{\geq0}{\frak e}_{2}$ as 
$\Im\langle{\frak f}_{2},z\rangle_{{\Bbb M}_{1}}\to+\infty$. 
We set $\Phi_{1}^{E}(\tau):=\lim_{\Im\langle z,{\frak f}_{2}\rangle_{{\Bbb M}_{1}}\to+\infty}\Phi_{1}(z)$.
Then we get by \eqref{eqn:Phi:expansion:level1}
\begin{equation}
\label{eqn:Siegel:operator:level:1:cusp}
\Phi_{1}^{E}(\tau)
=
\prod_{n\in{\bf Z}_{>0}}
\left(
\frac{1-e^{\pi i\langle n{\frak e}_{2},z\rangle_{{\Bbb M}_{1}}}}
{1+e^{\pi i\langle n{\frak e}_{2},z\rangle_{{\Bbb M}_{1}}}}
\right)^{c(\frac{(n{\frak e}_{2})^{2}}{2})}
=
\prod_{n>0}\left(\frac{1-e^{\pi in\tau}}{1+e^{\pi in\tau}}\right)^{8}
=
\frac{\eta(\tau/2)^{16}}{\eta(\tau)^{8}}.
\end{equation}
\par
We set $w_{{\Bbb E}_{8}(z)}:=w-\langle{\frak f}_{1},w\rangle_{{\Bbb M}_{2}}{\frak e}_{1}-\langle{\frak e}_{1},w\rangle_{{\Bbb M}_{2}}{\frak f}_{1}
\in{\Bbb E}_{8}(2)\otimes{\bf C}$.
Since $z^{2}=\langle z,{\frak e}_{2}\rangle_{{\Bbb M}_{1}}\langle z,{\frak f}_{2}\rangle_{{\Bbb M}_{1}}+z_{{\Bbb E}_{8}(2)}^{2}$,
we get $\Im\langle{\frak f}_{1},w\rangle_{{\Bbb M}_{2}}\to+\infty$ and $\langle{\frak e}_{1},w\rangle_{{\Bbb M}_{2}}=-1/\tau$,
$w_{{\Bbb E}_{8}(2)}=-z_{{\Bbb E}_{8}(2)}/\tau$ by \eqref{eqn:transformation}. 
We set $\sigma:=\langle{\frak e}_{1},w\rangle_{{\Bbb M}_{2}}$ and
$\Phi_{2}^{E}(\sigma):=\lim_{\Im\langle w,{\frak f}_{1}\rangle_{{\Bbb M}_{2}}\to+\infty}\Phi_{2}(w)$.
In the same way as above, we deduce from \eqref{eqn:Phi:expansion:level2:constant} that
\begin{equation}
\label{eqn:Siegel:operator:level:2:cusp}
\Phi_{2}^{E}(\sigma)
=
C\,e^{2\pi i\sigma}\prod_{n\in{\bf Z}_{>0}}
\left(
1-e^{2\pi i\langle n{\frak e}_{1},w\rangle_{{\Bbb M}_{2}}}
\right)^{(-1)^{n}c(\frac{(n{\frak e}_{1})^{2}}{2})}
=
C\,\eta(2\sigma)^{16}/\eta(\sigma)^{8}.
\end{equation}
Since $\sigma=-1/\tau$, we get the following by substituting \eqref{eqn:Siegel:operator:level:1:cusp} and
\eqref{eqn:Siegel:operator:level:2:cusp} into \eqref{eqn:relation:Phi1:Phi2:sec2.2}
\begin{equation}
\label{eqn:relation:Siegel:operators}
C\,\eta(-2/\tau)^{16}/\eta(-1/\tau)^{8}=\Phi_{2}^{E}(-1/\tau)=\tau^{4}\Phi_{1}^{E}(\tau)=\tau^{4}\cdot\eta(\tau/2)^{16}/\eta(\tau)^{8}.
\end{equation}
Since $\eta(-1/\tau)^{8}=\tau^{4}\eta(\tau)^{8}$, we get $C=2^{8}$ by \eqref{eqn:relation:Siegel:operators}.
\end{pf}

\subsection{Degenerations of Enriques surfaces and the Borcherds $\Phi$-function}
\label{subsec:degeneration}

\begin{theorem}
\label{thm:degeneration:Phi}
Let $\varDelta\subset{\bf C}$ be the unit disc.
Let $Z_{0}$ be a $K3$ surface with at most nodes as its singularities and 
let $f\colon(Z,Z_{0})\to(\varDelta,0)$ be a flat deformation of $Z_{0}$.
Let $\iota\colon Z\to Z$ be a holomorphic involution preserving the fibers of $f$.
For $t\in\varDelta$, set $Z_{t}:=f^{-1}(t)$, $\iota_{t}:=\iota|_{Z_{t}}$ and $S_{t}:=Z_{t}/\iota_{t}$. 
Assume the following:
\begin{itemize}
\item[(i)]
$Z$ is smooth.
\item[(ii)] 
The fixed-point-set of the $\iota$-action on $Z$ consists of one node of $Z_{0}$, say ${\frak o}$.  
\item[(iii)]
The map $f\colon Z\to\varDelta$ is projective.
\end{itemize}
Then $S_{t}$ is an Enriques surface for $t\not=0$ and the following holds
$$
\log\|\Phi(S_{t})\|^{2}
=
\frac{1}{2}\,\log|t|^{2}+O(1)
\qquad
(t\to0).
$$
\end{theorem}

\begin{pf}
{\em (Step 1) }
Let $\widetilde{\varDelta}$ be another disc and set $\widetilde{\varDelta}^{*}:=\widetilde{\varDelta}\setminus\{0\}$.
Let $Z\times_{\varDelta}\widetilde{\varDelta}$ be the family over $\widetilde{\varDelta}$ 
induced from $f\colon Z\to\varDelta$ by the map $\widetilde{\varDelta}\ni t\to t^{2}\in\varDelta$.
A resolution $\pi\colon\widetilde{Z}\to Z\times_{\varDelta}\widetilde{\varDelta}$ of the singularities of 
$Z\times_{\varDelta}\widetilde{\varDelta}$ is called a {\em simultaneous resolution} of $f\colon Z\to\varDelta$ if the following property is satisfied:
Set $\widetilde{\pi}:={\rm pr}_{1}\circ\pi\colon\widetilde{Z}\to Z$
and $\widetilde{f}:={\rm pr}_{2}\circ\pi\colon\widetilde{Z}\to\widetilde{\varDelta}$. 
For $t\in\widetilde{\varDelta}$, we set $\widetilde{Z}_{t}:=\widetilde{f}^{-1}(t)$ and 
$\widetilde{\pi}_{t}:=\widetilde{\pi}|_{\widetilde{Z}_{t}}\colon\widetilde{Z}_{t}\to Z_{t^{2}}$.
Then $\widetilde{\pi}_{t}$ is an isomorphism for $t\in\widetilde{\varDelta}^{*}$ and is the minimal resolution for $t=0$. 
\par
There exists a simultaneous resolution $\widetilde{f}:\widetilde{Z}\to\widetilde{\varDelta}$ of $f\colon Z\to\varDelta$.
Then $\widetilde{f}$ is a smooth (possibly non-projective) morphism and $Z_{t}$ is a smooth $K3$ surface for $t\not=0$.
By (ii), $S_{t}$ is an Enriques surface for $t\not=0$.
\par{\em (Step 2) }
We recall the simultaneous resolution of $f\colon Z\to\varDelta$.
By (ii), we can write ${\rm Sing}(Z_{0})=\{{\frak o},{\frak p}_{1},{\frak q}_{1},\ldots,{\frak p}_{m},{\frak q}_{m}\}$
with $\iota_{0}({\frak p}_{i})={\frak q}_{i}$.
By (iii) and an argument using morsification, there is a system of coordinates $(O,(z_{1},z_{2},z_{3}) )$ near ${\frak o}$ such that
\begin{equation}
\label{eqn:local:model:projection:involution}
f(z)=z_{1}^{2}+z_{2}^{2}+z_{3}^{2},
\qquad
\iota(z)=-z.
\end{equation}
Similarly, there is a system of coordinates $(U_{\alpha},(u_{1},u_{2},u_{3}))$ (resp. $(V_{\alpha},(v_{1},v_{2},v_{3}))$)
around ${\frak p}_{\alpha}$ (resp. ${\frak q}_{\alpha}$) such that
\begin{equation}
\label{eqn:local:model:projection:involution:2}
f(u)=u_{1}^{2}+u_{2}^{2}+u_{3}^{2},
\qquad
f(v)=v_{1}^{2}+v_{2}^{2}+v_{3}^{2},
\qquad
\iota(U_{\alpha})=V_{\alpha},
\qquad
\iota^{*}v_{\alpha}^{i}=u_{\alpha}^{i}.
\end{equation}
\par
Define
$O':=O\times_{\varDelta}\widetilde{\varDelta}=
\{(z,t)\in O\times\widetilde{\varDelta};\,z_{1}^{2}+z_{2}^{2}+z_{3}^{2}-t^{2}=0\}$
and
$$
\widetilde{O}:=\{((z,t),(\zeta_{0}:\zeta_{1}))\in O'\times{\bf P}^{1};\,(\zeta_{0}:\zeta_{1})=(z_{1}+\sqrt{-1}z_{2}:z_{3}+t)\}.
$$
We define $\widetilde{U}_{\alpha}$ and $\widetilde{V}_{\alpha}$ in the same manner.
By \cite[Sect.\,2.7]{Brieskorn66}, $\widetilde{O}$ is a complex manifold and is the closure of the graph of the rational map 
$$
O\times_{\varDelta}\widetilde{\varDelta}\ni(z,t)\dashrightarrow(z_{1}+\sqrt{-1}z_{2}:z_{3}+t)=(t-z_{3}:z_{1}-\sqrt{-1}z_{2})\in{\bf P}^{1}.
$$
Moreover, the obvious projection $\pi:={\rm pr}_{1}\colon\widetilde{O}\to O\times_{\varDelta}\widetilde{\varDelta}$ 
gives a resolution of the singularity of $O\times_{\varDelta}\widetilde{\varDelta}$ with $\pi^{-1}({\frak o})\cong{\bf P}^{1}$. 
Similarly, $\pi\colon\widetilde{U}_{\alpha}\to U_{\alpha}\times_{\varDelta}\widetilde{\varDelta}$ and
$\pi\colon\widetilde{V}_{\alpha}\to V_{\alpha}\times_{\varDelta}\widetilde{\varDelta}$ are resolutions with
$\pi^{-1}({\frak p}_{\alpha})\cong{\bf P}^{1}$ and $\pi^{-1}({\frak q}_{\alpha})\cong{\bf P}^{1}$.
Then $\widetilde{Z}$ is obtained from $Z\times_{\varDelta}\widetilde{\varDelta}$
by replacing $O\times_{\varDelta}\widetilde{\varDelta}$, $U_{\alpha}\times_{\varDelta}\widetilde{\varDelta}$, 
$V_{\alpha}\times_{\varDelta}\widetilde{\varDelta}$
by $\widetilde{O}$, $\widetilde{U}_{\alpha}$, $\widetilde{V}_{\alpha}$, respectively.
We set 
$$
E_{0}:=\widetilde{\pi}_{0}^{-1}({\frak o}),
\qquad
F_{\alpha}:=\widetilde{\pi}_{0}^{-1}({\frak p}_{\alpha}),
\qquad
F'_{\alpha}:=\widetilde{\pi}_{0}^{-1}({\frak q}_{\alpha}).
$$
Then $E_{0},F_{1},F'_{1},\ldots,F_{m},F'_{m}$ are mutually disjoint $(-2)$-curves of $\widetilde{Z}_{0}$.
\par{\em (Step 3) }
Let $\widetilde{\iota}$ be the {\em meromorphic} involution on $\widetilde{Z}$ induced by the involution $\iota\times{\rm id}_{\varDelta}$
on $Z\times_{\varDelta}\widetilde{\varDelta}$. 
Then $\widetilde{\iota}$ is a holomorphic involution on $\widetilde{Z}\setminus E_{0}$ exchanging $F_{\alpha}$ and $F'_{\alpha}$.
For $t\in\widetilde{\varDelta}^{*}$, set $\widetilde{\iota}_{t}:=\widetilde{\iota}|_{\widetilde{Z}_{t}}$ and
$\widetilde{S}_{t}:=\widetilde{Z}_{t}/\widetilde{\iota}_{t}$. 
Since the family $\widetilde{f}\colon\widetilde{Z}\to\widetilde{\varDelta}$ is differentiably trivial, 
it admits a marking $\mu$, i.e., a trivialization of the local system
$\mu\colon R^{2}\widetilde{f}_{*}{\bf Z}\cong{\Bbb L}_{K3}$
such that the condition \eqref{eqn:alpha:+:-} is satisfied for all $t\in\widetilde{\varDelta}^{*}$. 
Let 
$$
\varpi\colon\widetilde{\varDelta}\ni t\to\varpi(\widetilde{S}_{t},\mu)\in\Omega_{\LAM}^{+}
$$ 
be the period map for the marked family $(\widetilde{f}\colon(\widetilde{Z},\widetilde{\iota})\to\widetilde{\varDelta},\mu)$. 
Let $\varPi\colon\Omega_{\LAM}^{+}\to{\mathcal M}$ be the projection. 
Since $\widetilde{Z}_{t}=Z_{t^{2}}$, $\widetilde{\iota}_{t}=\iota_{t^{2}}$ and hence $\widetilde{S}_{t}=S_{t^{2}}$, we have 
\begin{equation}
\label{eqn:relation:period:Griffiths:period}
\varPi\circ\varpi(t)=\overline{\varpi}(S_{t^{2}}).
\end{equation}
\par{\em (Step 4) }
Let $\varSigma:=\{(x_{1},x_{2},x_{3})\in{\bf R}^{3};\,x_{1}^{2}+x_{2}^{2}+x_{3}^{2}=1\}$ be the unit sphere of ${\bf R}^{3}$.
Define the embedding 
$i\colon\varSigma\times\widetilde{\varDelta}\hookrightarrow\widetilde{O}\subset\widetilde{Z}$ by
\begin{equation}
\label{eqn:vanishing:cycle}
i((x_{1},x_{2},x_{3}),t):=((tx_{1},tx_{2},tx_{3},t),(x_{1}+\sqrt{-1}x_{2}:x_{3}+1)).
\end{equation}
We define the embeddings $j_{\alpha}\colon\varSigma\times\varDelta\hookrightarrow\widetilde{U}_{\alpha}$ and
$j'_{\alpha}\colon\varSigma\times\varDelta\hookrightarrow\widetilde{V}_{\alpha}$ in the same way.
For $t\in\widetilde{\varDelta}$, define submanifolds $E_{t},F_{\alpha,t},F'_{\alpha,t}\subset\widetilde{Z}_{t}$ diffeomorphic to $\varSigma$ as
$$
E_{t}:=i(\varSigma\times\{t\}),
\qquad
F_{\alpha,t}:=j_{\alpha}(\varSigma\times\{t\}),
\qquad
F'_{\alpha,t}:=j'_{\alpha}(\varSigma\times\{t\}).
$$
Then the $2$-spheres $E_{t}$, $F_{\alpha,t}$, $F'_{\alpha,t}$ are transcendental cycles for $t\not=0$, while
$E_{0},F_{\alpha,0}=F_{\alpha},F'_{\alpha,0}=F'_{\alpha}\subset\widetilde{Z}_{0}$ are $(-2)$-curves contracted by $\widetilde{f}_{0}$.
By \eqref{eqn:local:model:projection:involution:2} and the definitions of $F_{\alpha,t}$, $F'_{\alpha,t}$, we get for all $t\in\widetilde{\varDelta}$
\begin{equation}
\label{eqn:exceptional:curves:involution}
\widetilde{\iota}_{t}(F_{\alpha,t})=F'_{\alpha,t},
\qquad
\widetilde{\iota}_{t}(F'_{\alpha,t})=F_{\alpha,t}.
\end{equation}
\par
Let $c_{1}(E_{t})\in H^{2}(\widetilde{Z}_{t},{\bf Z})$ be the Poincar\'e dual of $E_{t}$ and define
$$
\delta:=\mu(c_{1}(E_{0}))=\mu(c_{1}(E_{t}))\in\Delta_{{\Bbb L}_{K3}}.
$$
Similarly, we define roots $d_{\alpha},d'_{\alpha}\in\Delta_{{\Bbb L}_{K3}}$ as
$$
d_{\alpha}:=\mu(c_{1}(F_{\alpha,0}))=\mu(c_{1}(F_{\alpha,t})),
\qquad
d'_{\alpha}:=\mu(c_{1}(F'_{\alpha,0}))=\mu(c_{1}(F'_{\alpha,t})).
$$
Since $\widetilde{\iota}_{t}$ preserves $E_{t}$ and reverses its orientation for $t\not=0$ 
by \eqref{eqn:local:model:projection:involution}, \eqref{eqn:vanishing:cycle}, we get $\delta\in\Delta_{\LAM}$. 
Since $d_{\alpha}$ and $d'_{\alpha}$ are not eigenvectors of the involution 
$\mu\circ\widetilde{\iota}_{t}^{*}\circ\mu^{-1}\in O({\Bbb L}_{K3})$ by \eqref{eqn:exceptional:curves:involution},
we get $d_{\alpha},d'_{\alpha}\not\in\LAM$.
\par{\em (Step 5) }
Since $E_{0}$ is an algebraic cycle of $\widetilde{Z}_{0}$ and $\delta\in\Delta_{\LAM}$, we get $\varpi(0)\in H_{\delta}$.
Set $H_{\delta}^{o}:=H_{\delta}\setminus\bigcup_{d\in\Delta(\LAM)\setminus\{\pm\delta\}}H_{d}$.
Let us see that $\varpi(0)\in H_{\delta}^{o}$.
By (iii), there is an $\iota$-invariant relatively ample line bundle $L$ on $Z$. 
Set $\widetilde{L}:=\widetilde{\pi}^{*}L$ and $\widetilde{L}_{t}:=\widetilde{L}|_{\widetilde{Z}_{t}}$ for $t\in\widetilde{\varDelta}$. 
We get $l:=\mu(c_{1}(\widetilde{L}_{0}))=\mu(c_{1}(\widetilde{L}_{t}))\in\LAM^{\perp}$ by the $\widetilde{\iota}_{t}$-invariance of $\widetilde{L}_{t}$. 
\par
Let $d\in\Delta_{\LAM}$ be such that $\varpi(0)\in H_{d}$. 
By \cite[Chap.\,VIII, Prop.\,3.7 (i)]{BPV84}, either $d$ or $-d$ is effective.
For simplicity, assume that $d$ is effective. There is an effective divisor $\Gamma$ on $\widetilde{Z}_{0}$ such that
$\mu(c_{1}(\Gamma))=d\in\Delta_{\LAM}$. 
If  $\widetilde{\pi}_{0}(\Gamma)$ contains a curve, then we get
$0<\deg(L_{0}|_{\widetilde{\pi}_{0}(\Gamma)})=\deg(\widetilde{L}_{0}|_{\Gamma})=\langle\mu(c_{1}(\widetilde{L}_{0})),\mu(c_{1}(\Gamma))\rangle
=\langle l,d\rangle=0$, a contradiction. 
Hence $\dim\widetilde{\pi}_{0}(\Gamma)=0$ and ${\rm Supp}(\Gamma)$ is contained in the exceptional divisor of $\widetilde{\pi}_{0}$.
\par
Write $\Gamma=\nu E_{0}+\sum_{\alpha}\nu_{\alpha}F_{\alpha}+\sum_{\alpha}\nu'_{\alpha}F'_{\alpha}$, 
where $\nu,\nu_{\alpha},\nu'_{\alpha}\in{\bf Z}$.
Since $\mu(c_{1}(\Gamma))=d\in\Delta_{\LAM}$ and since $E_{0},F_{1},F'_{1},\ldots,F_{m},F'_{m}$ are mutually disjoint $(-2)$-curves,
we get $\nu^{2}+\sum_{\alpha}\nu_{\alpha}^{2}+\sum_{\alpha}(\nu'_{\alpha})^{2}=1$.
Hence $\Gamma=E_{0},F_{\alpha},F'_{\alpha}$. 
Since $d_{\alpha},d'_{\alpha}\not\in\LAM$, $\delta\in\LAM$ by Step 4 and since $d\in\LAM$ by assumption, 
we get $\Gamma\not=F_{\alpha},F'_{\alpha}$ and $\Gamma=E_{0}$.
This proves that, if $d\in\Delta_{\LAM}$ and $\varpi(0)\in H_{d}$, then $d=\pm\delta$.
Namely, $\varpi(0)\in H_{\delta}^{o}$.
\par{\em (Step 6) }
Let $K_{Z}$ be the canonical bundle of $Z$. Since $K_{Z}$ is a trivial line bundle on $Z$,
there is a nowhere vanishing holomorphic $3$-form $\xi$ on $Z$. 
For $t\in\varDelta$, we set 
$$
\eta_{t}:={\rm Res}_{Z_{t}}[\xi/(f(z)-t)]\in H^{0}(Z_{t},K_{Z_{t}})\setminus\{0\}.
$$
Then $\widetilde{\eta}_{t}:=\eta_{t^{2}}$ is regarded as a non-zero canonical form on $\widetilde{Z}_{t}$ for $t\not=0$. 
By \eqref{eqn:local:model:projection:involution}, we can express $\eta_{t}=e^{\psi(z)}dz_{2}\wedge dz_{3}/z_{1}$
on a neighborhood of ${\frak o}$, where $\psi(z)$ is a holomorphic function near ${\frak o}$.
By this expression of $\eta_{t}$ and \eqref{eqn:vanishing:cycle}, we get
$$
\langle\alpha(\widetilde{\eta}_{t}),\delta\rangle
=
\int_{E_{t}}\widetilde{\eta}_{t}
=
\int_{t\varSigma}e^{\psi(z)}\frac{dz_{2}\wedge dz_{3}}{z_{1}}
=
t\,\left\{e^{\psi(0)}\int_{\varSigma}\frac{dz_{2}\wedge dz_{3}}{z_{1}}\right\}+O(t^{2}),
$$
where $t\varSigma:=\{(tx_{1},tx_{2},tx_{3})\in{\bf C}^{3};\,(x_{1},x_{2},x_{3})\in\varSigma\}$.
This proves that $\varpi(t)$ intersects $H_{\delta}^{o}$ transversally at $\varpi(0)$.
Since $\Phi(z)$ vanishes of order one on $H_{\delta}^{o}$,
$$
\log\|\Phi(\varpi(\widetilde{S}_{t},\alpha))\|^{2}=\log\|\Phi(\varpi(t))\|^{2}=\log|t|^{2}+O(1)
\qquad
(t\to0).
$$ 
This, together with \eqref{eqn:relation:period:Griffiths:period}, yields the result.
\end{pf}

\section{An algebraic expression of $\|\Phi\|$}
\label{sec:expression:Phi:generic}

\subsection
{The $(2,2,2)$-model of an Enriques surface}
\label{sect:3.1}
\par
Let ${\rm Sym}(3,{\bf C})$ be the set of complex $3\times3$-symmetric matrices.
For $S=(s_{ij})\in{\rm Sym}(3,{\bf C})$, we set
$$
Q(x;S):=\sum_{i,j=1}^{3}s_{ij}x_{i}x_{j}\in{\bf C}[x_{1},x_{2},x_{3}].
$$
For vectors $A=(A_{1},A_{2},A_{3})$ and $B=(B_{1},B_{2},B_{3})$ in ${\rm Sym}(3,{\bf C})\otimes{\bf C}^{3}$, we define
$$
X_{(A,B)}
:=
\{(x_{1}:x_{2}:x_{3}:y_{1}:y_{2}:y_{3})\in{\bf P}^{5};\,
Q(x;A_{i})+Q(y;B_{i})=0
\quad(i=1,2,3)
\}.
$$
If $X_{(A,B)}$ is a smooth complex surface, then $X_{(A,B)}$ is a $K3$ surface by the adjunction formula. 
We say that $(A,B)$ is {\em admissible} if $R(A)R(B)\not=0$ and $X_{(A,B)}$ is a reduced complex surface with at most  RDP's. 
When $(A,B)$ is admissible, $X_{(A,B)}$ is a $K3$ surface with at most RDP, i.e., the minimal resolution of $X_{(A,B)}$ is a $K3$ surface.
\par
Let $R(A)=R(A_{1},A_{2},A_{3})$ be the {\em resultant} of $Q(x;A_{1})$, $Q(x;A_{2})$, $Q(x;A_{3})$.
Namely, $R(A)$ is the unique {\em irreducible} polynomial in the entries of $A_{1},A_{2},A_{3}$ 
satisfying the following conditions (cf. \cite[Chap.\,3, (2.3)]{CLO05}):
\begin{itemize}
\item[(i)]
The system of equations $Q(x;A_{1})=Q(x;A_{2})=Q(x;A_{3})=0$ has a non-trivial solution if and only if $R(A)=0$.
\item[(ii)]
If $A_{1}={\rm diag}(1,0,0)$, $A_{2}={\rm diag}(0,1,0)$, $A_{3}={\rm diag}(0,0,1)$, then $R(A)=1$.
\end{itemize}
By \cite[Chap.\,3]{CLO05}, $R(A)$ is a homogeneous polynomial of degree $4$ in the entries of $A_{i}$, $1\leq i\leq3$, 
so that $\deg R(A)=12$. See \cite[p.215, Table 1]{KSZ92} for an explicit formula.
\par
By definition, $X_{(A,B)}$ is preserved by the involution on ${\bf P}^{5}$
\begin{equation}
\label{eqn:involution:P5}
\iota\colon(x_{1}:x_{2}:x_{3}:y_{1}:y_{2}:y_{3})\to(x_{1}:x_{2}:x_{3}:-y_{1}:-y_{2}:-y_{3}).
\end{equation}
When $\iota$ preserves a subset $V\subset{\bf P}^{5}$, we define $V^{\iota}:=\{p\in V;\,\iota(p)=p\}$.
We have $({\bf P}^{5})^{\iota}=P_{1}\amalg P_{2}$, where $P_{1}=\{x_{1}=x_{2}=x_{3}=0\}$ and $P_{2}=\{y_{1}=y_{2}=y_{3}=0\}$ 
are projective planes.
Since the three conics $Q(x;A_{i})=0$ $(i=1,2,3)$ have no points in common if and only if $R(A)\not=0$, 
$X_{(A,B)}^{\iota}=\emptyset$ if and only if $R(A)R(B)\not=0$. Hence
$$
Y_{(A,B)}:=X_{(A,B)}/\iota
$$ 
is an Enriques surface with at most RDP's (i.e., the minimal resolution of $Y_{(A,B)}$ is an Enriques surface) for admissible $(A,B)$.
Let ${\mathcal L}_{(A,B)}$ be the ample line bundle of degree $4$ on $Y_{(A,B)}$ induced from the $\iota$-invariant ample line bundle
${\mathcal O}_{{\bf P}^{5}}(1)$ on $X_{(A,B)}$.
When $(A,B)$ is admissible, the polarized variety $(Y_{(A,B)},{\mathcal L}_{(A,B)})$ or equivalently 
$(X_{(A,B)},\iota,{\mathcal O}_{{\bf P}^{5}}(1))$ is called a {\em $(2,2,2)$-model} in this paper. 
For simplicity, we often omit ${\mathcal L}_{(A,B)}$, ${\mathcal O}_{{\bf P}^{5}}(1)$.
Since every Enriques surface is birational to some $Y_{(A,B)}$ with admissible $(A,B)$ by \cite{Verra83}, \cite{Cossec85},
every point of ${\mathcal M}^{o}$ admits a $(2,2,2)$-model.
\par
For an admissible $(A,B)$, let $K_{X_{(A,B)}}$ be the dualizing sheaf on $X_{(A,B)}$. 
We define $\omega_{(A,B)}\in H^{0}(X_{(A,B)},K_{X_{(A,B)}})$ as the residue of $Q(x;A_{i})+Q(y;B_{i})$ $(i=1,2,3)$.
Let $q_{i}$ $(i=1,2,3)$ be the inhomogeneous equation of $Q(x;A_{i})+Q(y;B_{i})$ obtained by setting $y_{3}=1$.
Then $dx_{1}\wedge\cdots\wedge dy_{2}/q_{1}q_{2}q_{3}$ is a canonical form on ${\bf C}^{5}\setminus{\rm Sing}\,X_{(A,B)}$
with logarithmic poles along ${\rm div}\,q_{1}\cup{\rm div}\,q_{2}\cup{\rm div}\,q_{3}$, which extends to canonical form on
${\bf P}^{5}\setminus{\rm Sing}\,X_{(A,B)}$ with logarithmic poles along $\cup_{i=1}^{3}{\rm div}(Q(x;A_{i})+Q(y;B_{i}))$.
Then there exists a locally defined $2$-form $\Upsilon$ on ${\bf P}^{5}\setminus{\rm Sing}\,X_{(A,B)}$ such that
\begin{equation}
\label{eqn:formula:varXi}
dx_{1}\wedge\cdots\wedge dy_{2}/q_{1}q_{2}q_{3}
=
(dq_{1}/q_{1})\wedge(dq_{2}/q_{2})\wedge(dq_{3}/q_{3})\wedge\Upsilon
\end{equation}
Let $j\colon X_{(A,B)}\setminus{\rm Sing}\,X_{(A,B)}\to X_{(A,B)}$ be the inclusion. We define
\begin{equation}
\label{eqn:2-form:residue}
\omega_{(A,B)}=j_{*}(\Upsilon|_{X_{(A,B)}\setminus{\rm Sing}\,X_{(A,B)}}),
\end{equation}
which is independent of $\Upsilon$.
In this section, we prove the following theorem.

\begin{theorem}
\label{thm:valuePhi:general}
For every admissible $(A,B)\in{\rm Sym}(3,{\bf C})\otimes{\bf C}^{6}$,
$$
\|\Phi(Y_{(A,B)})\|^{2}
=
|R(A)R(B)|\cdot
\left(
\frac{2}{\pi^{4}}
\left|\int_{X_{(A,B)}}\omega_{(A,B)}\wedge\overline{\omega}_{(A,B)}\right|
\right)^{4}.
$$
\end{theorem}

\subsection
{A Grassmannian and a family of $K3$ surfaces with involution}
\label{subsec:Grassmannian}
\par
Let $V$ be the complex vector space of quadratic forms in the variables $x_{1},x_{2},x_{3},y_{1},y_{2},y_{3}$ invariant under
the involution \eqref{eqn:involution:P5}
$$
\begin{aligned}
V
&:=
{\bf C}x_{1}^{2}+{\bf C}x_{2}^{2}+{\bf C}x_{3}^{2}+{\bf C}x_{1}x_{2}+{\bf C}x_{1}x_{3}+{\bf C}x_{2}x_{3}
\\
&\qquad
+{\bf C}y_{1}^{2}+{\bf C}y_{2}^{2}+{\bf C}y_{3}^{2}+{\bf C}y_{1}y_{2}+{\bf C}y_{1}y_{3}+{\bf C}y_{2}y_{3}.
\end{aligned}
$$
Let $\Psi\colon{\bf P}^{5}\to{\bf P}^{11}$ be the morphism defined by
$$
\Psi(x:y):=(x_{1}^{2}:x_{2}^{2}:x_{3}^{2}:x_{1}x_{2}:x_{1}x_{3}:x_{2}x_{3}:y_{1}^{2}:y_{2}^{2}:y_{3}^{2}:y_{1}y_{2}:y_{1}y_{3}:y_{2}y_{3}).
$$
Then $\Psi$ induces the embedding ${\bf P}^{5}/\iota\hookrightarrow{\bf P}^{11}$. Hence ${\rm Sing}(\Psi({\bf P}^{5}))=\Psi(({\bf P}^{5})^{\iota})$.
\par
Let $S$ be the Grassmann variety of $3$-dimensional linear subspaces of $V$:
$$
S:={\rm Gr}(3,V).
$$
We identify $S$ with the Grassmann variety of $9$-dimensional subspaces of $V^{\lor}={\bf C}^{12}$.
Let $L\subset{\bf C}^{12}\times{\rm Gr}(9,{\bf C}^{12})$ be  the tautological vector bundle of rank $9$ over ${\rm Gr}(9,{\bf C}^{12})$
and let
$$
\pi\colon{\bf P}(L)\to{\rm Gr}(9,{\bf C}^{12})
$$ 
be the projective-space bundle associated to $L$. We regard ${\bf P}(L)$ as a complex submanifold of codimension $3$
of ${\bf P}^{11}\times{\rm Gr}(9,{\bf C}^{12})$. Then $\pi={\rm pr}_{2}|_{{\bf P}(L)}$.
By the canonical identification between $S$ and ${\rm Gr}(9,{\bf C}^{12})$, 
the fiber $L_{s}:=p^{-1}(s)\subset{\bf P}^{11}$ is the linear subspace of codimension $3$ corresponding to $s\in S$.
\par
In ${\bf P}^{11}\times S$, we have two subvarieties $\Psi({\bf P}^{5})\times S$ and ${\bf P}(L)$. 
We define ${\mathcal Y}\subset{\bf P}^{11}\times S$ and ${\mathcal X}\subset{\bf P}^{5}\times S$ by
$$
\begin{array}{ll}
{\mathcal Y}
&:=
(\Psi({\bf P}^{5})\times S)\cap{\bf P}(L)
\\
{\mathcal X}
&:=
(\Psi\times{\rm id}_{S})^{-1}({\mathcal Y}),
\end{array}
$$
which are equipped with the projections $\pi\colon{\mathcal Y}\to S$ and $\pi\colon{\mathcal X}\to S$.
For $s\in S$, set
$$
Y_{s}:=\pi^{-1}(s)=\Psi({\bf P}^{5})\cap{\bf P}(L_{s}),
\qquad
X_{s}:=\pi^{-1}(s)=\Psi^{-1}(Y_{s}).
$$
If $\{F_{s},G_{s},H_{s}\}\subset V$ is a basis of the $3$-dimensional subspace of $V$ corresponding to $s\in S$,
then we have the expressions
$$
X_{s}=\{(x,y)\in{\bf P}^{5};\,F_{s}(x,y)=G_{s}(x,y)=H_{s}(x,y)=0\},
\qquad
Y_{s}=X_{s}/\iota.
$$
Here $X_{s}$ is regarded as the scheme defined by the ideal sheaf generated by $F_{s},G_{s},H_{s}$.
It is clear that for every $s\in S$, there exists $(A,B)\in{\rm Sym}(3,{\bf C})\otimes{\bf C}^{6}$ such that $X_{s}=X_{(A,B)}$.
Conversely, $X_{(A,B)}=X_{s}$ for some $s\in S$ if and only if the three vectors 
$(A_{1},B_{1}),(A_{2},B_{2}),(A_{3},B_{3})\in{\rm Sym}(3,{\bf C})\otimes{\bf C}^{2}$ are linearly independent. 
When $(A,B)\in{\rm Sym}(3,{\bf C})\otimes{\bf C}^{6}$ satisfies this condition, we set
$$
s(A,B):={\rm span}\{Q(x;A_{i})+Q(y;B_{i})\}_{i=1,2,3}\in S.
$$
Hence $X_{(A,B)}=X_{s(A,B)}$ in this case.
\par
We define
$$
\begin{array}{ll}
S^{*}
&:=
\{s\in S;\,\pi\colon{\mathcal X}\to S\hbox{ is flat at any }x\in X_{s}\},
\\
S^{o}
&:=
\{s\in S^{*};\,{\rm Sing}(X_{s})=\emptyset\}.
\end{array}
$$
Since the non-flat locus and the critical locus of $\pi$ are Zariski closed subsets of ${\mathcal X}$,
$S^{*}$ and $S^{o}$ are Zariski open subsets of $S$.
Then $\pi\colon{\mathcal X}|_{S^{*}}\to S^{*}$ is a flat family of projective surfaces, which is smooth over $S^{o}$.
We define $D^{*}:=S^{*}\setminus S^{o}$ and
$$
\begin{array}{ll}
D'
&:=
\{s\in D^{*};\,\dim{\rm Sing}(X_{s})=0\},
\\
D''
&:=
\{s\in D^{*};\,{\rm emb}(X_{s},x)=3\,(\forall\,x\in{\rm Sing}(X_{s}))\},
\\
D^{o}
&:=
\{s\in D^{*};\,{\rm Sing}(X_{s})\hbox{ consists of nodes}\},
\\
E
&:=
S\setminus(S^{o}\cup D^{o})=(S\setminus S^{*})\cup(D^{*}\setminus D^{o}),
\end{array}
$$
where ${\rm emb}(X_{s},x)$ denotes the embedding dimension of the germ $(X_{s},x)$.
Then $D'$ is a Zariski open subset of $D^{*}$. 
By Jacobi's criterion, $D''$ is a Zariski open subset of $D^{*}$. 
Hence $D'\cap D''$ is a Zariski open subset of $D^{*}$.

\begin{lemma}
\label{lemma:Grassmann:discriminant}
The following hold:
\begin{itemize}
\item[(1)]
$D^{o}$ is a Zariski open subset of $D^{*}$.
\item[(2)]
$E$ is a Zariski closed subset of $S$ with $\dim E\leq\dim S-2$.
\end{itemize}
\end{lemma}

\begin{pf}
{\bf (1) }
Let $s\in D'\cap D''$. Then ${\rm Sing}(X_{s})$ consists of isolated hypersurface singularities.
Let $(Z,0)\in{\rm Sing}(X_{s})$. Let ${\rm Def}(Z,0)$ be its miniversal deformation space
and let $f\colon({\mathcal Z},(Z,0))\to({\rm Def}(Z,0),0)$ be the miniversal deformation of $(Z,0)$.
There is a map $\mu_{Z}\colon(S,s)\to({\rm Def}(Z,0),0)$ such that $\pi\colon({\mathcal X},(Z,0))\to(S,s)$
is induced from $f\colon({\mathcal Z},(Z,0))\to({\rm Def}(Z,0),0)$ by $\mu_{Z}$.
Let $(D_{Z},0)$ be the discriminant locus of $f$. By \cite[Props.\,4.10 and 6.11, Cor.\,4.11]{Looijenga84}, 
there is a proper Zariski closed subset $(F_{Z},0)\subset(D_{Z},0)$ such that 
${\rm Sing}(f^{-1}(t))$ consists of a unique node for $t\in D_{Z}\setminus F_{Z}$.
Then $D^{o}=(D'\cap D'')\setminus\bigcup_{(Z,0)\in{\rm Sing}(X_{s})}\mu_{Z}^{-1}(F_{Z})$ on a small neighborhood $s$ in $S$.
By this expression, $D^{o}$ is a Zariski open subset of $D^{*}$.
\par{\bf (2) }
Let $s\in S^{o}\cup D^{o}$ be an arbitrary point. Then $X_{s}$ has at most nodes as its singular points.
Since any deformation of a node is either a smoothing or a trivial deformation,
$s$ is an interior point of $S^{o}\cup D^{o}$. Hence $E$ is a closed subset of $S$ with respect to the Euclidean topology. 
We get $E=\overline{E}=(S\setminus S^{*})\cup\overline{(D^{*}\setminus D^{o})}$, where $\overline{\mathcal A}$ denotes 
the closure of a subset ${\mathcal A}\subset S$ with respect to the Euclidean topology on $S$. 
Since $\overline{D^{*}\setminus D^{o}}$ is a Zariski closed subset of $S$ by (1), so is $E$.
\par
Assume that there is a component $E'\subset E$ with $\dim E'=\dim S-1$. Regard $E'$ as a reduced effective divisor of $S$.
Since $S$ is a Grassmann variety and hence ${\rm Pic}(S)={\bf Z}$, $E'$ must be an ample divisor of $S$.
Then, for any irreducible subvariety $T\subset S$ with $\dim T>0$, $T\cap E'$ is either a divisor of $T$ or $T$ itself.
In particular, we always have $\dim[T\cap E]\geq\dim T-1$ for any proper subvariety $T\subset S$ with $\dim T>0$.
However, by Lemma~\ref{lemma:family:K3:configuration} below, there exists a subGrassmannian $T\subset S$ of dimension $9$
such that $\dim[T\cap E]=\dim[T\setminus(S^{o}\cup D^{o})]\leq\dim T-2$. This contradicts the inequality 
$\dim[T\cap E]\geq\dim T-1$. This proves that $\dim E\leq\dim S-2$.
\end{pf}

\subsection
{A subfamily parametrized by a subGrassmannian}
\label{subsect:subGrassmannian}
\par
In this subsection, we use the coordinates $(x_{1},\ldots,x_{6})$ instead of $(x_{1},x_{2},x_{3},y_{1},y_{2},y_{3})$.
Write $M_{p,q}({\bf C})$ for the set of complex $p\times q$-matrices.
For $N=({\bf n}_{1},\ldots,{\bf n}_{6})\in M_{3,6}({\bf C})$, ${\bf n}_{i}\in{\bf C}^{3}$, we set 
$
\Delta_{ijk}(N):=\det({\bf n}_{i},{\bf n}_{j},{\bf n}_{k})
$ 
and we define
$$
M^{o}_{3,6}({\bf C})
:=
\{N\in M(3,6;{\bf C});\,\Delta_{ijk}(N)\not=0\quad(\forall\,1\leq i<j<k\leq6)\}.
$$
The Grassmann variety of $3$-dimensional subspaces of ${\bf C}^{6}$ is defined by
$$
T:={\rm Gr}(3,{\bf C}^{6})=GL({\bf C}^{3})\backslash\{N\in M(3,6;{\bf C});\,{\rm rk}(N)=3\}.
$$
For $N\in M(3,6;{\bf C})$ with ${\rm rk}(N)=3$, the corresponding point of $T$ is denoted by $[N]$.
For $1\leq i<j<k\leq6$, we set
$$
\begin{array}{ll}
T_{ijk}
&:=
GL({\bf C}^{3})\backslash\{N\in M(3,6;{\bf C});\,{\rm rk}(N)=3,\,\Delta_{ijk}(N)=0\},
\\
T^{o}
&:=
GL({\bf C}^{3})\backslash M^{o}_{3,6}({\bf C})=T\setminus\bigcup_{i<j<k}T_{ijk}.
\end{array}
$$
Then $T_{ijk}$ is a divisor of $T$ and $T^{o}$ is a non-empty Zariski open subset of $T$. We set
$$
\begin{array}{ll}
T_{ijk}^{o}
&:=
T_{ijk}\setminus\bigcup_{l<m<n,\,\{l,m,n\}\not=\{i,j,k\}}T_{lmn},
\\
D_{T}^{o}
&:=
\bigcup_{i<j<k}T^{o}_{ijk}.
\end{array}
$$
Then $T^{o}_{ijk}$ is a non-empty Zariski open subset of $T_{ijk}$ and hence
$T\setminus(T^{o}\cup D_{T}^{o})=\bigcup_{\{i,j,k\}\not=\{l,m,n\}}T_{ijk}\cap T_{lmn}$ is a Zariski closed subset of $T$.
\par
We regard $T$ as a subGrassmannian of $S$ by the embedding
$$
i\colon
T\ni[N]=[(n_{ij})]\mapsto
[{\bf C}\sum_{j=1}^{6}n_{1j}x_{j}^{2}+{\bf C}\sum_{j=1}^{6}n_{2j}x_{j}^{2}+{\bf C}\sum_{j=1}^{6}n_{3j}x_{j}^{2}]
\in S.
$$
Then we get for $N=({\bf n}_{1},\ldots,{\bf n}_{6})$ with ${\rm rk}(N)=3$
$$
X_{N}=X_{[N]}:=X_{i(N)}=\{[x]\in{\bf P}^{5};\,\sum_{k=1}^{6}x_{k}^{2}{\bf n}_{k}={\bf 0}\}.
$$

\begin{lemma}
\label{lemma:family:K3:configuration}
The following hold.
\begin{itemize}
\item[(1)]
For $[N]\in T^{o}$, $X_{[N]}$ is smooth. In particular, $T^{o}\subset S^{o}$.
\item[(2)]
For $[N]\in D_{T}^{o}$, ${\rm Sing}(X_{[N]})$ consists of nodes. In particular, $D_{T}^{o}\subset D^{o}$.
\item[(3)]
The Zariski closed subset $T\setminus(T^{o}\cup D_{T}^{o})$ has codimension at least $2$ in $T$.
\end{itemize}
\end{lemma}

\begin{pf}
Let $N=(n_{ij})=({\bf n}_{1},\ldots,{\bf n}_{6})\in M(3,6;{\bf C})$ and assume that $[N]\in T^{o}\cup D_{T}^{o}$.
Then there exists $i<j<k$ such that $\Delta_{lmn}(N)\not=0$ for any $\{l,m,n\}\not=\{i,j,k\}$.
In particular, any two column vectors ${\bf n}_{l}$, ${\bf n}_{m}$ of $N$ are linearly independent for $l\not=m$.
\par
Let $x=(x_{1}:\cdots:x_{6})\in X_{[N]}$. Then $\sum_{l}{\bf n}_{l}x_{l}^{2}={\bf 0}$.
If $4$ of the coordinates of $x$ vanish, then there are $l<m$ such that ${\bf n}_{l}x_{l}^{2}+{\bf n}_{m}x_{m}^{2}={\bf 0}$.
Since ${\bf n}_{l}$, ${\bf n}_{m}$ are linearly independent, we get $x_{l}=x_{m}=0$ and hence all coordinates of $x$ vanish. 
This contradicts $x\in{\bf P}^{5}$. Thus the condition $x=(x_{i})\in X_{[N]}$, $[N]\in T^{o}\cup D_{T}^{o}$ implies the existence of
$1\leq l<m<n\leq 6$ with $x_{l}x_{m}x_{n}\not=0$.
\par
Since the Jacobian matrix of the defining equations of $X_{[N]}$ at $x$ is given by the $3\times6$-matrix 
$2N\cdot{\rm diag}(x_{1},\ldots,x_{6})$, since any two column vectors of $N$ are linearly independent, 
and since $3$ of $x_{1},\ldots,x_{6}$ are non-zero, there is a pair $1\leq p<q\leq3$ such that
$df_{p}\wedge df_{q}\not=0$ at $x$. 
For simplicity, assume that $df_{1}\wedge df_{2}\not=0$ at $x$.
Then $f_{1}$, $f_{2}$ are part of local coordinates of ${\bf P}^{5}$ around $x$.
Choose a non-vanishing $2\times2$-minor $\det\binom{n_{1i}\,n_{1j}}{n_{2i}\,n_{2j}}=n_{1i}n_{2j}-n_{1j}n_{2i}\not=0$. 
Then we get 
$$
(n_{1i}n_{2j}-n_{1j}n_{2i})\cdot f_{3}
\equiv
\sum_{k\not=i,j,\,1\leq k\leq6}\Delta_{ijk}(N)\,x_{k}^{2}
\mod f_{1},f_{2}.
$$
Hence, on a neighborhood of $(x,[N])\in{\mathcal X}|_{T^{o}\cup D_{T}^{o}}\cap\{x_{i}x_{j}\not=0\}$, 
the total space ${\mathcal X}|_{T^{o}\cup D_{T}^{o}}$ is locally isomorphic to the hypersurface of ${\bf P}^{3}\times(T^{o}\cup D_{T}^{o})$
defined by
\begin{equation}
\label{eqn:localequation:X}
\Delta_{k_{1}ij}(N)\,z_{1}^{2}+\Delta_{k_{2}ij}(N)\,z_{2}^{2}+\Delta_{k_{3}ij}(N)\,z_{3}^{2}+\Delta_{k_{4}ij}(N)\,z_{4}^{2}=0,
\end{equation}
where $\{i,j\}\cup\{k_{1},k_{2},k_{3},k_{4}\}=\{1,\ldots,6\}$.
Since $(\Delta_{ijk}(N))_{i<j<k}$ are the Pl\"ucker coordinates of $T$ and since at most one of $\Delta_{ijk}(N)$ can vanish
on $T^{o}\cup D_{T}^{o}$, we deduce from \eqref{eqn:localequation:X} the assertions (1), (2).
Since $T_{ijk}$ is a Schubert cycle and hence is an irreducible ample divisor on $T$, we get $\dim[T_{ijk}\cap T_{lmn}]=\dim T-2$ 
when $T_{ijk}\not=T_{lmn}$. This proves (3).
\end{pf}

\subsection
{A natural pluricanonical differential associated to a $(2,2,2)$-model}
\label{sect:3.2}
\par

\subsubsection
{Equivalence of $(2,2,2)$-models}
\par
We set 
$$
GL({\bf C}^{6})^{\iota}:=\{\gamma\in GL({\bf C}^{6});\,\gamma\iota=\iota\gamma\}= GL({\bf C}^{3})\times GL({\bf C}^{3}).
$$
For $A=(A_{1},A_{2},A_{3})\in{\rm Sym}(3,{\bf C})\otimes{\bf C}^{3}$ and $P=(p_{ij})\in GL({\bf C}^{3})$, we define
$$
\begin{array}{ll}
A\cdot P
&:=
A\cdot(I_{{\rm Sym}(3,{\bf C})}\otimes P)
=
(\sum_{i=1}^{3}A_{i}p_{i1},\sum_{i=1}^{3}A_{i}p_{i2},\sum_{i=1}^{3}A_{i}p_{i3}),
\\
A^{P}
&:=
({}^{t}PA_{1}P,{}^{t}PA_{2}P,{}^{t}PA_{3}P).
\end{array}
$$
By definition, $s(A,B)=s(A',B')$ in $S$ if and only if $(A',B')=(A\cdot P,B\cdot P)$ for some $P\in GL({\bf C}^{3})$. 
Since every isomorphism $(X_{(A,B)},\iota,{\mathcal O}_{{\bf P}^{5}}(1))\cong(X_{(A',B')},\iota,{\mathcal O}_{{\bf P}^{5}}(1))$ 
is induced by a projective transform of ${\bf P}^{5}$ associated to an element of $GL({\bf C}^{6})^{\iota}$, we see that
$(Y_{(A,B)},{\mathcal L}_{(A,B)})\cong(Y_{(A',B')},{\mathcal L}_{(A',B')})$ if and only if 
$Y_{(A',B')}=Y_{(A^{P},B^{P'})}$ for some $(P,P')\in GL({\bf C}^{6})^{\iota}$.

\subsubsection
{A pluricanonical differential on $X_{s}$ and its invariance property}
\par
Define
$$
D^{+}:=\{s\in D^{*};\,{\rm Sing}(X_{s})\hbox{ consists of RDP's}\}.
$$
If $s\in D^{+}$, then $X_{s}$ is a $K3$ surface with RDP's. 
Since the singular fiber of any flat deformation of a RDP contains at most RDP's, $S^{o}\cup D^{+}$ is a Zariski open subset of $S$.
\par
For an admissible $(A,B)\in{\rm Sym}(3,{\bf C})\otimes{\bf C}^{6}$, we define 
$$
\varXi_{X_{(A,B)}}
:=
R(A)R(B)\,\omega_{(A,B)}^{\otimes8}
\in 
H^{0}\left(X_{(A,B)},K_{X_{(A,B)}}^{\otimes8}\right).
$$
Then $\varXi_{X_{(A,B)}}$ is an $\iota$-invariant pluricanonical form of weight $8$ on $X_{(A,B)}$
and is identified with a pluricanonical form of weight $8$ on $Y_{(A,B)}$.

\begin{lemma}
\label{lemma:pluricanonical:differential:1}
If $s(A,B)=s(A',B')\in S^{o}\cup D^{+}$, then
$$
\varXi_{X_{(A,B)}}=\varXi_{X_{(A',B')}}.
$$
\end{lemma}

\begin{pf}
We can write $(A',B')=(A\cdot P,B\cdot P)$, $P\in GL({\bf C}^{3})$.
By \eqref{eqn:2-form:residue}, \eqref{eqn:formula:varXi}, we get 
$$
\omega_{(A\cdot P,B\cdot P)}=(\det P)^{-1}\omega_{(A,B)}.
$$
By the explicit formula \cite[p.215, Table 1]{KSZ92} for $R(A)$, we get
$$
R(A\cdot P)=(\det P)^{4}R(A).
$$
The result follows from these two equalities.
\end{pf}

After Lemma~\ref{lemma:pluricanonical:differential:1}, it makes sense to define for $s=s(A,B)\in S^{o}\cup D^{+}$
$$
\varXi_{s}:=\varXi_{X_{(A,B)}}.
$$
Then we get a section $\varXi\in H^{0}(S^{o}\cup D^{+},\pi_{*}K_{{\mathcal X}/S}^{\otimes8})$ such that
$\varXi(s):=\varXi_{s}$. Here $K_{{\mathcal X}/S}$ denotes the relative dualizing sheaf of the family
$\pi\colon{\mathcal X}|_{S^{o}\cup D^{+}}\to S^{o}\cup D^{+}$.
\par
The group $GL({\bf C}^{6})^{\iota}=GL({\bf C}^{3})\times GL({\bf C}^{3})$ acts on ${\mathcal X}$ and $S$ by
$$
\begin{array}{ll}
(P,P')\cdot(\binom{x}{y},s)
&:=
(\binom{Px}{P'y},(P,P')\cdot s),
\\
(P,P')\cdot s(A,B)
&:=
s(A^{P},B^{P'}).
\end{array}
$$
Then the projection $\pi\colon{\mathcal X}\to S$ is $GL({\bf C}^{6})^{\iota}$-equivariant,
and $\pi_{*}K_{{\mathcal X}/S}^{\otimes8}$ is equipped with the structure of a $GL({\bf C}^{6})^{\iota}$-equivariant
line bundle on $S^{o}\cup D^{+}$.

\begin{lemma}
\label{lemma:pluricanonical:differential:2}
The section $\varXi\in H^{0}(S^{o}\cup D^{+},\pi_{*}K_{{\mathcal X}/S}^{\otimes8})$ is $GL({\bf C}^{6})^{\iota}$-invariant.
\end{lemma}

\begin{pf}
Let $g=(P,P')\in GL({\bf C}^{6})^{\iota}$ and $s=s(A,B)\in S^{o}\cup D^{+}$.
The projective transform $T_{g}\binom{x}{y}:=\binom{Px}{P'y}$ of ${\bf P}^{5}$ induces
an isomorphism $T_{g}\colon X_{g\cdot s}\to X_{s}$.
Since $X_{s}$ and $X_{g\cdot s}$ are $K3$ surfaces with possibly RDP's, 
there exists $\chi(g)\in{\bf C}^{*}$ with $T_{g}^{*}\omega_{(A,B)}=\chi(g)\,\omega_{(A^{P},B^{P'})}$.
Since $T_{gg'}=T_{g}T_{g'}$, we see that $\chi\colon GL({\bf C}^{6})^{\iota}\to{\bf C}^{*}$ is a character. 
Letting $g$ be a scalar matrix, we get $\chi(g)=\det(g)$. Hence
$$
(T_{(P,P')})^{*}\omega_{(A,B)}=\det(P)\det(P')\,\omega_{(A^{P},B^{P'})}.
$$
It is classical that
$$
R(A^{P})=\det(P)^{8}R(A).
$$
By these two equalities, we get $T_{g}^{*}\varXi_{g\cdot s}=\varXi_{s}$. This proves the lemma.
\end{pf}

We define the $GL({\bf C}^{6})^{\iota}$-invariant divisor ${\frak R}$ on $S^{o}\cup D^{o}$ by
$$
{\frak R}:={\rm div}(\varXi)=\{s=s(A,B)\in S^{o}\cup D^{o};\,R(A)R(B)=0\}.
$$
By the definition of resultants, $R(A)R(B)=0$ if and only if $X_{(A,B)}^{\iota}\not=\emptyset$.
Hence $\varXi_{s}\not=0$ if and only if $\iota|_{X_{s}}$ is free from fixed points, i.e., 
$Y_{s}=X_{s}/\iota$ is an Enriques surface (with possibly RDP's). 
By this interpretation, we get
$$
{\frak R}=\{s\in S^{o}\cup D^{o};\,X_{s}^{\iota}\not=\emptyset\}.
$$ 
We define 
$$
{\frak R}^{o}:=\{s\in{\frak R};\,\#X_{s}^{\iota}=1\}.
$$

\begin{lemma}
\label{lemma:zero:divisor:pluricanonical:differential}
The following hold:
\begin{itemize}
\item[(1)]
${\frak R}\subset D^{o}$.
\item[(2)]
${\frak R}^{o}$ is a dense Zariski open subset of ${\frak R}$.
\end{itemize}
\end{lemma}

\begin{pf}
{\bf (1) }
Since $X_{s}^{\iota}\subset{\rm Sing}(X_{s})$ for any $s\in S$, we get (1).
\par{\bf (2) }
It is clear that ${\frak R}^{o}$ is a Zariski open subset of ${\frak R}$. We see the density of ${\frak R}^{o}$ in ${\frak R}$.
Let ${\frak o}\in{\frak R}$ be an arbitrary point. It suffices to see that there is a curve 
$\gamma\colon(\varDelta,0)\to({\frak R},{\frak o})$ such that $\gamma(\varDelta\setminus\{0\})\subset{\frak R}^{o}$.
Write ${\frak o}=s(A,B)$. Since $X_{(A,B)}^{\iota}\not=\emptyset$, either the system of equations
$Q(x;A_{1})=Q(x;A_{2})=Q(x;A_{3})=0$ or $Q(y;B_{1})=Q(y;B_{2})=Q(y;B_{3})=0$ has a non-trivial solution.
Assume that the former has a non-trivial solution $c\in{\bf P}^{2}$.
Since generic three conics of ${\bf P}^{2}$ have no points in common, 
we can find a vector-valued holomorphic function $\varDelta\ni t\to(A(t),B(t))\in{\rm Sym}(3,{\bf C})\otimes{\bf C}^{6}$
satisfying the following conditions:
\begin{itemize}
\item[(i)]
$A_{i}(0)=A_{i}$ and $B_{i}(0)=B_{i}$ for $i=1,2,3$.
\item[(ii)]
For $t\not=0$, $\{x\in{\bf P}^{2};\,Q(x;A_{1}(t))=Q(x;A_{2}(t))=Q(x;A_{3}(t))=0\}=\{c\}$.
\item[(iii)]
For $t\not=0$, $\{x\in{\bf P}^{2};\,Q(x;B_{1}(t))=Q(x;B_{2}(t))=Q(x;B_{3}(t))=0\}=\emptyset$.
\end{itemize}
Set $\gamma(t):=s(A(t),B(t))\in S$. Since $S^{o}\cup D^{o}=S\setminus E$ is a Zariski open subset of $S$
by Lemma~\ref{lemma:Grassmann:discriminant} (2), we may assume $\gamma(\varDelta)\subset S^{o}\cup D^{o}$. 
We get $\gamma(0)={\frak o}$ by (i) and
$\gamma(\varDelta\setminus\{0\})\subset{\frak R}^{o}$ by (ii), (iii). This proves (2).
\end{pf}

\subsubsection
{A canonical Hermitian metric on the space of pluricanonical forms}
\label{subsubsect:3.4.3}
\par
Let $X$ be a $K3$ surface with possibly RDP's. 
For every $\nu\in{\bf Z}_{>0}$, $H^{0}\left(X,K_{X}^{\otimes\nu}\right)$ is equipped with
the Hermitian structure $\|\cdot\|_{L^{2/\nu}}$, which depends only on the complex structure on $X$, such that
$$
\|\xi\|_{L^{2/\nu}}:=\left(\int_{X}|\xi\wedge\overline{\xi}|^{1/\nu}\right)^{\nu/2},
\qquad
\xi\in H^{0}\left(X,K_{X}^{\otimes\nu}\right).
$$
When $s=s(A,B)\in S^{o}\cup D^{+}$, we get
$$
\|\varXi_{s}\|_{L^{1/4}}=\|\varXi_{X_{(A,B)}}\|_{L^{1/4}}
=
|R(A)R(B)|\cdot\left|\int_{X_{(A,B)}}\omega_{(A,B)}\wedge\overline{\omega}_{(A,B)}\right|^{4}.
$$
The Hermitian structure on the the invertible sheaf $\pi_{*}K_{{\mathcal X}/S}^{\otimes\nu}|_{S^{o}\cup D^{+}}$ induced from 
this fiberwise Hermitian structure is also denoted by $\|\cdot\|_{L^{2/\nu}}$.
\par
By Lemma~\ref{lemma:pluricanonical:differential:2},
the norm $\|\varXi\|_{L^{1/4}}$ is a $GL({\bf C}^{6})^{\iota}$-invariant, nowhere-vanishing $C^{\infty}$ function on $S^{o}$.
By the existence of simultaneous resolution of the family $\pi\colon{\mathcal X}|_{S^{o}\cup D^{+}}\to S^{o}\cup D^{+}$,
$\|\varXi\|_{L^{1/4}}$ extends to a $C^{0}$ function on $S^{o}\cup D^{+}$.

\subsection
{\bf A comparison of $\|\Phi\|$ and $\|\varXi\|_{L^{1/4}}$ over $S$}
\label{subsect:Comparison:Phi:varXi}
\par
Recall that ${\mathcal Y}={\mathcal X}/\iota$. Let 
$$
\overline{\varpi}\colon S^{o}\to{\mathcal M}^{o}
$$ 
be the period map for the family of Enriques surfaces $\pi\colon{\mathcal Y}|_{S^{o}}\to S^{o}$.
Then $\overline{\varpi}$ is $GL({\bf C}^{6})^{\iota}$-equivariant and $\overline{\varpi}$ is dominant.
By the existence of simultaneous resolution of the family $\pi\colon{\mathcal X}|_{S^{o}\cup D^{+}}\to S^{o}\cup D^{+}$,
$\overline{\varpi}$ extends to a holomorphic map from $S^{o}\cup D^{+}$ to ${\mathcal M}$. 
This extension is again denoted by $\overline{\varpi}$. Then, for $s\in D^{+}$, we have
$$
\|\Phi(Y_{s})\|:=\|\Phi(\varpi(s))\|=\|\Phi(\widetilde{Y}_{s})\|,
$$
where $\widetilde{Y}_{s}\to Y$ is the minimal resolution.
Since $s\in{\frak R}$ if and only if $Y_{s}$ is not an Enriques surface, i.e., $\overline{\varpi}(s)\in\overline{\mathcal D}$,
we get on $S^{o}\cup D^{o}$
\begin{equation}
\label{eqn:support:zero:Phi}
{\rm Supp}(\overline{\varpi}^{*}{\rm div}(\Phi))={\rm Supp}(\overline{\varpi}^{*}\overline{\mathcal D})={\rm Supp}({\frak R}).
\end{equation}

\begin{theorem}
\label{thm:comparison:Phi:pluricanonical:differential}
The following equality of $GL({\bf C}^{6})^{\iota}$-invariant functions on $S^{o}\cup D^{+}$ holds
$$
\overline{\varpi}^{*}\|\Phi\|^{2}=(2\pi^{-4})^{4}\,\|\varXi\|_{L^{1/4}}.
$$
\end{theorem}

\begin{pf}
{\em (Step 1) }
Recall that $\kappa$ is the K\"ahler form of the Bergman metric on ${\mathcal M}$.
Let $\delta_{\frak R}$ be the Dirac $\delta$-current on $S^{o}\cup D^{o}$ associated to the divisor ${\frak R}$ on $S^{o}\cup D^{o}$. 
Since $\log\|\omega_{(A,B)}\|_{L^{2}}^{2}$ is a local potential function of $\overline{\varpi}^{*}\kappa$,
we get the following equation of currents on $S^{o}\cup D^{o}$ by the Poincar\'e-Lelong formula
\begin{equation}
\label{eqn:curvature:log:G}
-dd^{c}\log\|\varXi\|_{L^{1/4}}=4\overline{\varpi}^{*}\kappa-\frac{1}{2}\delta_{\frak R}.
\end{equation}
\par{\em (Step 2) }
Let ${\frak o}\in{\frak R}^{o}\setminus{\rm Sing}({\frak R}^{o})$ be an arbitrary point.
Let $C\subset S$ be a compact Riemann surface intersecting ${\frak R}^{o}$ transversally at ${\frak o}$.
Let $(\varDelta,t)$ be a coordinate neighborhood of $C$ centered at ${\frak o}$ such that $\varDelta\setminus\{0\}\subset C\cap S^{o}$.
For $t\in\varDelta$, let $s(t)\in S$ be the corresponding point in $S$. Write $s(t)=s(A(t),B(t))$, 
where the map $\varDelta\ni t\mapsto(A(t),B(t))\in{\rm Sym}(3,{\bf C})\otimes{\bf C}^{6}$ is holomorphic. 
Since $s(0)={\frak o}\in{\frak R}^{o}$ and $s(\varDelta\setminus\{0\})\subset S^{o}$, 
the flat family $\pi\colon{\mathcal X}|_{\varDelta}\to\varDelta$ has the following properties: 
\begin{itemize}
\item[(i)]
$X_{s(t)}$ is smooth for $t\in\varDelta\setminus\{0\}$ and ${\rm Sing}(X_{\frak o})$ consists of nodes.
\item[(ii)]
${\mathcal X}^{\iota}|_{\varDelta}=X^{\iota}_{\frak o}\subset{\rm Sing}(X_{\frak o})$ and $\#X^{\iota}_{\frak o}=1$.
\end{itemize}
We see that ${\mathcal X}|_{C}$ is smooth on a neighborhood of $X_{\frak o}$ 
if the curve germ $(C,0)\subset(S,{\frak o})$ is generic. It suffices to prove the smoothness of ${\mathcal X}$ at ${\rm Sing}(X_{\frak o})$.
\par
Let ${\frak p}=(x_{0},y_{0})\in{\rm Sing}(X_{\frak o})$.
Set $f_{i}(x,y;t):=Q(x;A_{i}(t))+Q(y;B_{i}(t))$ for $i=1,2,3$.
For $(\widetilde{A},\widetilde{B})\in{\rm Sym}(3,{\bf C})\otimes{\bf C}^{6}$, we define 
$$
J_{{\frak p},(\widetilde{A},\widetilde{B})}
=
(\partial_{x_{j}}f_{i}(x_{0},y_{0};0),\partial_{y_{k}}f_{i}(x_{0},y_{0};0),Q(x_{0};\widetilde{A}_{i})+Q(y_{0};\widetilde{B}_{i}))_{i,j,k=1,2,3}.
$$ 
Since $X_{s(t)}$ is defined by the system of equations $f_{i}(x,y;t)=Q(x;A_{i}(t))+Q(y;B_{i}(t))$ $(i=1,2,3)$ and hence
the corresponding Jacobian matrix at ${\frak p}=(x_{0},y_{0})$ is given by the $3\times7$-matrix $J_{{\frak p},(\partial_{t}A(0),\partial_{t}B(0))}$,
${\mathcal X}|_{C}$ is smooth at ${\frak p}=(x_{0},y_{0})$ if ${\rm rank}\,J_{{\frak p},(\partial_{t}A(0),\partial_{t}B(0))}=3$.
\par
Since ${\frak p}=(x_{0},y_{0})\in X_{\frak o}$ is a node by (i) and hence ${\rm emb}(X_{\frak o},{\frak p})=3$, we get 
\begin{equation}
\label{eqn:Jacobi:criterion:node}
{\rm rank}(\partial_{x_{j}}f_{i}(x_{0},y_{0};0),\partial_{y_{k}}f_{i}(x_{0},y_{0};0))=2
\end{equation}
by Jacobi's criterion.
Since $(x_{0},y_{0})\not=(0,0)$ and hence the linear map
$$
{\rm Sym}(3,{\bf C})\otimes{\bf C}^{6}\ni(\widetilde{A},\widetilde{B})\to(Q(x_{0};\widetilde{A}_{i})+Q(y_{0};\widetilde{B}_{i}))_{i=1,2,3}\in{\bf C}^{3}
$$
is surjective, we deduce from \eqref{eqn:Jacobi:criterion:node} that ${\rm rank}\,J_{{\frak p},(\widetilde{A},\widetilde{B})}=3$ for generic 
$(\widetilde{A},\widetilde{B})\in{\rm Sym}(3,{\bf C})\otimes{\bf C}^{6}$.
This implies that the subset of ${\rm Sym}(3,{\bf C})\otimes{\bf C}^{6}$ defined by
$$
{\frak V}_{\frak p}
:=
\{(A,B)\in{\rm Sym}(3,{\bf C})\otimes{\bf C}^{6};\,{\rm rank}\,J_{{\frak p},(\widetilde{A},\widetilde{B})}<3
\}
$$
is a proper Zariski closed subset.
Hence ${\mathcal X}$ is smooth at every point of ${\rm Sing}(X_{\frak o})$ 
if $(\partial_{t}A(0),\partial_{t}B(0))\in[{\rm Sym}(3,{\bf C})\otimes{\bf C}^{6}]\setminus\bigcup_{{\frak p}\in{\rm Sing}(X_{\frak o})}{\frak V}_{p}$.
This proves that ${\mathcal X}|_{C}$ is smooth on a neighborhood of  $X_{\frak o}$ for a generic curve germ $(C,0)\subset(S,{\frak o})$.
\par{\em (Step 3) }
Let $(C,0)\subset(S,{\frak o})$ be a generic curve germ such that ${\mathcal X}|_{C}$ is smooth on a neighborhood of  $X_{\frak o}$.
By Step 2 (i), (ii), we may apply Theorem~\ref{thm:degeneration:Phi} to the deformation germ of $K3$ surfaces with involution 
$f\colon({\mathcal X}|_{C},\iota)\to C$. Then we get
\begin{equation}
\label{eqn:asymptitic:Phi}
\log\|\Phi\|^{2}|_{C}(t)=\frac{1}{2}\log|t|^{2}+O(1)
\qquad
(t\to0).
\end{equation}
Since ${\frak R}^{o}$ is a dense Zariski open subset of ${\frak R}$ by 
Lemma~\ref{lemma:zero:divisor:pluricanonical:differential} (2), we get the following equation of currents on $S^{o}\cup D^{o}$
by \eqref{eqn:Poincare:Lelong:Phi:moduli}, \eqref{eqn:support:zero:Phi}, \eqref{eqn:asymptitic:Phi}
\begin{equation}
\label{eqn:current:equation:norm:Phi}
-dd^{c}\log(\overline{\varpi}^{*}\|\Phi\|^{2})=4\overline{\varpi}^{*}\kappa-\frac{1}{2}\delta_{\frak R}.
\end{equation}
Comparing \eqref{eqn:curvature:log:G} and \eqref{eqn:current:equation:norm:Phi}, we get the following equation of currents 
on $S^{o}\cup D^{o}$
$$
dd^{c}\log(\overline{\varpi}^{*}\|\Phi\|^{2}/\|\varXi\|_{L^{1/4}})=0,
$$
which implies that $\partial\log(\overline{\varpi}^{*}\|\Phi\|^{2}/\|\varXi\|_{L^{1/4}})$ is a holomorphic $1$-form 
on $S^{o}\cup D^{o}$.
Since $\dim[S\setminus(S^{o}\cup D^{o})]=\dim E\leq\dim S-2$ by Lemma~\ref{lemma:Grassmann:discriminant} (2),
it follows from the Hartogs extension principle that $\partial\log(\overline{\varpi}^{*}\|\Phi\|^{2}/\|\varXi\|_{L^{1/4}})$ 
extends to a holomorphic $1$-form on the Grassmann variety $S$. 
Since $H^{0}(S,\Omega^{1}_{S})=0$ by the rationality of $S$, we get 
$\partial\log(\overline{\varpi}^{*}\|\Phi\|^{2}/\|\varXi\|_{L^{1/4}})=0$ on $S$. 
Since $\overline{\varpi}^{*}\|\Phi\|^{2}/\|\varXi\|_{L^{1/4}}$ is real-valued, 
we get $d\log(\overline{\varpi}^{*}\|\Phi\|^{2}/\|\varXi\|_{L^{1/4}})=0$ on $S$. 
This proves the existence of a constant ${\frak C}\in{\bf R}$ such that the following equality of functions on $S$ holds
\begin{equation}
\label{eqn:constant:Phi:generic}
\overline{\varpi}^{*}\|\Phi\|^{2}={\frak C}(2\pi^{-4})^{4}\,\|\varXi\|_{L^{1/4}}.
\end{equation}
Restricting \eqref{eqn:constant:Phi:generic} to a certain subset of $S$, we get ${\frak C}=1$ by 
Theorem~\ref{theorem:constant:Thomae:Kummer:product:type} below.
\end{pf}

\par
Let $\lambda$ be the Hodge bundle on ${\mathcal M}$. 
By Lemma~\ref{cor:character:Phi:square}, $\Phi^{2}$ is a holomorphic section of $\lambda^{\otimes8}$. 
Regarding $\Phi^{2}\in H^{0}({\mathcal M},\lambda^{\otimes8})$, we get ${\rm div}(\Phi^{2})=\overline{\mathcal D}$.
To emphasize, denote by $\|\cdot\|_{\rm Pet}$ the Petersson norm on $\lambda$.
Theorem~\ref{thm:comparison:Phi:pluricanonical:differential} implies the following.

\begin{corollary}
\label{cor:Phi:pluri-canonical:form}
There exists a $GL({\bf C}^{6})^{\iota}$-equivariant holomorphic isometry 
$$
f\colon\left(\overline{\varpi}^{*}(\lambda^{\otimes8}),\overline{\varpi}^{*}\|\cdot\|_{\rm Pet}^{\otimes8}\right)
\cong 
\left(\pi_{*}(K_{{\mathcal X}/S}^{\otimes8}),\|\cdot\|_{L^{1/4}}\right)
$$ 
of holomorphic Hermitian line bundles on $S^{o}\cup D^{+}$ such that
$$
f(\overline{\varpi}^{*}\Phi^{2})=(2\pi^{-4})^{4}\,\varXi.
$$
\end{corollary}

\begin{pf}
By Theorem~\ref{thm:comparison:Phi:pluricanonical:differential}, $(\overline{\varpi}^{*}\Phi^{2})^{-1}\otimes\varXi$
is a nowhere vanishing $GL({\bf C}^{6})^{\iota}$-invariant holomorphic section of 
$\overline{\varpi}^{*}(\lambda^{\otimes8})^{-1}\otimes\pi_{*}(K_{{\mathcal X}/S}^{\otimes8})$ on $S^{o}\cup D^{+}$.
Defining $f(\xi):=(\overline{\varpi}^{*}\Phi^{2})^{-1}\otimes\varXi\otimes\xi$, we see that
$\overline{\varpi}^{*}(\lambda^{\otimes8})\cong\pi_{*}(K_{{\mathcal X}/S}^{\otimes8})$ on $S^{o}\cup D^{+}$ via $f$. 
The desired equality follows from Theorem~\ref{thm:comparison:Phi:pluricanonical:differential}. 
\end{pf}

{\bf Proof of Theorem~\ref{thm:valuePhi:general} }
If $(A,B)\in{\rm Sym}(3,{\bf C})\otimes{\bf C}^{6}$ is admissible, then $s(A,B)\in S^{o}\cup D^{+}$ and $X_{(A,B)}=X_{s(A,B)}$.
Theorem~\ref{thm:valuePhi:general} follows from Theorem~\ref{thm:comparison:Phi:pluricanonical:differential}.
\qed

\section
{An algebraic expression of the Borcherds $\Phi$-function}
\label{sect:Thomae:Phi}
\par
In this section, we fix the following notation.
Let $(A,B)\in{\rm Sym}(3,{\bf C})\otimes{\bf C}^{6}$ be admissible. Hence $s(A,B)\in S^{o}\cup D^{+}$.
Let ${\bf v}\in H^{2}(X_{(A,B)},{\bf Z})_{-}$ be a primitive isotropic vector of level $\ell$. 
Here the level of ${\bf v}$ is the positive generator of $\langle{\bf v},\LAM\rangle\subset{\bf Z}$.
Then $\ell=1$ or $\ell=2$ by \cite[Prop.\,4.5]{Sterk91}.
Let ${\bf v}'\in H^{2}(X_{(A,B)},{\bf Z})_{-}$ be another primitive isotropic vector of level $\ell$ such that $\langle{\bf v},{\bf v}'\rangle=\ell$
and set 
$$
{\frak L}_{{\bf v},{\bf v}'}:=({\bf Z}{\bf v}+{\bf Z}{\bf v}')^{\perp_{H^{2}(X_{(A,B)},{\bf Z})_{-}}}.
$$
We get the orthogonal decomposition $H^{2}(X_{(A,B)},{\bf Z})_{-}=({\bf Z}{\bf v}+{\bf Z}{\bf v}')\oplus{\frak L}_{{\bf v},{\bf v}'}$, from which 
${\frak L}_{{\bf v},{\bf v}'}$ is $2$-elementary with ${\frak L}_{{\bf v},{\bf v}'}\cong{\Bbb M}_{\ell}={\Bbb U}(2/\ell)\oplus{\Bbb E}_{8}(2)$.
We define
$$
(-1)^{2/\ell}z_{A,B,{\bf v},{\bf v}'}
:=
\frac{\omega_{(A,B)}
-\langle\omega_{(A,B)},{\bf v}'/\ell\rangle{\bf v}
-\langle\omega_{(A,B)},{\bf v}\rangle({\bf v}'/\ell)}
{\langle\omega_{(A,B)},{\bf v}\rangle}.
$$
Then $z_{A,B,{\bf v},{\bf v}'}\in{\frak L}_{{\bf v},{\bf v}'}\otimes{\bf R}+i\,\Ccal_{{\frak L}_{{\bf v},{\bf v}'}}$,
and the following equality holds:
$$
\omega_{(A,B)}/\langle\omega_{(A,B)},{\bf v}\rangle=-\frac{(z_{A,B,{\bf v},{\bf v}'})^{2}}{2}{\bf v}+\frac{{\bf v}'}{\ell}+(-1)^{2/\ell}z_{A,B,{\bf v},{\bf v}'}.
$$
Let $\alpha_{{\bf v},{\bf v}'}\colon{\frak L}_{{\bf v},{\bf v}'}\cong{\Bbb M}_{\ell}$ be an isometry of lattices, called a marking of ${\frak L}_{{\bf v},{\bf v}'}$, 
such that $\alpha_{{\bf v},{\bf v}'}(z_{A,B,{\bf v},{\bf v}'})\in{\Bbb M}_{\ell}\otimes{\bf R}+i\,{\mathcal C}_{{\Bbb M}_{\ell}}^{+}$.
(If $\alpha_{{\bf v},{\bf v}'}(z_{A,B,{\bf v},{\bf v}'})\in{\Bbb M}_{\ell}\otimes{\bf R}-i\,{\mathcal C}_{{\Bbb M}_{\ell}}^{+}$, 
then we replace $\alpha_{{\bf v},{\bf v}'}$ by $-\alpha_{{\bf v},{\bf v}'}$.)
\par
Recall that $\{{\frak e}_{\ell},{\frak f}_{\ell}\}$ is a basis of ${\Bbb U}(\ell)$ with ${\frak e}_{\ell}^{2}={\frak f}_{\ell}^{2}=0$ and 
$\langle{\frak e}_{\ell},{\frak f}_{\ell}\rangle=\ell$.
We extend $\alpha_{{\bf v},{\bf v}'}$ to an isometry $\widetilde{\alpha}_{{\bf v},{\bf v}'}\colon H^{2}(X_{(A,B)},{\bf Z})_{-}\cong\LAM$
by setting 
$$
\widetilde{\alpha}_{{\bf v},{\bf v}'}(m{\bf v}+n{\bf v}'+x):=m{\frak e}_{\ell}+n{\frak f}_{\ell}+\alpha_{{\bf v},{\bf v}'}(x),
$$
where $m,n\in{\bf Z}$, $x\in{\frak L}_{{\bf v},{\bf v}'}$. 
Hence $\widetilde{\alpha}_{{\bf v},{\bf v}'}({\bf v})={\frak e}_{\ell}$ and $\widetilde{\alpha}_{{\bf v},{\bf v}'}({\bf v}')={\frak f}_{\ell}$.
Since 
$$
\widetilde{\alpha}_{{\bf v},{\bf v}'}(\omega_{(A,B)})/\langle\omega_{(A,B)},{\bf v}\rangle
=
-\frac{({\alpha}_{{\bf v},{\bf v}'}(z))^{2}}{2}{\frak e}_{\ell}+\frac{{\frak f}_{\ell}}{\ell}+(-1)^{2/\ell}{\alpha}_{{\bf v},{\bf v}'}(z)
=
\iota_{\ell}(\alpha_{{\bf v},{\bf v}'}(z)),
$$
$\alpha_{{\bf v},{\bf v}'}(z_{A,B,{\bf v},{\bf v}'})$ is the period of the marked Enriques surface $(Y_{(A,B)},\widetilde{\alpha}_{{\bf v},{\bf v}'})$
under the isomorphism \eqref{eqn:tube:domain:level1}.
In this section, we prove the following theorem.

\begin{theorem}
\label{theorem:Thomae:Phi}
Let ${\bf v}^{\lor}\in H_{2}(X_{(A,B)},{\bf Z})$ be the Poincar\'e dual of ${\bf v}$.
Then the following equality holds:
$$
\Phi_{\ell}\left(\alpha_{{\bf v},{\bf v}'}(z_{A,B,{\bf v},{\bf v}'})\right)^{2} 
=
R(A)R(B)
\left(
\frac{2}{\pi^{2}}
\int_{{\bf v}^{\lor}}\omega_{(A,B)}
\right)^{8}. 
$$ 
\end{theorem}

For the proof, we need some intermediary results.

\begin{lemma}
\label{lemma:Thomae:Phi:absolute:value}
The following equality holds:
\begin{equation}
\label{eqn:Thomae:Phi:weak}
\left|
\Phi_{\ell}\left(\alpha_{{\bf v},{\bf v}'}(z_{A,B,{\bf v},{\bf v}'})\right)
\right|^{2} 
=
|R(A)R(B)|\cdot
\left|
\frac{2}{\pi^{2}}\int_{{\bf v}^{\lor}}\omega_{(A,B)}
\right|^{8}. 
\end{equation}
In particular, there exists a unique angle $\theta_{A,B,{\bf v},{\bf v}',\alpha_{{\bf v},{\bf v}'}}\in{\bf R}/2\pi{\bf Z}$ such that
\begin{equation}
\label{eqn:Thomae:Phi:constant}
\Phi_{\ell}\left(\alpha_{{\bf v},{\bf v}'}(z_{A,B,{\bf v},{\bf v}'})\right)^{2} 
=
e^{i\,\theta_{A,B,{\bf v},{\bf v}',\alpha_{{\bf v},{\bf v}'}}}\,R(A)R(B)\cdot
\left(
\frac{2}{\pi^{2}}\int_{{\bf v}^{\lor}}\omega_{(A,B)}
\right)^{8}. 
\end{equation}
\end{lemma}

\begin{pf}
For simplicity, write $z:=\alpha_{{\bf v},{\bf v}'}(z_{A,B,{\bf v},{\bf v}'})\in{\Bbb M}_{\ell}\otimes{\bf R}+i\,{\mathcal C}_{{\Bbb M}_{\ell}}^{+}$.
By the definition of the Petersson norm (cf. Sect~\ref{subsec:summary:Borcherds:Phi}), we get
\begin{equation}
\label{eqn:Petersson:norm:Phi}
\|\Phi(Y_{(A,B)})\|
=
\langle\Im z,\Im z\rangle_{{\Bbb M}_{\ell}}^{2}|\Phi_{\ell}(z)|.
\end{equation}
By \eqref{eqn:relation:bergman:tube:bounded} and the equality
$\langle\omega_{(A,B)},{\bf v}\rangle=\int_{{\bf v}^{\lor}}\omega_{(A,B)}$, we get
\begin{equation}
\label{eqn::L2norm:Petersson:norm}
\begin{aligned}
\int_{X_{(A,B)}}\omega_{(A,B)}\wedge\overline{\omega_{(A,B)}}
&=
|\langle\omega_{(A,B)},{\bf v}\rangle|^{2}
\left(
\int_{X_{(A,B)}}
\frac{\omega_{(A,B)}}
{\langle\omega_{(A,B)},{\bf v}\rangle}\wedge
\overline{\frac{\omega_{(A,B)}}{\langle\omega_{(A,B)},{\bf v}\rangle}}
\right)
\\
&=
|\langle\omega_{(A,B)},{\bf v}\rangle|^{2}
\left\langle
\iota_{\ell}(z),\overline{\iota_{\ell}(z)}
\right\rangle_{\LAM}
=
2\langle\Im z,\Im z\rangle_{{\Bbb M}_{\ell}}
\left|\int_{{\bf v}^{\lor}}\omega_{(A,B)}\right|^{2}.
\end{aligned}
\end{equation}
Substituting \eqref{eqn:Petersson:norm:Phi} and \eqref{eqn::L2norm:Petersson:norm} into the formula in Theorem~\ref{thm:valuePhi:general},
we get \eqref{eqn:Thomae:Phi:weak}. 
\end{pf}

\begin{lemma}
\label{lemma:transformation:rule:Phi}
Let ${\bf v}''\in H^{2}(X_{(A,B)},{\bf Z})_{-}$ be another primitive isotropic vector of level $\ell$ such that $\langle{\bf v},{\bf v}''\rangle=\ell$
and let $\alpha_{{\bf v},{\bf v}''}\colon{\frak L}_{{\bf v},{\bf v}''}\cong{\Bbb M}_{\ell}$ be an isometry of lattices such that
$\alpha_{{\bf v},{\bf v}''}(z_{A,B,{\bf v},{\bf v}''})\in{\Bbb M}_{\ell}\otimes{\bf R}+i\,{\mathcal C}_{{\Bbb M}_{\ell}}^{+}$. 
Then
$$
\Phi_{\ell}\left(\alpha_{{\bf v},{\bf v}'}(z_{A,B,{\bf v},{\bf v}'})\right)^{2}
=
\Phi_{\ell}\left(\alpha_{{\bf v},{\bf v}''}(z_{A,B,{\bf v},{\bf v}''})\right)^{2}.
$$
\end{lemma}

\begin{pf}
For simplicity, write $z_{{\bf v},{\bf v}'}$ (resp. $z_{{\bf v},{\bf v}''}$) for $z_{A,B,{\bf v},{\bf v}'}$ (resp. $z_{A,B,{\bf v},{\bf v}''}$).
Set $g_{{\bf v}',{\bf v}''}:=\widetilde{\alpha}_{{\bf v},{\bf v}'}\circ\widetilde{\alpha}_{{\bf v},{\bf v}''}^{-1}\in O^{+}(\LAM)$.
Since 
$$
\begin{array}{ll}
\widetilde{\alpha}_{{\bf v},{\bf v}'}(\omega_{(A,B)})
&=
g_{{\bf v}',{\bf v}''}(\widetilde{\alpha}_{{\bf v},{\bf v}''}(\omega_{(A,B)})),
\\
\widetilde{\alpha}_{{\bf v},{\bf v}'}(\omega_{(A,B)})/\langle\omega_{(A,B)},{\bf v}\rangle
&=
\iota_{\ell}(\alpha_{{\bf v},{\bf v}'}(z_{{\bf v},{\bf v}'})),
\\
\widetilde{\alpha}_{{\bf v},{\bf v}''}(\omega_{(A,B)})/\langle\omega_{(A,B)},{\bf v}\rangle
&=
\iota_{\ell}(\alpha_{{\bf v},{\bf v}''}(z_{{\bf v},{\bf v}''}))
\end{array}
$$
and hence 
$g_{{\bf v}',{\bf v}''}(\iota_{\ell}(\alpha_{{\bf v},{\bf v}''}(z_{{\bf v},{\bf v}''})))=\iota_{\ell}(\alpha_{{\bf v},{\bf v}'}(z_{{\bf v},{\bf v}'}))$,
it follows from the definition \eqref{eqn:group:action:domain:type:IV:level1} of the $O^{+}(\LAM)$-action on 
${\Bbb M}_{\ell}\otimes{\bf R}+i\,{\mathcal C}_{{\Bbb M}_{\ell}}^{+}$ that
$$
\alpha_{{\bf v},{\bf v}'}(z_{{\bf v},{\bf v}'})=g_{{\bf v}',{\bf v}''}\cdot\alpha_{{\bf v},{\bf v}''}(z_{{\bf v},{\bf v}''}).
$$
\par
Since $g_{{\bf v}',{\bf v}''}({\frak e}_{\ell})={\frak e}_{\ell}$, $g_{{\bf v}',{\bf v}''}$ induces an isometry
$\overline{g}_{{\bf v}',{\bf v}''}\in O({\frak e}_{\ell}^{\perp}/{\bf Z}{\frak e}_{\ell})$ of the lattice 
${\frak e}_{\ell}^{\perp}/{\bf Z}{\frak e}_{\ell}\cong{\Bbb M}_{\ell}$.
Since ${\frak f}_{\ell}/\ell\in{\Bbb M}_{\ell}^{\lor}$, there exists $\lambda\in{\Bbb M}_{\ell}^{\lor}$ such that
$$
g_{{\bf v}',{\bf v}''}
\left(
m{\frak e}_{\ell}+n\frac{{\frak f}_{\ell}}{\ell}+x
\right)
=
\left(
m-n\frac{(\lambda^{2})}{2}-\langle\lambda,\overline{g}_{{\bf v}',{\bf v}''}(x)\rangle
\right)
{\frak e}_{\ell}
+
n\frac{{\frak f}_{\ell}}{\ell}
+
\overline{g}_{{\bf v}',{\bf v}''}(x)
+
n\lambda
$$
for all $m,n\in{\bf Z}$ and $x\in{\Bbb M}_{\ell}$. 
By this expression and \eqref{eqn:automorphic:factor:domain:type:IV:level1}, 
we get $j_{\ell}(g_{{\bf v}',{\bf v}''},z)=1$ for all $z\in{\Bbb M}_{\ell}\otimes{\bf R}+i\,{\mathcal C}_{{\Bbb M}_{\ell}}^{+}$.
By \eqref{eqn:automorphic:property:Phi:level:N} and Lemma~\ref{cor:character:Phi:square}, we get
$$
\Phi_{\ell}(g_{{\bf v}',{\bf v}''}\cdot z)^{2}=\Phi_{\ell}(z)^{2}.
$$
Since $\Phi_{\ell}(\alpha_{{\bf v},{\bf v}'}(z_{{\bf v},{\bf v}'}))^{2}=\Phi_{\ell}(g_{{\bf v}',{\bf v}''}\cdot\alpha_{{\bf v},{\bf v}''}(z_{{\bf v},{\bf v}''}))^{2}$,
we get the desired equality $\Phi_{\ell}(\alpha_{{\bf v},{\bf v}''}(z_{{\bf v},{\bf v}''}))^{2}=\Phi_{\ell}(\alpha_{{\bf v},{\bf v}'}(z_{{\bf v},{\bf v}'}))^{2}$.
\end{pf}

After Lemma~\ref{lemma:transformation:rule:Phi} and \eqref{eqn:Thomae:Phi:constant}, 
$\theta_{{\bf v},{\bf v}',\alpha_{{\bf v},{\bf v}'}}\in{\bf R}/2\pi{\bf Z}$ depends only on $(A,B,{\bf v})$.
Set
$$
c(A,B,{\bf v})
:=
e^{2\pi i\,\theta_{{\bf v},{\bf v}',\alpha_{{\bf v},{\bf v}'}}}
=
\frac{\Phi_{\ell}\left(\alpha_{{\bf v},{\bf v}'}(z_{A,B,{\bf v},{\bf v}'})\right)^{2}}
{R(A)R(B)\left(2\pi^{-2}\int_{{\bf v}^{\lor}}\omega_{(A,B)}\right)^{8}}.
$$

\begin{lemma}
\label{lemma:independency:constant:level2:cocycle}
Let ${\bf w}\in H^{2}(X_{(A,B)},{\bf Z})_{-}$ be an arbitrary primitive isotropic vector of level $\ell$. 
Then 
$$
c(A,B,{\bf v})=c(A,B,{\bf w}).
$$
In particular, $c(A,B,{\bf v})$ depends only on $(A,B)$ and $\ell$.
\end{lemma}

\begin{pf}
Let ${\bf w}'\in H^{2}(X_{(A,B)},{\bf Z})_{-}$ be a primitive isotropic vectors of level $\ell$ with $\langle{\bf w},{\bf w}'\rangle=\ell$.
For simplicity, write $z_{{\bf v},{\bf v}'}$ (resp. $z_{{\bf w},{\bf w}'}$) for $z_{A,B,{\bf v},{\bf v}'}$ (resp. $z_{A,B,{\bf w},{\bf w}'}$).
Let $\alpha_{{\bf w},{\bf w}'}\colon{\frak L}_{{\bf w},{\bf w}'}\cong{\Bbb M}_{\ell}$ be a marking and
set $g_{{\bf v},{\bf v}';{\bf w},{\bf w}'}:=\widetilde{\alpha}_{{\bf v},{\bf v}'}\circ\widetilde{\alpha}_{{\bf w},{\bf w}'}^{-1}\in O^{+}(\LAM)$.
Since 
$$
\begin{array}{ll}
\widetilde{\alpha}_{{\bf v},{\bf v}'}(\omega_{(A,B)})
&=
g_{{\bf v},{\bf v}';{\bf w},{\bf w}'}(\widetilde{\alpha}_{{\bf w},{\bf w}'}(\omega_{(A,B)})),
\\
\widetilde{\alpha}_{{\bf v},{\bf v}'}(\omega_{(A,B)})/\langle\omega_{(A,B)},{\bf v}\rangle
&=
\iota_{\ell}(\alpha_{{\bf v},{\bf v}'}(z_{{\bf v},{\bf v}'})),
\\
\widetilde{\alpha}_{{\bf w},{\bf w}'}(\omega_{(A,B)})/\langle\omega_{(A,B)},{\bf w}\rangle
&=
\iota_{\ell}(\alpha_{{\bf w},{\bf w}'}(z_{{\bf w},{\bf w}'}))
\end{array}
$$
and hence
$g_{{\bf v},{\bf v}';{\bf w},{\bf w}'}(\iota_{\ell}(\alpha_{{\bf w},{\bf w}'}(z_{{\bf w},{\bf w}'})))=\iota_{\ell}(\alpha_{{\bf v},{\bf v}'}(z_{{\bf v},{\bf v}'}))$,
it follows from the definition \eqref{eqn:group:action:domain:type:IV:level1} of the $O^{+}(\LAM)$-action on 
${\Bbb M}_{\ell}\otimes{\bf R}+i\,{\mathcal C}_{{\Bbb M}_{\ell}}^{+}$ that $z_{{\bf v},{\bf v}'}=g_{{\bf v},{\bf v}';{\bf w},{\bf w}'}\cdot z_{{\bf w},{\bf w}'}$.
By \eqref{eqn:automorphic:factor:domain:type:IV:level1},
the automorphic factor for the Borcherds $\Phi$-function $\Phi_{\ell}$ is expressed as
$$
j_{\ell}(g,[\eta])=\langle g(\eta),{\frak e}_{\ell}\rangle/\langle\eta,{\frak e}_{\ell}\rangle,
\qquad
[\eta]\in\Omega^{+}_{\LAM},
\quad
g\in O^{+}(\LAM)
$$
under the identification $\Omega_{\LAM}^{+}\cong{\Bbb M}_{\ell}\otimes{\bf R}+i\,{\mathcal C}_{{\Bbb M}_{\ell}}^{+}$ given by \eqref{eqn:tube:domain:level1}.
We get
$$
\begin{aligned}
j_{\ell}\left(g_{{\bf v},{\bf v}';{\bf w},{\bf w}'},[\widetilde{\alpha}_{{\bf w},{\bf w}'}(\omega_{(A,B)})]\right)
&=
\frac
{\langle g_{{\bf v},{\bf v}';{\bf w},{\bf w}'}(\widetilde{\alpha}_{{\bf w},{\bf w}'}(\omega_{(A,B)})),{\frak e}_{\ell}\rangle}
{\langle\widetilde{\alpha}_{{\bf w},{\bf w}'}(\omega_{(A,B)}),{\frak e}_{\ell}\rangle}
\\
&=
\frac
{\langle\widetilde{\alpha}_{{\bf v},{\bf v}'}(\omega_{(A,B)}),{\frak e}_{\ell}\rangle}
{\langle\widetilde{\alpha}_{{\bf w},{\bf w}'}(\omega_{(A,B)}),{\frak e}_{\ell}\rangle}
=
\frac
{\langle\omega_{(A,B)},{\bf v}\rangle}
{\langle\omega_{(A,B)},{\bf w}\rangle}.
\end{aligned}
$$
By Lemma~\ref{cor:character:Phi:square} and the automorphic property of $\Phi_{\ell}^{2}$, we get
$$
\begin{aligned}
\Phi_{\ell}\left(
\alpha_{{\bf v},{\bf v}'}(z_{{\bf v},{\bf v}'})
\right)^{2}
&=
\Phi_{\ell}\left(
g_{{\bf v},{\bf v}';{\bf w},{\bf w}'}\cdot z_{{\bf w},{\bf w}'}
\right)^{2}
\\
&=
j_{\ell}\left(
g_{{\bf v},{\bf v}';{\bf w},{\bf w}'},[\widetilde{\alpha}_{{\bf w},{\bf w}'}(\omega_{(A,B)})]
\right)^{8}
\Phi_{\ell}\left(
\alpha_{{\bf w},{\bf w}'}(z_{{\bf w},{\bf w}'})
\right)^{2}
\\
&=
\left(
\frac{\langle\omega_{(A,B)},{\bf v}\rangle}{\langle\omega_{(A,B)},{\bf w}\rangle}
\right)^{8}
c(A,B,{\bf w})
R(A)R(B)
\left(
\frac{2}{\pi^{2}}
\int_{{\bf w}^{\lor}}\omega_{(A,B)}
\right)^{8}
\\
&=
c(A,B,{\bf w})
R(A)R(B)
\left(
\frac{2}{\pi^{2}}
\int_{{\bf v}^{\lor}}\omega_{(A,B)}
\right)^{8}.
\end{aligned}
$$
Comparing this with the definition of $c(A,B,{\bf v})$, we get $c(A,B,{\bf v})=c(A,B,{\bf w})$.
\end{pf}

After Lemma~\ref{lemma:independency:constant:level2:cocycle} and the relations in the proof of Lemma~\ref{lemma:pluricanonical:differential:1}, 
one can define a function $c_{\ell}(\cdot)$ on $S^{o}$ by
$$
c_{\ell}(s(A,B)):=c(A,B,{\bf v}).
$$

\begin{lemma}
\label{lemma:independency:constant:holomorphic:parameter}
The function $c_{\ell}(\cdot)$ on $S^{o}$ is constant. In particular, there exists a constant ${\frak c}_{\ell}$ with $|{\frak c}_{\ell}|={\frak C}$ such that
$$
\Phi_{\ell}\left(
\alpha_{{\bf v},{\bf v}'}(z_{A,B,{\bf v},{\bf v}'})
\right)^{2}
=
{\frak c}_{\ell}\,
R(A)R(B)
\left(
\frac{2}{\pi^{2}}
\int_{{\bf v}^{\lor}}\omega_{(A,B)}
\right)^{8}
$$
for all $(A,B)$ with $s(A,B)\in S^{o}$ and primitive isotropic vector ${\bf v}\in H^{2}(X_{(A,B)},{\bf Z})_{-}$ of level $\ell$.
Here ${\frak C}$ is the constant in \eqref{eqn:constant:Phi:generic}.
\end{lemma}

\begin{pf}
Since $S^{o}$ is connected, it suffices to show that $c_{\ell}(\cdot)$ is a locally constant function on $S^{o}$.
Let $s_{0}\in S^{o}$. There is a small neighborhood ${\mathcal U}$ of $s_{0}$ in $S^{o}$ such that 
the family of $K3$ surfaces $\pi\colon{\mathcal X}|_{\mathcal U}\to{\mathcal U}$ is topologically trivial. 
By the topological triviality, we can choose ${\bf v}$, ${\bf v}'$ and $\alpha_{{\bf v},{\bf v}'}$ to be constant under the identification of the fiber $X_{s}=\pi^{-1}(s)$ 
with $X_{s_{0}}$. Choosing ${\mathcal U}$ sufficiently small, we may assume that there exists a holomorphic section
$\sigma\colon{\mathcal U}\ni u\to(A(u),B(u))\in{\rm Sym}(3,{\bf C})\otimes{\bf C}^{6}$ such that $s(A(u),B(u))=u$ for all $u\in{\mathcal U}$.
Then both $\Phi_{\ell}\left(\alpha_{{\bf v},{\bf v}'}(z_{A(u),B(u),{\bf v},{\bf v}'})\right)^{2}$ and 
$R(A(u))R(B(u))\left(2\pi^{-2}\int_{{\bf v}^{\lor}}\omega_{(A(u),B(u))}\right)^{8}$ are holomorphic functions on ${\mathcal U}$.
As a result, $c_{\ell}(u):=c(A(u),B(u),{\bf v})$ is a holomorphic function on ${\mathcal U}$. 
Since $|c_{\ell}(u)|=1$, $c_{\ell}(u)$ must be a constant function on ${\mathcal U}$. This proves that $c_{\ell}(\cdot)$ is locally constant on $S^{o}$.
\end{pf}

\begin{lemma}
\label{lemma:Thomae:Phi:singular:case}
Let $(C,D)\in{\rm Sym}(3,{\bf C})\otimes{\bf C}^{6}$ be such that $s(C,D)\in D^{+}$.
Let $p\colon\widetilde{X}_{(C,D)}\to X_{(C,D)}$ be the minimal resolution.
Let ${\bf w}\in H^{2}_{c}(\widetilde{X}_{(C,D)},{\bf Z})_{-}$ be a primitive isotropic vector of level $\ell$ with compact support in
$\widetilde{X}_{(C,D)}\setminus p^{-1}({\rm Sing}\,X_{(C,D)})$. 
Let ${\bf w}'\in H^{2}(\widetilde{X}_{(C,D)},{\bf Z})_{-}$ be another primitive isotropic vector of level $\ell$ with $\langle{\bf w},{\bf w}'\rangle=\ell$.
Then the following equality holds
\label{cor:Thomae:Phi}
\begin{equation}
\label{eqn:Thomae:Phi:strong}
\Phi_{\ell}
\left(
\alpha_{{\bf w},{\bf w}'}(z_{C,D,{\bf w},{\bf w}'})
\right)^{2}
=
{\frak c}_{\ell}\,
R(C)R(D)
\left(
\frac{2}{\pi^{2}}
\int_{{\bf w}^{\lor}}p^{*}\omega_{(C,D)}
\right)^{8}. 
\end{equation}
\end{lemma}

\begin{pf}
Let $\gamma\colon\varDelta\ni t\to(A(t),B(t))\in{\rm Sym}(3,{\bf C})\otimes{\bf C}^{6}$ be a holomorphic map such that $s(A(t),B(t))\in S^{o}$ for $t\not=0$
and $(A(0),B(0))=(C,D)$. Let $\pi\colon{\mathcal X}\times_{S}\varDelta\to\varDelta$ be the family of $K3$ surfaces induced from
$\pi\colon{\mathcal X}\to S$ by $\gamma$. Then $\pi^{-1}(t)=X_{(A(t),B(t))}$ for $t\in\varDelta$. 
Let $\omega$ be the holomorphic section of the relative dualizing sheaf $K_{{\mathcal X}\times_{S}\varDelta/\varDelta}$
such that $\omega(t)=\omega_{(A(t),B(t))}$ for $t\in\varDelta$.
There exists a simultaneous resolution $\pi'\colon{\mathcal X}'\to\varDelta$ of the family $\pi\colon{\mathcal X}\times_{S}\varDelta\to\varDelta$ 
such that $(\pi')^{-1}(t)=X_{(A(t^{\nu}),B(t^{\nu}))}$ for $t\not=0$ and $(\pi')^{-1}(0)=\widetilde{X}_{(C,D)}$, where $\nu$ is a certain integer.
Then $\omega$ induces a nowhere vanishing holomorphic section $\omega'$ of the relative dualizing sheaf of the family
$\pi'\colon{\mathcal X}'\to\varDelta$ such that $\omega'(t)=\omega_{(A(t^{\nu}),B(t^{\nu}))}$ for $t\not=0$ and $\omega'(0)=p^{*}\omega_{(C,D)}$.
Since the family $\pi\colon{\mathcal X}'\to\varDelta$ is topologically trivial, ${\bf w}$ and ${\bf w}'$ can be regarded as $2$-cocycles of 
$(\pi')^{-1}(t)$ for all $t\in\varDelta$. Hence the functions $f(t):=\Phi_{\ell}\left(\alpha_{{\bf w},{\bf w}'}(z_{A(t),B(t),{\bf w},{\bf w}'})\right)^{2}$ and
$g(t):=R(A(t))R(B(t))\left(2\pi^{-2}\int_{{\bf w}^{\lor}}\omega_{(A(t),B(t))}\right)^{8}$ extends to holomorphic functions on $\varDelta$ such that
$f(0)=\Phi_{\ell}\left(\alpha_{{\bf w},{\bf w}'}(z_{C,D,{\bf w},{\bf w}'})\right)^{2}$ and 
$g(0)=R(C)R(D)\left(2\pi^{-2}\int_{{\bf w}^{\lor}}p^{*}\omega_{(C,D)}\right)^{8}$.
Since $f(t)/g(t)={\frak c}_{\ell}$ for $t\not=0$ by Lemma~\ref{lemma:independency:constant:holomorphic:parameter}, we get $f(0)/g(0)={\frak c}_{\ell}$.
This proves the result.
\end{pf}

\begin{proposition}
\label{proposition:Thomae:Phi:with:constant}
There exists a non-zero constant ${\frak c}_{\ell}$ with $|{\frak c}_{\ell}|={\frak C}$ such that
$$
\Phi_{\ell}\left(
\alpha_{{\bf v},{\bf v}'}(z_{A,B,{\bf v},{\bf v}'})
\right)^{2}
=
{\frak c}_{\ell}\,
R(A)R(B)
\left(
\frac{2}{\pi^{2}}
\int_{{\bf v}^{\lor}}p^{*}\omega_{(A,B)}
\right)^{8}
$$
for all admissible $(A,B)\in{\rm Sym}(3,{\bf C})\otimes{\bf C}^{6}$ and 
primitive isotropic vector ${\bf v}\in H^{2}(\widetilde{X}_{(A,B)},{\bf Z})_{-}$ of level $\ell$
with compact support in $\widetilde{X}_{(A,B)}\setminus p^{-1}({\rm Sing}\,X_{(A,B)})$.
\end{proposition}

\begin{pf}
The result follows from Lemmas~\ref{lemma:independency:constant:holomorphic:parameter} and \ref{lemma:Thomae:Phi:singular:case}.
\end{pf}

{\bf Proof of Theorem \ref{theorem:Thomae:Phi} }
In Theorem~\ref{theorem:constant:Thomae:Kummer:product:type} below, we shall prove that ${\frak c}_{\ell}=1$
in Proposition~\ref{proposition:Thomae:Phi:with:constant} for $\ell=1,2$.
The result follows from Proposition~\ref{proposition:Thomae:Phi:with:constant}.
\qed

\section
{Borcherds $\Phi$-function and theta function}
\label{sect:Phi:theta}
\par
In this section, we give a new relation between the Borcherds $\Phi$-function and certain theta functions closely related to the configuration space $X^{o}(3,6)$.

\subsection
{The configuration space $X^{o}(3,6)$}
\label{sec:configuration:space}
\par
Identify $\lambda\in({\bf C}^{*})^{6}$ with ${\rm diag}(\lambda)\in GL({\bf C}^{6})$ and
consider the $GL({\bf C}^{3})\times({\bf C}^{*})^{6}$-action on $M_{3,6}({\bf C})$ defined by $(g,\lambda)\cdot N:=gN\lambda^{2}$. 
We define
$$
X^{o}(3,6):=GL({\bf C}^{3})\backslash M^{o}_{3,6}({\bf C})/({\bf C}^{*})^{6}.
$$
By \cite[p.149]{Yoshida97}, $X^{o}(3,6)$ is a Zariski open subset of ${\bf C}^{4}$. The image of $N\in M^{o}(3,6)$ in $X^{o}(3,6)$ is denoted by $[N]$.
We consider the following group actions on $X^{o}(3,6)$:
\par
The symmetric group $S_{6}$ acts on $M(3,6)$ by $N^{\sigma}:=({\bf n}_{\sigma(1)},\ldots,{\bf n}_{\sigma(6)})$
for $N=({\bf n}_{1},\ldots,{\bf n}_{6})\in M(3,6)$ and $\sigma\in S_{6}$. This $S_{6}$-action descends to an $S_{6}$-action on $X^{o}(3,6)$.
Following \cite[Chap.\,VII Sect.\,3]{Yoshida97}, we define for $N=(N_{1},N_{2})\in M^{o}_{3,6}({\bf C})$ 
$$
N^{\lor}:=({}^{t}N_{1}^{-1},{}^{t}N_{2}^{-1})\in M^{o}_{3,6}({\bf C}).
$$
The involution $\lor\colon M^{o}(3,6)\ni N\to N^{\lor}\in M^{o}(3,6)$ on $M^{o}(3,6)$ descends to an involution $\lor\colon X^{o}(3,6)\to X^{o}(3,6)$.
By \cite[Chap.7 Prop.3.3]{Yoshida97}, the actions of $S_{6}$ and $\lor$ on $X^{o}(3,6)$ commute.
We refer to \cite[Chap.\,7-9]{Yoshida97} for more details about $X^{o}(3,6)$.

\subsection{Ten families of Enriques surfaces parametrized by $X^{o}(3,6)$}
\label{sec:expression:Phi:special}
\par
In this subsection, we use the coordinates $(x_{1},\ldots,x_{6})$ instead of $(x_{1},x_{2},x_{3},y_{1},y_{2},y_{3})$.
For $N=(n_{ij})\in M^{o}_{3,6}({\bf C})$, we set $f_{i}(x;N)= \sum_{j=1}^{6}n_{ij}x_{j}^{2}$.
Recall that
$$
X_{N}
=
\{(x_1:x_2:x_3:x_4:x_5:x_6)\in{\bf P}^{5};\,   f_1(x;N)=f_2(x;N)=f_3(x;N)=0\}.
$$
If $N=N'$ in $X^{o}(3,6)$, then $X_{N}\cong X_{N'}$.
\par
In this section, the holomorphic $2$-form on $X_{N}$ defined as the residue of $f_1(x;N)$, $f_2(x;N)$, $f_3(x;N)$
is denoted by $\omega_{N}\in H^{0}(X_{N},K_{X_{N}})$.
For a later use, let us give an explicit formula for $\omega_{N}$.
Setting $x_{3}=1$, we consider the affine coordinates $(x_{1},x_{2},x_{4},x_{5},x_{6})$ of ${\bf P}^{5}$.
Set
$$
N
=
\begin{pmatrix}
n_{11}&n_{12}&n_{13}&1&0&0
\\
n_{21}&n_{22}&n_{23}&0&1&0
\\
n_{31}&n_{32}&n_{33}&0&0&1
\end{pmatrix}
\qquad
\hbox{and}
\qquad
\begin{cases}
\begin{array}{ll}
f_{1}
&=
n_{11}x_{1}^{2}+n_{12}x_{2}^{2}+n_{13}+x_{4}^{2}
\\
f_{2}
&=
n_{21}x_{1}^{2}+n_{22}x_{2}^{2}+n_{23}+x_{5}^{2}
\\
f_{3}
&=
n_{31}x_{1}^{2}+n_{32}x_{2}^{2}+n_{33}+x_{6}^{2}.
\end{array}
\end{cases}
$$
Since 
$df_{1}\wedge df_{2}\wedge df_{3}
=
2^{3}x_{4}x_{5}x_{6}\,dx_{4}\wedge dx_{5}\wedge dx_{6}+\cdots$,
we deduce from the relation
$df_{1}\wedge df_{2}\wedge df_{3}\wedge\Upsilon=dx_{1}\wedge dx_{2}\wedge dx_{4}\wedge dx_{5}\wedge dx_{6}$ that
$$
\Upsilon
=
2^{-3}\frac{dx_{1}\wedge dx_{2}}{x_{4}x_{5}x_{6}}
=
(2i)^{-3}\frac{dx_{1}\wedge dx_{2}}
{\sqrt{\prod_{i=1}^{3}(n_{i1}x_{1}^{2}+n_{i2}x_{2}^{2}+n_{i3})}}
\mod
dx_{4},dx_{5},dx_{6}.
$$
By the definition $\omega_{N}=\Upsilon|_{X_{N}}$, we get the expression on $X_{N}\setminus\{x_{4}x_{5}x_{6}=0\}$
\begin{equation}
\label{eqn:formula:2-form:X}
\omega_{N}
=
\frac{dx_{1}\wedge dx_{2}}{2^{3}x_{4}x_{5}x_{6}}
=
(2i)^{-3}\frac
{dx_{1}\wedge dx_{2}}{\sqrt{\prod_{i=1}^{3}(n_{i1}x_{1}^{2}+n_{i2}x_{2}^{2}+n_{i3})}}.
\end{equation}

\par
For $J=\{j_{1},j_{2},j_{3}\}\subset\{1,\ldots,6\}$ with $j_{1}<j_{2}<j_{3}$, let $\langle J\rangle$ denote the partition $\{1,\ldots,6\}=J\amalg J^{c}$. 
Write $J^{c}=\{j_{4},j_{5},j_{6}\}$, $j_{4}<j_{5}<j_{6}$. Then $\langle J\rangle$ is also written as $\binom{j_{1}j_{2}j_{3}}{j_{4}j_{5}j_{6}}$. 
Hence $\langle J\rangle=\binom{j_{1}j_{2}j_{3}}{j_{4}j_{5}j_{6}}$. 
For a partition $\langle J\rangle=\binom{j_{1}j_{2}j_{3}}{j_{4}j_{5}j_{6}}$ and $N\in M(3,6)$, 
we set $\Delta_{\langle J\rangle}(N):=\Delta_{j_{1}j_{2}j_{3}}(N)\Delta_{j_{4}j_{5}j_{6}}(N)$. 
We define an involution $\iota_{\langle J\rangle}$ on ${\bf P}^{5}$ by
$$
\iota_{\langle J\rangle}(x_{j_{1}},x_{j_{2}},x_{j_{3}},x_{j_{4}},x_{j_{5}},x_{j_{6}})
:=
(x_{j_{1}},x_{j_{2}},x_{j_{3}},-x_{j_{4}},-x_{j_{5}},-x_{j_{6}}).
$$
Then $\iota_{\langle J\rangle}$ acts on $X_{N}$. 
For $N\in M_{3,6}^{o}({\bf C})$, we define an Enriques surface $Y_{N,\langle J\rangle}$ by
$$
Y_{N,\langle J\rangle}:=X_{N}/\iota_{\langle J\rangle}.
$$
Then $Y_{N,\langle J\rangle}\cong Y_{N',\langle J\rangle}$ if $[N]=[N']$. 
We define a map $\overline{\varpi}_{\langle J\rangle}\colon X^{o}(3,6)\to{\mathcal M}^{o}$ by
$\overline{\varpi}_{\langle J\rangle}([N]):=\overline{\varpi}(Y_{N,\langle J\rangle})$. 
For all $[N]\in X^{o}(3,6)$, we have
$$
(\overline{\varpi}_{\langle J\rangle}^{*}\|\Phi\|)([N])=\left\|\Phi(Y_{N,\langle J\rangle})\right\|.
$$

\begin{theorem}
\label{thm:valuePhi:special}
For all $N\in M^{o}_{3,6}({\bf C})$ and for all partitions $\langle J\rangle$,
$$
\left\|\Phi(Y_{N,\langle J\rangle})\right\|
=
|\Delta_{\langle J\rangle}(N)|^{2}\cdot
\left(
\frac{2}{\pi^{4}}
\left|\int_{X_{N}}\omega_{N}\wedge\overline{\omega}_{N}\right|
\right)^{2}.
$$
\end{theorem}

\begin{pf}
The symmetric group $S_{6}$ acts on ${\bf P}^{5}$ by $\sigma(x_{i})=x_{\sigma(i)}$, $(1\leq i\leq 6)$.
There is $\sigma\in S_{6}$ with $\sigma^{-1}\iota\sigma=\iota_{\langle J\rangle}$.
Since $\sigma(X_{N})=X_{N^{\sigma}}$ and $\sigma^{*}\omega_{N^{\sigma}}=\omega_{N}$, we get
$$
|\Delta_{\binom{123}{456}}(N^{\sigma})|^{2}
\left(
\frac{2}{\pi^{4}}
\left|
\int_{X_{N^{\sigma}}}\omega_{N^{\sigma}}\wedge\overline{\omega}_{N^{\sigma}}
\right|
\right)^{2}
=
|\Delta_{\langle J\rangle}(N)|^{2}
\left(
\frac{2}{\pi^{4}}
\left|
\int_{X_{N}}\omega_{N}\wedge\overline{\omega}_{N}
\right|
\right)^{2}.
$$
Since $R(Q)=\det(A)^{4}$ for $Q=(Q_{1},Q_{2},Q_{3})$, $Q_{i}(x)=\sum_{j=1}^{3}a_{ij}x_{j}^{2}$ and 
since $Y_{N,\langle J\rangle}=Y_{N^{\sigma}}$,
the left hand side is equal to $\left\|\Phi(Y_{N^{\sigma}})\right\|$ by Theorem~\ref{thm:valuePhi:general}.
\end{pf}

\subsection
{Borcherds $\Phi$-function and theta function}
\label{subsect:Phi:Freitag:theta}
\par
We consider another families of $K3$ surfaces over $X^{o}(3,6)$ (cf. \cite{Yoshida97}, \cite{MatsumotoTerasoma12}).
For $N=(n_{ij})\in M^{o}_{3,6}({\bf C})$, define
$$
Z_{N}
:=
\{((z_{1}:z_{2}:z_{3}),w)\in{\mathcal O}_{{\bf P}^{2}}(3);\,
w^{2}=\prod_{i=1}^{6}(n_{1i}z_{1}+n_{2i}z_{2}+n_{3i}z_{3})\},
$$
where $w$ denotes the coordinate of the fibers of ${\mathcal O}_{{\bf P}^{2}}(3)$. 
When $N\in M^{o}_{3,6}({\bf C})$, $Z_{N}$ is the singular $K3$ surface with $15$ nodes defined as the double covering of ${\bf P}^{2}$, 
whose branch divisor is the union of $6$ lines $n_{1i}z_{1}+n_{2i}z_{2}+n_{3i}z_{3}=0$ $(i=1,\ldots,6)$.
Let $\pi\colon W_{N}\to Z_{N}$ be the minimal resolution of $Z_{N}$.
Then $W_{N}$ is a smooth $K3$ surface. 
By construction, $W_{N}\cong W_{N'}$ if $[N]=[N']$ in $X^{o}(3,6)$.
\par
We define a holomorphic section of the dualizing sheaf of $Z_{N}$ by
$$
\eta_{N}
:=
\frac
{z_{1}dz_{2}\wedge dz_{3}-z_{2}dz_{1}\wedge dz_{3}+z_{3}dz_{1}\wedge dz_{2}}{w}.
$$
On the open subset $Z_{N}\cap\{z_{3}\not=0\}$, we get the expression
\begin{equation}
\label{eqn:formula:eta}
\eta_{N}
=
\left.
\frac{dy_{1}\wedge dy_{2}}{\sqrt{\prod_{i=1}^{6}(n_{1i}y_{1}+n_{2i}y_{2}+n_{3i})}}
\right|_{Z_{N}\cap\{z_{3}\not=0\}},
\end{equation}
where $(y_{1},y_{2})$ is the inhomogeneous coordinates of ${\bf P}^{2}$.
\par
Let $H_{2}(W_{N},{\bf Z})_{-}$ be the anti-invariant part of $H_{2}(W_{N},{\bf Z})$ with respect to the involution induced by the one 
$(z,w)\mapsto(z,-w)$ on $Z_{N}$. 
In \cite[Sect.\,3.1]{MatsumotoTerasoma12}, Matsumoto-Terasoma constructed a basis
$\{\gamma'_{12},\gamma'_{13},\gamma'_{14},\gamma'_{23},\gamma'_{24},\gamma'_{34}\}$ of $H^{2}(W_{N},{\bf Z})_{-}$. 
(See also Sect.~\ref{subsect:Periods:Jacobian:Kummer} below.)
Let $2H$ be the Gram matrix of the intersection form on $H^{2}(W_{N},{\bf Z})_{-}$ with respect to 
$\{\gamma'_{12},\gamma'_{13},\gamma'_{14},\gamma'_{23},\gamma'_{24},\gamma'_{34}\}$. 
By \cite[Prop.\,3]{MatsumotoTerasoma12}, we have
$$
\left(
\int_{W_{N}}\gamma'_{ij}\wedge\gamma'_{kl}
\right)_{1\leq i<j\leq 4,\,1\leq k<l\leq 4}
=
2H,
\qquad
H:=
-
\begin{pmatrix}
0&0&0&0&0&1
\\
0&1&0&0&0&0
\\
0&0&0&1&0&0
\\
0&0&1&0&0&0
\\
0&0&0&0&1&0
\\
1&0&0&0&0&0
\end{pmatrix}.
$$
We set $\eta_{ij}(N):=\langle\gamma'_{ij},\eta_{N}\rangle/2$ for $1\leq i<j\leq 4$ and
$$
\eta(N)
:=
(\eta_{12}(N),\eta_{13}(N),\eta_{14}(N),\eta_{23}(N),\eta_{24}(N),\eta_{34}(N))\in{\bf C}^{6}.
$$
Since 
$-\eta_{N}=\eta_{34}(N)\gamma'_{12}+\eta_{13}(N)\gamma'_{13}+\eta_{23}(N)\gamma'_{14}
+\eta_{14}\gamma'_{23}+\eta_{24}(N)\gamma'_{24}+\eta_{12}(N)\gamma'_{34}$,
it follows from the Riemann-Hodeg bilinear relations that
\begin{equation}
\label{eqn:Riemann:Hodge}
\eta(N)H{}^{t}\eta(N)=0,
\qquad
\int_{W_{N}}\eta_{N}\wedge\overline{\eta}_{N}
=
\eta(N)\cdot 2H\cdot{}^{t}\overline{\eta(N)}.
\end{equation}
\par
Following \cite[Eq.(15)]{MatsumotoTerasoma12}, we define the normalized period of $W_{N}$ by
$$
\varOmega_{N}
:=
\frac{1}{\eta_{34}(N)}
\begin{pmatrix}
\eta_{14}(N)&-\frac{\eta_{13}(N)-i\eta_{24}(N)}{1+i}
\\
-\frac{\eta_{13}(N)+i\eta_{24}(N)}{1-i}&-\eta_{23}(N)
\end{pmatrix}
\in{\Bbb D},
$$
where
$$
{\Bbb D}:=\{\tau\in M_{2,2}({\bf C});\,(\tau-{}^{t}\overline{\tau})/2\sqrt{-1}\hbox{ is positive-definite}\}
$$ 
is a symmetric bounded domain of type $I_{2,2}$. 
Notice that ${\Bbb D}$ is isomorphic to a symmetric bounded domain of type IV of dimension $4$ (cf. \cite[Sect.\,4.1]{MatsumotoTerasoma12}).
\par
Let $U_{22}^{M}(1+i)$ be the discrete group acting on ${\Bbb D}$ defined in \cite[Prop.\,6]{MatsumotoTerasoma12}.
For $\varOmega\in{\Bbb D}$, write $[\varOmega]\in{\Bbb D}/U_{22}^{M}(1+i)$ for the $U_{22}^{M}(1+i)$-orbit of $\varOmega$.
By \cite[Prop.\,6]{MatsumotoTerasoma12}, the period map for the family of $K3$ surfaces $\bigcup_{N\in M^{o}(3,6)}W_{N}\to M^{o}(3,6)$ 
induces a holomorphic map
$$
{\mathcal P}\colon X^{o}(3,6)\ni[N]\to[\varOmega_{N}]\in{\Bbb D}/U_{22}^{M}(1+i).
$$
\par
Set ${\bf e}(x):=\exp(2\pi ix)$. Following \cite[p.137]{MatsumotoTerasoma12}, we set
$$
Ev:=\{({\bf a},{\bf b})\in({\bf Z}[i]/(1+i){\bf Z}[i])^{4};\,{\bf a}{}^{t}\bar{\bf b}\equiv0\mod(1+i)\}.
$$
For $({\bf a},{\bf b})\in Ev$, the theta function $\Theta_{[{\bf a},{\bf b}]}(\varOmega)\in{\mathcal O}({\Bbb D})$ is defined as 
(cf. \cite[Sect.4.4]{MatsumotoTerasoma12})
$$
\Theta_{[{\bf a},{\bf b}]}(\varOmega)
:=
\sum_{n\in{\bf Z}[i]^{2}}
{\bf e}
\left[
\left(n+\frac{\bf a}{1+i}\right)
\frac{\varOmega}{2}
{}^{t}
\overline{\left(n+\frac{\bf a}{1+i}\right)}
+\Re\left(n+\frac{\bf a}{1+i}\right){}^{t}\overline{\left(\frac{\bf b}{1+i}\right)}
\right].
$$
By Matsumoto-Terasoma \cite[p.138]{MatsumotoTerasoma12},
there is a one-to-one correspondence between the set of partitions $\{\langle J\rangle\}$ and 
$Ev=\{\binom{\bf a}{\bf b}=\binom{{\bf a}_{1}\,{\bf a}_{2}}{{\bf b}_{1}\,{\bf b}_{2}}\}$ as follows:
$$
\begin{array}{cccccccccccc}
\binom{\bf a}{\bf b}
&
\vline
&
\binom{0\,0}{0\,0}
&
\binom{i\,0}{0\,0}
&
\binom{0\,i}{0\,0}
&
\binom{0\,0}{i\,0}
&
\binom{0\,0}{0\,i}
&
\binom{i\,i}{0\,0}
&
\binom{i\,0}{0\,i}
&
\binom{0\,0}{i\,i}
&
\binom{0\,i}{i\,0}
&
\binom{i\,i}{i\,i}
\\
\hline
\langle J\rangle
&
\vline
&
\binom{135}{246}
&
\binom{146}{235}
&
\binom{136}{245}
&
\binom{125}{346}
&
\binom{134}{256}
&
\binom{145}{236}
&
\binom{156}{234}
&
\binom{124}{356}
&
\binom{126}{345}
&
\binom{123}{456}
\end{array}
$$
When a partition $\langle J\rangle$ corresponds to $({\bf a},{\bf b})\in Ev$ by the above rule, we define
$$
\Theta_{\langle J\rangle}(\varOmega)
:=
\Theta_{[{\bf a},{\bf b}]}(\varOmega).
$$
\par
The Petersson norm of $\Theta_{[{\bf a},{\bf b}]}(\varOmega)$ is the function on ${\Bbb D}$ defined as
$$
\left\|
\Theta_{[{\bf a},{\bf b}]}(\varOmega)
\right\|^{2}
:=
\det\left(\frac{\varOmega-{}^{t}\overline{\varOmega}}{2i}\right)
\left|
\Theta_{[{\bf a},{\bf b}]}(\varOmega)
\right|^{2}.
$$
Since $\|\Theta_{[{\bf a},{\bf b}]}\|$ is an $U_{22}^{M}(1+i)$-invariant function on ${\Bbb D}$ by \cite[Prop.\,7]{MatsumotoTerasoma12},
we regard $\|\Theta_{[{\bf a},{\bf b}]}\|$ as a function on ${\Bbb D}/U_{22}^{M}(1+i)$. 
Then ${\mathcal P}^{*}\|\Theta_{[{\bf a},{\bf b}]}\|$ is a function on $X^{o}(3,6)$.

\begin{theorem}
\label{thm:Phi:Freitag:theta}
For all $N\in M^{o}_{3,6}({\bf C})$ and for all partitions $\langle J\rangle$, 
$$
\left\|\Phi\left(Y_{N,\langle J\rangle}\right)\right\|
=
\left\|
\Theta_{\langle J\rangle}\left(\varOmega_{N}\right)
\right\|^{4}.
$$
Namely, one has the equality of functions $\overline{\varpi}_{\langle J\rangle}^{*}\|\Phi\|={\mathcal P}^{*}\left\|\Theta_{\langle J\rangle}\right\|^{4}$
on $X^{o}(3,6)$.
\end{theorem}

\begin{pf}
{\em (Step 1) }
Firstly, we set $\langle J\rangle=\binom{123}{456}$ and prove the assertion in this case.
For $K=(k_{ij})\in M_{3,3}({\bf C})$ with $(K,I)\in M^{o}_{3,6}({\bf C})$, set 
$$
V_{K}
:=
\left\{
(y_{1},y_{2},y_{4},y_{5},y_{6})\in{\bf C}^{5};\,
k_{i1}y_{1}+k_{i2}y_{2}+k_{i3}+y_{i+3}=0
\quad
(i=1,2,3)
\right\}.
$$
Then $V_{K}$ is an affine subspace of ${\bf C}^{5}$, whose closure in ${\bf P}^{5}$ is denoted by $\overline{V}_{K}$. 
Set
$$
\xi_{K}
=
\left.
2^{-5}
\frac{dy_{1}\wedge dy_{2}}
{\sqrt{y_{1}y_{2}\prod_{i=1,2,3}(k_{i1}y_{1}+k_{i2}y_{2}+k_{i3})}}
\right|_{V_{K}}.
$$
Let $\varphi\colon{\bf P}^{5}\to{\bf P}^{5}$ be the map defined as $\varphi(x_{1}:\ldots:x_{6})=(x_{1}^{2}:\ldots:x_{6}^{2})$. 
Since $\overline{V}_{K}=\varphi(X_{(K,I)})$ and $\omega_{(K,I)}=i^{-3}\varphi^{*}\xi_{K}$, 
we get by Theorem~\ref{thm:valuePhi:special} and \eqref{eqn:formula:2-form:X}
\begin{equation}
\label{eqn:Phi:theta:1}
\begin{aligned}
\,&
\left\|\Phi(Y_{(K,I)})\right\|^{2}
=
|\Delta_{123}(K,I)|^{4}\cdot
\left(
\frac{2}{\pi^{4}}
\left|
\int_{X_{(K,I)}}\varphi^{*}(\xi_{K}\wedge\overline{\xi}_{K})
\right|
\right)^{4}
\\
&=
|\det(K)|^{4}
\cdot
\left(
\frac{2\deg(\varphi)}{\pi^{4}}
\left|\int_{V_{K}}\xi_{K}\wedge\overline{\xi}_{K}\right|
\right)^{4}
\\
&=
\left(
\frac{2^{-9}\deg(\varphi)|\det(K)|}{\pi^{4}}
\int_{\overline{V}_{K}}
\frac{
(\sqrt{-1})^{2}dy_{1}\wedge d\bar{y}_{1}\wedge dy_{2}\wedge d\bar{y}_{2}}
{|y_{1}y_{2}\prod_{i=1}^{3}(k_{i1}y_{1}+k_{i2}y_{2}+k_{i3})|}
\right)^{4}
\\
&=
\left(
\frac{2^{-4}|\det(K)|}{\pi^{4}}
\int_{{\bf P}^{2}}
\eta_{({}^{t}K,I)}\wedge\overline{\eta}_{({}^{t}K,I)}
\right)^{4}.
\end{aligned}
\end{equation}
Here $\eta_{({}^{t}K,I)}\wedge\overline{\eta}_{({}^{t}K,I)}$ is regarded as a volume form on ${\bf P}^{2}$ and
the last line follows from \eqref{eqn:formula:eta} and the formula $\deg(\varphi)=2^{5}$.
\par
For $g\in GL({\bf C}^{3})$, let $f_{g}\in PGL({\bf C}^{3})$ be the projective transform $f_{g}[z]=[{}^{t}g^{-1}z]$.
Regarding $\eta_{N}\wedge\overline{\eta}_{N}$ as a volume form on ${\bf P}^{2}$ for $N\in M^{o}_{3,6}({\bf C})$, we get
\begin{equation}
\label{eqn:eta:transformation}
f_{g}^{*}(\eta_{gN}\wedge\overline{\eta}_{gN})
=
|\det(g)|^{-2}\eta_{N}\wedge\overline{\eta}_{N}.
\end{equation}
Since $X_{gN}=X_{N}$ and thus $\|\Phi(Y_{gN})\|=\|\Phi(Y_{N})\|$ for all $g\in GL({\bf C}^{3})$, $N\in M^{o}_{3,6}({\bf C})$,
we deduce from \eqref{eqn:Phi:theta:1}, \eqref{eqn:eta:transformation} that
for any $N=(N_{1},N_{2})=N_{2}(N_{2}^{-1}N_{1},I)\in M^{o}_{3,6}({\bf C})$
\begin{equation}
\label{eqn:Phi:theta:2}
\begin{aligned}
\left\|\Phi(Y_{N})\right\|^{2}
&=
\left(
\frac{2^{-4}|\det(N_{2}^{-1}N_{1})|}{\pi^{4}}
\int_{{\bf P}^{2}}
\eta_{({}^{t}N_{1}{}^{t}N_{2}^{-1},I)}\wedge\overline{\eta}_{({}^{t}N_{1}{}^{t}N_{2}^{-1},I)}
\right)^{4}
\\
&=
\left(
\frac{2^{-4}|\det({}^{t}N_{1}^{-1})\det({}^{t}N_{2}^{-1})|}{\pi^{4}}
\int_{{\bf P}^{2}}
\eta_{({}^{t}N_{2}^{-1},{}^{t}N_{1}^{-1})}\wedge\overline{\eta}_{({}^{t}N_{2}^{-1},{}^{t}N_{1}^{-1})}
\right)^{4}
\\
&=
\left(
\frac{2^{-4}|\Delta_{\langle J\rangle}(N^{\lor})|}{\pi^{4}}\cdot
\frac{1}{2}
\int_{Z_{N^{\lor}}}\eta_{N^{\lor}}\wedge\overline{\eta}_{N^{\lor}}
\right)^{4}.
\end{aligned}
\end{equation}
Here $\eta_{N^{\lor}}$ is regarded as a canonical form on $Z_{N^{\lor}}$ in the last line.
By \eqref{eqn:Phi:theta:2}, we get
\begin{equation}
\label{eqn:Phi:theta:3}
\begin{aligned}
\,&
\|\Phi(Y_{N})\|^{2}
=
\left\{
\frac{2^{-3}|\Delta_{\langle J\rangle}(N^{\lor})\eta_{34}(N^{\lor})^{2}|}{4\pi^{4}}\,
\int_{Z_{N^{\lor}}}\frac{\eta_{N^{\lor}}}{\eta_{34}(N^{\lor})}\wedge\overline{\frac{{\eta}_{N^{\lor}}}{\eta_{34}(N^{\lor})}}
\right\}^{4}
\\
&=
\left\{
2^{-3}
\left|
\Theta_{\langle J\rangle}\left(\varOmega_{N^{\lor}}\right)
\right|^{2}
8\det\left(
\frac{\varOmega_{N^{\lor}}-{}^{t}\overline{\varOmega}_{N^{\lor}}}{2i}
\right)
\right\}^{4}
=
\left\|
\Theta_{\langle J\rangle}\left(\varOmega_{N^{\lor}}\right)
\right\|^{8}.
\end{aligned}
\end{equation}
Here the second equality follows from \cite[Eq.(4)]{MatsumotoTerasoma12} and the following identity
\begin{equation}
\label{eqn:Bergman:kernal:D}
\begin{aligned}
\,&
\det\left(\frac{\varOmega_{N^{\lor}}-{}^{t}\overline{\varOmega}_{N^{\lor}}}{2i}\right)
=
-(\Im\frac{\eta_{14}}{\eta_{34}})\cdot(\Im\frac{\eta_{23}}{\eta_{34}})
-\frac{1}{2}
\left\{
(\Im\frac{\eta_{13}}{\eta_{34}})^{2}+(\Im\frac{\eta_{24}}{\eta_{34}})^{2}
\right\}
=
\\
&
\frac{1}{2}\Im\left(\frac{\eta}{\eta_{34}}\right)
\cdot H\cdot
{}^{t}\left\{
\Im\left(\frac{\eta}{\eta_{34}}\right)
\right\}
=
\frac{1}{8}\left(\frac{\eta}{\eta_{34}}\right)
\cdot 2H\cdot
{}^{t}\overline{\left(\frac{\eta}{\eta_{34}}\right)}
=
\frac{1}{8}\int_{Z_{N^{\lor}}}\frac{\eta_{N^{\lor}}}{\eta_{34}}\wedge\overline{\frac{{\eta}_{N^{\lor}}}{\eta_{34}}},
\end{aligned}
\end{equation}
where we wrote $\eta_{ij}=\eta_{ij}(N^{\lor})$ and $\eta=\eta(N^{\lor})$ and used \eqref{eqn:Riemann:Hodge} to get the last equality. 
Since $[\varOmega_{N^{\lor}}]=[{}^{t}\varOmega_{N}]$ by the equality $\psi(*x)=T\cdot\psi(x)$ in \cite[p.260]{Yoshida97}
and since $\Theta_{[{\bf a},{\bf b}]}(\varOmega)^{2}=\Theta_{[{\bf a},{\bf b}]}({}^{t}\varOmega)^{2}$ by \cite[Lemma 3.1.3]{Matsumoto93}, 
we deduce from \eqref{eqn:Phi:theta:3} that
\begin{equation}
\label{eqn:Phi:theta:4}
\|\Phi(Y_{N})\|
=
\left\|\Theta_{\langle J\rangle}(\varOmega_{N^{\lor}})\right\|^{4}
=
\left\|\Theta_{\langle J\rangle}({}^{t}\varOmega_{N})\right\|^{4}
=
\left\|\Theta_{\langle J\rangle}(\varOmega_{N})\right\|^{4}.
\end{equation}
\par{\em (Step 2) }
Let $\langle J\rangle$ be an arbitrary partition and let $\sigma\in S_{6}$ be a permutation such that
$\langle J\rangle=\binom{\sigma(1)\sigma(2)\sigma(3)}{\sigma(4)\sigma(5)\sigma(6)}$.
By \cite[Th.\,2]{MatsumotoTerasoma12} and \eqref{eqn:Bergman:kernal:D}, we get for all $N\in M^{o}(3,6)$,
\begin{equation}
\label{eqn:Phi:theta:5}
\begin{aligned}
\left\|\Theta_{\langle J\rangle}(\varOmega_{N})^{2}\right\|^{2}
&=
\left|
\frac{1}{4\pi^{2}}\Delta_{\langle J\rangle}(N)\eta_{34}(N)^{2}
\right|^{2}
\left(
\frac{1}{8}\int_{Z_{N}}\frac{\eta_{N}}{\eta_{34}(N)}\wedge\overline{\frac{\eta_{N}}{\eta_{34}(N)}}
\right)^{2}
\\
&=
\left|
\frac{1}{4\pi^{2}}\Delta_{\langle J\rangle}(N)
\right|^{2}
\left(
\frac{1}{8}\int_{Z_{N}}\eta_{N}\wedge\overline{\eta}_{N}
\right)^{2}
\\
&=
\left|
\frac{1}{4\pi^{2}}\Delta_{\binom{123}{456}}(N^{\sigma})\eta_{34}(N^{\sigma})^{2}
\right|^{2}
\left(
\frac{1}{8}\int_{Z_{N^{\sigma}}}\frac{\eta_{N^{\sigma}}}{\eta_{34}(N^{\sigma})}\wedge\overline{\frac{{\eta}_{N^{\sigma}}}{\eta_{34}(N^{\sigma})}}
\right)^{2}
\\
&=
\left\|\Theta_{\binom{123}{456}}(\varOmega_{N^{\sigma}})^{2}\right\|^{2}.
\end{aligned}
\end{equation}
Since $Y_{N,\langle J\rangle}=Y_{N^{\sigma}}$, the result follows from \eqref{eqn:Phi:theta:5} and \eqref{eqn:Phi:theta:4} applied to $N^{\sigma}$.
\end{pf}

\subsection
{The case where $X_{N}$ is a Jacobian Kummer surface}
\label{subsect:JacobanKummerCase}
\par
Let $\lambda:=(\lambda_{k})\in{\bf C}^{6}$ be such that $\lambda_{k}\not=\lambda_{l}$ for all $k\not=l$.
Then
$$
C_{\lambda}
:=
\{
(x,y)\in{\bf C}^{2};\,
y^{2}=(x-\lambda_{1})(x-\lambda_{2})(x-\lambda_{3})(x-\lambda_{4})(x-\lambda_{5})(x-\lambda_{6})
\}
$$
is a curve of genus $2$ with the {\em ordered} set of branch points $\lambda=(\lambda_{k})$. 
We identify $C_{\lambda}$ with the corresponding projective curve. Then $H^{0}(C_{\lambda},\Omega_{C_{\lambda}}^{1})$ is equipped with
the basis $\{\omega_{1}=dx/y,\omega_{2}=xdx/y\}$. 
Let $\{A_{1},A_{2},B_{1},B_{2}\}$ be the canonical basis of $H_{1}(C_{\lambda},{\bf Z})$ as in \cite[Sect.\,3.1]{MatsumotoTerasoma12}.
Since $\{\lambda_{k}\}$ is an ordered set, $C_{\lambda}$ is equipped with a level $2$-structure.
The period of $C_{\lambda}$ is defined by
$$
[T_{\lambda}]
:=
\left[
\begin{pmatrix}
\int_{B_{1}}\omega_{1}&\int_{B_{2}}\omega_{1}
\\
\int_{B_{1}}\omega_{2}&\int_{B_{2}}\omega_{2}
\end{pmatrix}^{-1}
\begin{pmatrix}
\int_{A_{1}}\omega_{1}&\int_{A_{2}}\omega_{1}
\\
\int_{A_{1}}\omega_{2}&\int_{A_{2}}\omega_{2}
\end{pmatrix}
\right]
\in{\frak S}_{2}/\Gamma(2),
$$
where ${\frak S}_{2}$ is the Siegel upper half-space of degree $2$ and $\Gamma(2)\subset Sp_{4}({\bf Z})$ is the principal congruence subgroup
of level $2$.
\par
Let $K(C_{\lambda})$ be the Kummer surface associated to the Jacobian variety of $C_{\lambda}$. 
There are two models of $K(C_{\lambda})$.
For each $\lambda_{i}$, the point $(1:\lambda_{i}:\lambda_{i}^{2})$ lies on the conic $x_{0}x_{2}-x_{1}^{2}=0$
with tangent line $\ell_{i}=\{\lambda_{i}^{2}x_{0}-2\lambda_{i}x_{1}+x_{2}=0\}$. 
We set
$$
N_{\lambda}
:=
\begin{pmatrix}
1&1&1&1&1&1
\\
\lambda_{1}&\lambda_{2}&\lambda_{3}&\lambda_{4}
&\lambda_{5}&\lambda_{6}
\\
\lambda_{1}^{2}&\lambda_{2}^{2}&\lambda_{3}^{2}&\lambda_{4}^{2}
&\lambda_{5}^{2}&\lambda_{6}^{2}
\end{pmatrix},
\qquad
N'_{\lambda}
:=
\begin{pmatrix}
1&0&0
\\
0&-2&0
\\
0&0&1
\end{pmatrix}N_{\lambda}.
$$
Then $K(C_{\lambda})$ is isomorphic to the minimal resolution of the double covering of ${\bf P}^{2}$ with branch divisor $\ell_{1}\cup\cdots\cup\ell_{6}$. 
Hence
$$
K(C_{\lambda})\cong W_{N'_{\lambda}}\cong W_{N_{\lambda}}.
$$
On the other hand, it follows from \cite[p.769--770, p.789]{GriffithsHarris78} that
$$
K(C_{\lambda})\cong X_{N_{\lambda}}.
$$
By the second expression of $K(C_{\lambda})$, we have ten free involutions $\{\iota_{\langle J\rangle}\}$ on $K(C_{\lambda})$.
\par
Recall that, for $a,b\in\{0,\frac{1}{2}\}^{2}$, the Riemann theta constant $\theta_{a,b}(T)$ is defined as
$$
\theta_{a,b}(T)
:=
\sum_{n\in{\bf Z}^{2}}
{\bf e}
\left[
\left(n+a\right)(T/2){}^{t}\left(n+a\right)+{}^{t}\left(n+a\right)b
\right],
\qquad
T\in{\frak S}_{2},
$$
whose Petersson norm is defined as $\|\theta_{a,b}(T)\|^{2}:=(\det\Im T)^{\frac{1}{2}}|\theta_{a,b}(T)|^{2}$ and whose parity is defined as $4{}^{t}ab\mod2$.
By \cite[Lemma 2.1.1 (vi) and p.399 l.2-4]{Matsumoto93}, we get
\begin{equation}
\label{eqn:theta:Freitag:Riemann}
\Theta_{[{\bf a},{\bf b}]}(T)
=
\theta_{\Re(\frac{\bf a}{1+i}),\Re(\frac{\bf b}{1+i})}(T)^{2}
=
\theta_{\Im(\frac{\bf a}{1+i}),\Im(\frac{\bf b}{1+i})}(T)^{2}.
\end{equation}

\begin{theorem}
\label{thm:Phi:theta:Kummer} 
Let $\langle J\rangle$ be a partition corresponding to $({\bf a},{\bf b})\in Ev$. Then
$$
\left\|
\Phi\left(K(C_{\lambda})/\iota_{\langle J\rangle}\right)
\right\|
=
\left\|
\theta_{\Re(\frac{\bf a}{1+i})\Re(\frac{\bf b}{1+i})}\left(T_{\lambda}\right)
\right\|^{8}
=
\left\|
\theta_{\Im(\frac{\bf a}{1+i})\Im(\frac{\bf b}{1+i})}\left(T_{\lambda}\right)
\right\|^{8}.
$$
\end{theorem}

\begin{pf}
Since $[\varOmega_{N_{\lambda}}]=[T_{\lambda}]$ by \cite[Chap.\,IX, Remark 10.2]{Yoshida97}, 
the result follows from Theorem~\ref{thm:Phi:Freitag:theta} and \eqref{eqn:theta:Freitag:Riemann}.
\end{pf}

\section
{An infinite product expansion of theta constants of genus two}
\label{sect:theta:constant:genus2}
\par
\subsection
{Periods of principally polarized Abelian surfaces}
\label{subsect:Periods:ppas}
\par
Let ${\bf e}_{1}=\binom{1}{0},{\bf e}_{2}=\binom{0}{1}\in{\bf C}^{2}$.
For $T=({\bf t}_{1},{\bf t}_{2})\in{\frak S}_{2}$ with ${\bf t}_{1}=\binom{T_{11}}{T_{21}},{\bf t}_{2}=\binom{T_{12}}{T_{22}}\in{\bf C}^{2}$, 
we define 
$$
A_{T}:={\bf C}^{2}/{\bf Z}{\bf t}_{1}+{\bf Z}{\bf t}_{2}+{\bf Z}{\bf e}_{1}+{\bf Z}{\bf e}_{2}.
$$
Then $A_{T}$ is an Abelian surface with period matrix $(T,I_{2})$.
Let $\theta(z,T)\in{\mathcal O}({\bf C}^{2}\times{\frak S}_{2})$ be the Riemann theta function and let $C_{T}$ be the theta divisor of $A_{T}$:
$$
C_{T}:=\{[z]\in A_{T};\,\theta(z,T)=0\}
$$ 
\par
Let ${\frak S}_{2}^{o}$ be the complement of the $Sp_{4}({\bf Z})$-orbit of the diagonal locus in ${\frak S}_{2}$.
Let $T\in{\frak S}_{2}^{o}$. Then $C_{T}$ is a curve of genus $2$ with Jacobian variety $A_{T}$. 
Let $i_{T}\colon C_{T}\hookrightarrow A_{T}$ be the inclusion, which  induces an isomorphism $H^{1}(A_{T},{\bf Z})\cong H^{1}(C_{T},{\bf Z})$.
Let $\{\alpha_{1},\alpha_{2},\beta_{1},\beta_{2}\}$ be the canonical basis of $H_{1}(C_{T},{\bf Z})$ (cf. \cite[pp.227-228]{GriffithsHarris78}) such that
$(i_{T})_{*}\alpha_{i}$, $(i_{T})_{*}\beta_{i}$ correspond to the cycles $\vec{\bf e}_{i},\vec{\bf t}_{i}\in H_{1}(A_{T},{\bf Z})$. Then
$$
\int_{(i_{T})_{*}\alpha_{i}}dz_{j}=T_{ij},
\qquad
\int_{(i_{T})_{*}\beta_{i}}dz_{j}=\delta_{ij}.
$$
Let $\alpha_{1}^{\lor},\alpha_{2}^{\lor},\beta_{1}^{\lor},\beta_{2}^{\lor}\in H^{1}(C_{T},{\bf Z})$ be the Poincar\'e duals of 
$\alpha_{1},\alpha_{2},\beta_{1},\beta_{2}$, respectively. 
Then $\{\alpha_{1}^{\lor},\alpha_{2}^{\lor},\beta_{1}^{\lor},\beta_{2}^{\lor}\}$ is a symplectic basis of $H^{1}(C_{T},{\bf Z})$.
\par
Let $\{{\bf a}_{1}^{\lor},{\bf a}_{2}^{\lor},{\bf b}_{1}^{\lor},{\bf b}_{2}^{\lor}\}\subset H^{1}(A_{T},{\bf Z})$ be the basis such that
$i_{T}^{*}{\bf a}_{i}^{\lor}=\alpha_{i}^{\lor}$, $i_{T}^{*}{\bf b}_{i}^{\lor}=\beta_{i}^{\lor}$.
By  \cite[p.310 Lemma]{GriffithsHarris78}, 
$[C_{T}]={\bf a}_{1}^{\lor}\wedge{\bf b}_{1}^{\lor}+{\bf a}_{2}^{\lor}\wedge{\bf b}_{2}^{\lor}\in H^{2}(A_{T},{\bf Z})$ 
is the Poincar\'e dual of $C_{T}\subset A_{T}$. Set
$$
{\Bbb A}:=[C_{T}]^{\perp}=\{\alpha\in H^{2}(A_{T},{\bf Z});\,\langle[C_{T}],\alpha\rangle=0\}\subset H^{2}(A_{T},{\bf Z}).
$$
Equipped with the basis
$\{{\bf a}_{1}^{\lor}\wedge{\bf a}_{2}^{\lor},{\bf b}_{1}^{\lor}\wedge{\bf b}_{2}^{\lor},
{\bf a}_{1}^{\lor}\wedge{\bf b}_{2}^{\lor},{\bf a}_{2}^{\lor}\wedge{\bf b}_{1}^{\lor},
{\bf a}_{1}^{\lor}\wedge{\bf b}_{1}^{\lor}-{\bf a}_{2}^{\lor}\wedge{\bf b}_{2}^{\lor}\}$,
${\Bbb A}$ is regarded as the lattice ${\Bbb U}(-1)\oplus{\Bbb U}(-1)\oplus\langle-2\rangle$.
\par
The domain of type $IV$ associated to the lattice ${\Bbb A}$ is defined as
$$
\Omega_{\Bbb A}:=\{[x]\in{\bf P}({\Bbb A}\otimes{\bf C});\,\langle x,x\rangle_{\Bbb A}=0,\,\langle x,\overline{x}\rangle_{\Bbb A}>0\}.
$$
By the Riemann-Hodge bilinear relations, $[dz_{1}\wedge dz_{2}]\in\Omega_{\Bbb A}^{+}$, 
where $\Omega_{\Bbb A}^{+}$ is the component of $\Omega_{\Bbb A}$ containing $[dz_{1}\wedge dz_{2}]$.
We define an isomorphism of complex manifolds $\varpi_{\Bbb A}\colon{\frak S}_{2}\to\Omega_{\Bbb A}^{+}$ by
$\varpi_{\Bbb A}(T):=[dz_{1}\wedge dz_{2}]$. 
Since $[dz_{j}]={\bf a}_{j}^{\lor}-{\bf b}_{1}^{\lor}T_{1j}-{\bf b}_{2}^{\lor}T_{2j}$, we have
$$
\varpi_{\Bbb A}(T)
=
{\bf a}_{1}^{\lor}\wedge{\bf a}_{2}^{\lor}
+
\det T\,{\bf b}_{1}^{\lor}\wedge{\bf b}_{2}^{\lor}
+
T_{11}{\bf a}_{2}^{\lor}\wedge{\bf b}_{1}^{\lor}
-
T_{22}{\bf a}_{1}^{\lor}\wedge{\bf b}_{2}^{\lor}
-
T_{12}({\bf a}_{1}^{\lor}\wedge{\bf b}_{1}^{\lor}-{\bf a}_{2}^{\lor}\wedge{\bf b}_{2}^{\lor}).
$$

\subsection
{$2$-cycles on Jacobian Kummer surfaces}
\label{subsect:Periods:Jacobian:Kummer}
\par
Let $A_{T}[2]$ be the points of order $2$ of $A_{T}$.
Let $p\colon\widetilde{A}_{T}\to A_{T}$ be the blowing-up of $A_{T}[2]$.
Let $\widetilde{-1}$ be the involution on $\widetilde{A}_{T}$ induced by  $-1(z)=-z$ on $A_{T}$. 
Then $K_{T}:=\widetilde{A}_{T}/\widetilde{-1}$ is a Kummer surface associated to $A_{T}$. Since $A_{T}$ is the Jacobian variety of $C_{T}$, 
we have $K_{T}=K(C_{T})$.
Let $q\colon\widetilde{A}_{T}\to K_{T}$ be the projection. By \cite[VIII, Prop.\,5.1]{BPV84}, the injective homomorphism
\begin{equation}
\label{eqn:morphism:Abel:K3}
\phi:=q_{!}p^{*}\colon
H^{2}(A_{T},{\bf Z})\hookrightarrow H^{2}(K_{T},{\bf Z})
\end{equation}
satisfies the equality $\langle\phi(l),\phi(m)\rangle=2\langle l,m\rangle$ for all $l,m\in H^{2}(K_{T},{\bf Z})$.
If $T\in{\frak S}_{2}^{o}$ is generic enough, then $\phi({\Bbb A})\subset H^{2}(K_{T},{\bf Z})$ is the transcendental lattice of $K_{T}$.
Set
$$
\Gamma_{12}^{\lor}:=\phi({\bf a}_{1}^{\lor}\wedge{\bf a}_{2}^{\lor}),
\quad
\Gamma_{13}^{\lor}:=\phi({\bf a}_{1}^{\lor}\wedge{\bf b}_{1}^{\lor}),
\quad
\Gamma_{14}^{\lor}:=\phi({\bf a}_{1}^{\lor}\wedge{\bf b}_{2}^{\lor}),
$$
$$
\Gamma_{23}^{\lor}:=\phi({\bf a}_{2}^{\lor}\wedge{\bf b}_{1}^{\lor}),
\quad
\Gamma_{24}^{\lor}:=\phi({\bf a}_{2}^{\lor}\wedge{\bf b}_{2}^{\lor}),
\quad
\Gamma_{34}^{\lor}:=\phi({\bf b}_{1}^{\lor}\wedge{\bf b}_{2}^{\lor}).
$$
In what follows, $\phi({\Bbb A})$ is equipped with the basis
$\{\Gamma^{\lor}_{12},\Gamma^{\lor}_{34},\Gamma^{\lor}_{14},\Gamma^{\lor}_{23},\Gamma^{\lor}_{13}-\Gamma^{\lor}_{24}\}$.
Hence $\phi({\Bbb A})={\Bbb A}(2)$. Let us see the relation between the cycles $\gamma'_{ij}$ and $\Gamma^{\lor}_{ij}$. 
\par
Let $T\in{\frak S}^{o}$.
Define the map $\sigma\colon C_{T}\times C_{T}\to A_{T}$ by $\sigma(x,y):=i_{T}(x)+i_{T}(y)$
and let $\varPi\colon A_{T}\to A_{T}/\pm1$ be the projection.
We define ${\rm sym}:=\varPi\circ\sigma\colon C_{T}\times C_{T}\to A_{T}/\pm1$ (cf. \cite[p.128]{MatsumotoTerasoma12}).
Set $\gamma_{1}:=\alpha_{1}$, $\gamma_{2}:=\alpha_{2}$, $\gamma_{3}:=\beta_{1}$, $\gamma_{4}:=\beta_{2}$.
Let $\overline{\gamma}_{ij}\in H^{2}(K_{T},{\bf Z})$ be the Poincar\'e dual of the proper inverse image of ${\rm sym}_{*}(\gamma_{i}\times\gamma_{j})$ 
with respect to the resolution $r\colon K_{T}\to A_{T}/\pm1$. Namely, if $E_{1}\amalg\ldots\amalg E_{16}=r^{-1}(A_{T}[2])\subset K_{T}$ 
are the exceptional curves of $r\colon K_{T}\to A_{T}/\pm1$,
then $\overline{\gamma}_{ij}$ is the unique element of
$H^{2}_{c}(K_{T}\setminus\amalg_{k=1}^{16}E_{k},{\bf Z})=H^{2}(K_{T},{\bf Z})\cap c_{1}(E_{1})^{\perp}\cap\ldots\cap c_{1}(E_{16})^{\perp}$ satisfying
\begin{equation}
\label{eqn:integral:2-cycle:0}
\langle r^{*}y,\overline{\gamma}_{ij}\rangle
=
\int_{{\rm sym}_{*}(\gamma_{i}\times\gamma_{j})}y,
\qquad
\forall\,y\in H^{2}_{c}(\{A_{T}\setminus A_{T}[2]\}/\pm1,{\bf R}).
\end{equation}
\par
Recall that if $C_{T}\cong C_{\lambda}$ with some $\lambda\in{\bf C}^{6}$, then $K_{T}\cong K(C_{\lambda})\cong W_{N_{\lambda}}$.
Set ${\Bbb T}:=H^{2}(W_{N_{\lambda}},{\bf Z})_{-}$ and regard ${\Bbb T}$ as a primitive sublattice of $H^{2}(K_{T},{\bf Z})$ via the isomorphism
$H^{2}(K_{T},{\bf Z})\cong H^{2}(W_{N_{\lambda}},{\bf Z})$, 
where $H^{2}(W_{N_{\lambda}},{\bf Z})_{-}$ was introduced in Sect.~\ref{subsect:Phi:Freitag:theta}. 
By \cite[Sect.\,3.1]{MatsumotoTerasoma12} $\gamma'_{ij}$ is defined as the orthogonal projection of $\overline{\gamma}_{ij}$ to ${\Bbb T}$.

\begin{lemma}
\label{lemma:2-cycles}
One has $\Gamma^{\lor}_{ij}-\gamma'_{ij}\in\phi({\Bbb A})^{\perp}$.
\end{lemma}

\begin{pf}
If $T\in{\frak S}_{2}^{o}$ is generic,
$\phi({\Bbb A})$ is the smallest sublattice of $H^{2}(K_{T},{\bf Z})$ satisfying $\phi({\Bbb A})\otimes{\bf C}\supset H^{0}(K_{T},\Omega^{2}_{K_{T}})$.
Since $H^{0}(W_{N_{\lambda}},\Omega^{2}_{W_{N_{\lambda}}})\subset{\Bbb T}\otimes{\bf C}$, this implies $\phi({\Bbb A})\subset{\Bbb T}$.
\par
Set ${\bf c}_{1}^{\lor}:={\bf a}_{1}^{\lor}$, ${\bf c}_{2}^{\lor}:={\bf a}_{2}^{\lor}$, ${\bf c}_{3}^{\lor}:={\bf b}_{1}^{\lor}$, ${\bf c}_{4}^{\lor}:={\bf b}_{2}^{\lor}$
and $\gamma_{1}^{\lor}:=\alpha_{1}^{\lor}$, $\gamma_{2}^{\lor}:=\alpha_{2}^{\lor}$, $\gamma_{3}^{\lor}:=\beta_{1}^{\lor}$, $\gamma_{4}^{\lor}:=\beta_{2}^{\lor}$.
For every $x\in H^{2}(A_{T},{\bf Z})\cap[C_{T}]^{\perp}$, we get
\begin{equation}
\label{eqn:integral:2-cycle:1}
\begin{aligned}
\,&
\langle\phi(x),\Gamma_{ij}^{\lor}\rangle
=
2\langle x,{\bf c}_{i}^{\lor}\wedge{\bf c}_{j}^{\lor}\rangle_{H^{2}(A_{T},{\bf Z})}
=
\frac{2}{\deg\sigma}\int_{C_{T}\times C_{T}}\sigma^{*}x\wedge\sigma^{*}{\bf c}_{i}^{\lor}\wedge\sigma^{*}{\bf c}_{j}^{\lor}
\\
&\quad
=\int_{C_{T}\times C_{T}}
\sigma^{*}x
\wedge
({\rm pr}_{1}^{*}\gamma_{i}^{\lor}\wedge{\rm pr}_{2}^{*}\gamma_{j}^{\lor}
-
{\rm pr}_{1}^{*}\gamma_{j}^{\lor}\wedge{\rm pr}_{2}^{*}\gamma_{i}^{\lor})
=
2\int_{\gamma_{i}\times\gamma_{j}}\sigma^{*}x.
\end{aligned}
\end{equation}
By the Mayer-Vietoris exact sequence, we get $H^{2}(A_{T},{\bf Z})=H^{2}_{c}(A_{T}\setminus A_{T}[2],{\bf Z})$.
Let $x\in{\Bbb A}=H^{2}_{c}(A_{T}\setminus A_{T}[2],{\bf Z})\cap[C_{T}]^{\perp}$. Since $x=-x$, we get
\begin{equation}
\label{eqn:integral:2-cycle:2}
\int_{\gamma_{i}\times\gamma_{j}}\sigma^{*}x
=
\int_{\sigma_{*}(\gamma_{i}\times\gamma_{j})}x
=
\int_{\sigma_{*}(\gamma_{i}\times\gamma_{j})}\frac{x+(-1)^{*}x}{2}
=
\frac{1}{2}\int_{{\rm sym}_{*}(\gamma_{i}\times\gamma_{j})}\varPi_{*}x.
\end{equation}
Since $r^{*}\varPi_{*}(x)=q_{*}p^{*}(x)$ for all $x\in H^{2}(A_{T}\setminus A_{T}[2],{\bf Z})$,
we get by \eqref{eqn:integral:2-cycle:0}, \eqref{eqn:integral:2-cycle:1}, \eqref{eqn:integral:2-cycle:2}
\begin{equation}
\label{eqn:integral:2-cycle:3}
\langle\phi(x),\Gamma_{ij}^{\lor}\rangle
=
\int_{{\rm sym}_{*}(\gamma_{i}\times\gamma_{j})}\varPi_{*}x
=
\langle r^{*}\varPi_{*}x,\overline{\gamma}_{ij}\rangle
=
\langle\phi(x),\overline{\gamma}_{ij}\rangle.
\end{equation}
Since $\gamma'_{ij}-\overline{\gamma}_{ij}\in{\Bbb T}^{\perp}\subset\phi({\Bbb A})^{\perp}$ by the definition of $\gamma'_{ij}$
and $\Gamma_{ij}^{\lor}-\overline{\gamma}_{ij}\in\phi({\Bbb A})^{\perp}$ by \eqref{eqn:integral:2-cycle:3}, we get the result.
\end{pf}

\subsection
{Switches and periods of Jacobian Kummer surfaces}
\label{subsect:Switch:Periods:Jacobian:Kummer}
\par
Let $\pi\colon{\mathcal K}\to{\frak S}_{2}$ be the universal family of Kummer surfaces with $\pi^{-1}(T)=K_{T}=K(C_{T})$ for all $T\in{\frak S}_{2}$.
By fixing an order of odd characteristics $(a,b)$, $a,b\in\{0,1\}^{2}$, the branch points of the canonical map 
$\Phi_{|K_{C_{T}}|}\colon C_{T}\to{\bf P}^{1}$ are also ordered, because they are given by the $6$ points
$\{(\frac{\partial\theta}{\partial z_{1}}(\frac{a+Tb}{2},T):\frac{\partial\theta}{\partial z_{2}}(\frac{a+Tb}{2},T))\}_{(a,b):{\rm odd}}$.
Hence $K(C_{T})\cong X_{N_{\lambda}}$ with
$\lambda=\{\frac{\partial\theta}{\partial z_{1}}(\frac{a+Tb}{2},T)/\frac{\partial\theta}{\partial z_{2}}(\frac{a+Tb}{2},T)\}_{(a,b):{\rm odd}}$.
By using this realization, $\pi\colon{\mathcal K}|_{{\frak S}_{2}^{o}}\to{\frak S}_{2}^{o}$ is equipped with the $10$ fixed-point-free involutions $\{\iota_{\langle J\rangle}\}$.
\par
Let $T_{0}\in{\frak S}_{2}^{o}$.
Since $\phi({\Bbb A})$ is the smallest sublattice of $H^{2}(K_{T_{0}},{\bf Z})$ whose complexification contains $\phi(dz_{1}\wedge dz_{2})$ 
if $T_{0}\in{\frak S}_{2}^{o}$ is generic, we get the inclusion
$\phi({\Bbb A})\subset H^{2}(K_{T_{0}},{\bf Z})_{-}=\{l\in H^{2}(K_{T_{0}},{\bf Z});\,\iota_{\langle J\rangle}^{*}l=-l\}$.
Let $\alpha_{\langle J\rangle,T_{0}}\colon H^{2}(K_{T_{0}},{\bf Z})\cong{\Bbb L}_{K3}$ be a marking satisfying 
\eqref{eqn:alpha:+:-} with respect to $\iota_{\langle J\rangle}$.
Then $\alpha_{\langle J\rangle,T_{0}}(\phi({\Bbb A}))\subset\LAM$. 
Since ${\frak S}_{2}$ is contractible, $\alpha_{\langle J\rangle,T_{0}}$ extends to a marking
$\alpha_{\langle J\rangle}\colon R^{2}\pi_{*}{\bf Z}\cong{\Bbb L}_{K3}$ for the family $\pi\colon{\mathcal K}\to{\frak S}_{2}$.
Then $\alpha_{\langle J\rangle}|_{T}(\phi({\Bbb A}))=\alpha_{\langle J\rangle,T_{0}}(\phi({\Bbb A}))$ for all $T\in{\frak S}_{2}$.
\par
The period map for the marked family of Enriques surfaces
$(\pi\colon{\mathcal K}|_{{\frak S}_{2}^{o}}/\iota_{\langle J\rangle}\to{\frak S}_{2}^{o},\alpha_{\langle J\rangle})$
is given by
\begin{equation}
\label{eqn:period:Kummer:switch}
\begin{aligned}
\alpha_{\langle J\rangle}(\phi(\varpi_{\Bbb A}(T)))
&=
[
\alpha_{\langle J\rangle}(\Gamma^{\lor}_{12})
+
\det T\,\alpha_{\langle J\rangle}(\Gamma^{\lor}_{34})+
\\
&\qquad
T_{11}\alpha_{\langle J\rangle}(\Gamma^{\lor}_{23})
-
T_{22}\alpha_{\langle J\rangle}(\Gamma^{\lor}_{14})
-
T_{12}\alpha_{\langle J\rangle}(\Gamma^{\lor}_{13}-\Gamma^{\lor}_{24})
],
\end{aligned}
\end{equation}
which extends to an embedding $\varpi_{\langle J\rangle}\colon{\frak S}_{2}\hookrightarrow\Omega_{\LAM}^{+}$ defined by
$$
\varpi_{\langle J\rangle}:=\alpha_{\langle J\rangle}\circ\phi\circ\varpi_{\Bbb A}.
$$

\subsection
{Infinite product expansion of theta constants of genus $2$}
\label{subsect:Infinite:product:theta:constants}
\par 
Let $\ell\in\{1,2\}$ be the level of the primitive isotropic vector $\alpha_{\langle J\rangle}(\Gamma_{34}^{\lor})$ in $\LAM$. 
Choosing $\alpha_{\langle J\rangle}$ suitably, we may and will assume by \cite[Prop.\,4.5]{Sterk91} that 
$\alpha_{\langle J\rangle}(\Gamma_{34}^{\lor})={\frak e}_{\ell}$.
Let $\Gamma_{34}\in H_{2}(K_{T},{\bf Z})$ be the Poincar\'e dual of $\Gamma_{34}^{\lor}$.
We define $z_{\langle J\rangle}(T)\in{\Bbb M}_{\ell}\otimes{\bf R}+i\,{\mathcal C}_{{\Bbb M}_{\ell}}$ by the equation
\begin{equation}
\label{eqn:definition:period:z}
\varpi_{\langle J\rangle}(T)
=
\left[
-(z_{\langle J\rangle}(T)^{2}/2){\frak e}_{\ell}+({\frak f}_{\ell}/\ell)+z_{\langle J\rangle}(T)
\right].
\end{equation}

\begin{lemma}
\label{lemma:expression:z}
There exist constants $A,B,C,D\in{\Bbb M}_{\ell}$ depending on $\langle J\rangle$ such that 
$$
z_{\langle J\rangle}(T)=(A+BT_{11}+CT_{12}+DT_{22})/2,
\qquad
\forall\,T\in{\frak S}_{2}.
$$
\end{lemma}

\begin{pf}
By \eqref{eqn:period:Kummer:switch} and \eqref{eqn:definition:period:z}, there exists $\lambda\in{\bf C}^{*}$ with
$$
\begin{aligned}
\,&
\alpha_{\langle J\rangle}(\Gamma^{\lor}_{12})
+
\det T\,{\frak e}_{\ell}
+
T_{11}\,\alpha_{\langle J\rangle}(\Gamma^{\lor}_{23})
-
T_{22}\,\alpha_{\langle J\rangle}(\Gamma^{\lor}_{14})
-
T_{12}\,\alpha_{\langle J\rangle}(\Gamma^{\lor}_{13}-\Gamma^{\lor}_{24})
\\
&=
\lambda
\left\{
-(z_{\langle J\rangle}(T)^{2}/2){\frak e}_{\ell}+({\frak f}_{\ell}/\ell)+z_{\langle J\rangle}(T)
\right\}.
\end{aligned}
$$
Comparing the inner products with ${\frak e}_{\ell}$, we get $\lambda=2$.
Set
$$
{\frak f}':=\alpha_{\langle J\rangle}(\Gamma^{\lor}_{12}),
\quad
{\frak a}:=\alpha_{\langle J\rangle}(\Gamma^{\lor}_{23}),
\quad
{\frak b}:=\alpha_{\langle J\rangle}(\Gamma^{\lor}_{14}),
\quad
{\frak c}:=\alpha_{\langle J\rangle}(\Gamma^{\lor}_{13}-\Gamma^{\lor}_{24}).
$$
Since 
$\langle{\frak f}_{\ell}/\ell,{\frak a}\rangle, \langle{\frak f}_{\ell}/\ell,{\frak b}\rangle, \langle{\frak f}_{\ell}/\ell,{\frak c}\rangle, \langle{\frak f}_{\ell}/\ell,{\frak f}'\rangle\in{\bf Z}$,
we get the desired expression
$$
\begin{aligned}
z_{\langle J\rangle}(T)
&=
\frac{{\frak f}'}{2}
-
\langle{\frak f}_{\ell}/\ell,{\frak f}'+T_{11}{\frak a}-T_{22}{\frak b}-T_{12}{\frak c}\rangle\frac{{\frak e}_{\ell}}{2}
-\frac{{\frak f}_{\ell}}{\ell}
+
\frac{1}{2}(T_{11}{\frak a}-T_{22}{\frak b}-T_{12}{\frak c})
\\
&=(A+BT_{11}+CT_{12}+DT_{22})/2
\end{aligned}
$$
with $A={\frak f}'-(2/\ell){\frak f}_{\ell}$, $B=-\langle{\frak f}_{\ell}/\ell,{\frak a}\rangle{\frak e}_{\ell}+{\frak a}$, 
$C=\langle{\frak f}_{\ell}/\ell,{\frak c}\rangle{\frak e}_{\ell}-{\frak c}$, $D=\langle{\frak f}_{\ell}/\ell,{\frak b}\rangle{\frak e}_{\ell}-{\frak b}$.
\end{pf}

\begin{theorem}
\label{thm:infinite:product:theta:constant}
Let $({\bf a},{\bf b})\in Ev$ correspond to $\langle J\rangle$. Then for all $T\in {\frak S}_{2}$,
$$
(\varpi_{\langle J\rangle}^{*}\Phi_{\ell})(T)
=
\Phi_{\ell}\left(z_{\langle J\rangle}(T)\right)
=
\pm\theta_{\Re(\frac{\bf a}{1+i})\Re(\frac{\bf b}{1+i})}(T)^{8}
=
\pm\theta_{\Im(\frac{\bf a}{1+i})\Im(\frac{\bf b}{1+i})}(T)^{8}.
$$
\end{theorem}

\begin{pf}
For short, write $z(T)$ for $z_{\langle J\rangle}(T)$ and $\alpha$ for $\alpha_{\langle J\rangle}$.
Let $\lambda\in{\bf C}^{6}$ be such that $K_{T}=K(C_{\lambda})$.
Since $K_{T}\cong W_{N_{\lambda}}$, $\eta_{N_{\lambda}}/\int_{\Gamma_{34}}\eta_{N_{\lambda}}$ is a canonical form on $K_{T}$ 
whose integration over $\Gamma_{34}$ is normalized to $1$. 
Hence we get the equality
$$
\alpha
\left(
\eta_{N_{\lambda}}/\langle\Gamma^{\lor}_{34},\eta_{N_{\lambda}}\rangle
\right)
=
-(z(T)^{2}/2){\frak e}_{\ell}+({\frak f}_{\ell}/\ell)+z(T),
$$
from which we get
\begin{equation}
\label{eqn:norm:Im:z:1}
\begin{aligned}
\,&
2\langle\Im z(T),\Im z(T)\rangle
=
\left\langle
-(z(T)^{2}/2){\frak e}_{\ell}+({\frak f}_{\ell}/\ell)+z(T)
,
\overline{-(z(T)^{2}/2){\frak e}_{\ell}+({\frak f}_{\ell}/\ell)+z(T)}
\right\rangle
\\
&=
\langle\alpha(\eta_{N_{\lambda}}),\alpha(\overline{\eta}_{N_{\lambda}})\rangle\cdot\left|\int_{\Gamma_{34}}\eta_{N_{\lambda}}\right|^{-2}
=
\left(
\int_{Z_{N_{\lambda}}}\eta_{N_{\lambda}}\wedge\overline{\eta}_{N_{\lambda}}
\right)
\cdot
\left|\int_{\Gamma_{34}}\eta_{N_{\lambda}}\right|^{-2}.
\end{aligned}
\end{equation}
Since 
$\eta_{34}(N_{\lambda})=\langle\gamma'_{34},\eta_{N_{\lambda}}\rangle/2=
\langle\Gamma^{\lor}_{34},\eta_{N_{\lambda}}\rangle/2=(\int_{\Gamma_{34}}\eta_{N_{\lambda}})/2$
by Lemma~\ref{lemma:2-cycles}, we get 
\begin{equation}
\label{eqn:norm:Im:z:2}
\det\Im T
=
\frac{1}{8}
\int_{Z_{N_{\lambda}}}\left(\frac{\eta_{N_{\lambda}}}{\eta_{34}(N_{\lambda})}\right)\wedge\overline{\left(\frac{\eta_{N_{\lambda}}}{\eta_{34}(N_{\lambda})}\right)}
=
\langle\Im z(T),\Im z(T)\rangle
\end{equation}
by  \eqref{eqn:Bergman:kernal:D}, \eqref{eqn:norm:Im:z:1}.
We deduce from Theorem~\ref{thm:Phi:Freitag:theta} and \eqref{eqn:norm:Im:z:2} that
\begin{equation}
\label{eqn:norm:Phi:Theta:1}
\|\Phi(Y_{N_{\lambda},\langle J\rangle})\|
=
\|\Theta_{\langle J\rangle}(T)\|^{4}
=
\langle\Im z(T),\Im z(T)\rangle^{2}
\left|
\Theta_{\langle J\rangle}(T)
\right|^{4}.
\end{equation}
Comparing \eqref{eqn:norm:Phi:Theta:1} with the definition
$\|\Phi(Y_{N_{\lambda},\langle J\rangle})\|=\langle\Im z(T),\Im z(T)\rangle^{2}|\Phi_{\ell}(z(T))|$, 
we get the equality of functions on ${\frak S}_{2}$:
$$
|\Phi_{\ell}(z(T))|=\left|\Theta_{\langle J\rangle}(T)\right|^{4}.
$$
Hence there exists by \eqref{eqn:theta:Freitag:Riemann} a constant $C\in{\bf C}$ with $|C|=1$ such that for all $T\in{\frak S}_{2}$
\begin{equation}
\label{eqn:Phi:Riemann:theta:1}
\Phi_{\ell}(z(T))^{2}
=
C\,\Theta_{\langle J\rangle}(T)^{8}
=
C\,\theta_{\Re(\frac{\bf a}{1+i})\Re(\frac{\bf b}{1+i})}(T)^{16}
=
C\,\theta_{\Im(\frac{\bf a}{1+i})\Im(\frac{\bf b}{1+i})}(T)^{16}.
\end{equation}
\par
Set $q_{mn}:=\exp(\pi iT_{mn})$ and $(a,b):=(\Re(\frac{\bf a}{1+i}),\Re(\frac{\bf b}{1+i}))$.
By \eqref{eqn:Phi:expansion:level1}, \eqref{eqn:Phi:expansion:level2} and Lemma~\ref{lemma:expression:z},
we get $\Phi_{\ell}(z(T))\in{\bf Z}\{q_{11},q_{12},q_{22}\}$. 
By definition, we get $\theta_{a,b}(T)\in{\bf Z}\{q_{11}^{1/4},q_{12}^{1/4},q_{22}^{1/4}\}$ for any even $(a,b)$. 
Since $\theta_{a,b}(T)^{8}$ is a Siegel modular form for the principal congruence subgroup of  level $2$ by \cite[p.190, Eq.(5.2)]{Mumford83},
we get $\theta_{a,b}(T+2B)^{8}=\theta_{a,b}(T)^{8}$ for any integral $2\times2$-matrix $B$.
Hence $\theta_{a,b}(T)^{8}\in{\bf Z}\{q_{11},q_{12},q_{22}\}$ for any even $(a,b)$. 
Since $\Phi_{\ell}(z(T))\in{\bf Z}\{q_{11},q_{12},q_{22}\}$ and $\theta_{a,b}(T)^{8}\in{\bf Z}\{q_{11},q_{12},q_{22}\}$,
we get $C\in{\bf Q}$ by \eqref{eqn:Phi:Riemann:theta:1}. 
Since $|C|=1$, this implies $C=\pm1$. 
By \eqref{eqn:Phi:Riemann:theta:1} again, $\Phi_{\ell}(z(T))^{2}=\pm\theta_{a,b}(T)^{16}$. 
If $\Phi_{\ell}(z(T))^{2}=-\theta_{a,b}(T)^{16}$, then $\Phi_{\ell}(z(T))=\pm\sqrt{-1}\theta_{a,b}(T)^{8}\in\sqrt{-1}{\bf Z}\{q_{11},q_{12},q_{22}\}$. 
Since $\Phi_{\ell}(z(T))\in{\bf Z}\{q_{11},q_{12},q_{22}\}$, we get a contradiction $\Phi_{\ell}(z(T))\equiv0$. This proves $C=1$.
\end{pf}

\begin{corollary}
\label{cor:infinite:product:theta:genus:2}
For all even $(a,b)\in\{0,\frac{1}{2}\}^{4}$, $\theta_{a,b}(T)^{8}$ admits an infinite product expansion of Borcherds type.
\end{corollary}

Let $\Delta_{5}(T):=\prod_{(a,b)\,{\rm even}}\theta_{a,b}(T)$ be the product of all even theta constants.
By Corollary~\ref{cor:infinite:product:theta:genus:2}, $\Delta_{5}(T)^{4}$ is expressed as an infinite product of Borcherds type.
In fact, Gritsenko-Nikulin \cite{GritsenkoNikulin97} proved that $\Delta_{5}(T)$ is a Borcherds product.

\section{Kummer surfaces of product type and involutions of even type}
\label{sec:Kummer}
In this section, we prove ${\frak C}={\frak c}_{1}={\frak c}_{2}=1$ in \eqref{eqn:constant:Phi:generic} and 
Proposition~\ref{proposition:Thomae:Phi:with:constant} by studying $\Phi$ for those Enriques surfaces 
which are the quotient of a Kummer surface of product type by involutions of {\it even} type. 
Let ${\frak H}$ be the complex upper-half plane.

\subsection
{Elliptic functions}
\label{subsect:elliptic:functions}
\par
For $\tau\in{\frak H}$, set $E_{\tau}:={\bf C}/2{\bf Z}+2\tau{\bf Z}$.
Let $\wp(u;1,\tau)$ be the Weierstrass $\wp$-function on the $u$-plane associated to the lattice ${\bf Z}+\tau{\bf Z}$.
The map $E_{\tau}={\bf C}/2{\bf Z}+2\tau{\bf Z}\ni[u]\mapsto(\wp(\frac{u}{2};1,\tau):\wp'(\frac{u}{2};1,\tau):1)\in{\bf P}^{2}$ is an isomorphism 
from $E_{\tau}$ to the cubic curve defined by the affine equation 
$
w^{2}=4(z-e_{1})(z-e_{2})(z-e_{3}),
$
where
$e_{1}=\wp(\frac{1}{2};1,\tau)$,
$e_{2}=\wp(\frac{1+\tau}{2};1,\tau)$,
$e_{3}=\wp(\frac{\tau}{2};1,\tau)$.
\par
By \cite[p.213 Eq.(3)]{HurwitzCourant29}, Jacobi's elliptic functions and the $\wp$-function are related as follows:
$$
\begin{array}{ll}
{\rm sn}(\sqrt{e_{1}-e_{3}}u)
&=
\sqrt{(e_{1}-e_{3})/(\wp(u;1,\tau)-e_{3})},
\\
{\rm cn}(\sqrt{e_{1}-e_{3}}u)
&=
\sqrt{(\wp(u;1,\tau)-e_{1})/(\wp(u;1,\tau)-e_{3})},
\\
{\rm dn}(\sqrt{e_{1}-e_{3}}u)
&=
\sqrt{(\wp(u;1,\tau)-e_{2})/(\wp(u;1,\tau)-e_{3})}.
\end{array}
$$
By \cite[p.215, Tabelle II]{HurwitzCourant29},
${\rm sn}(\sqrt{e_{1}-e_{3}}u)$ is an odd function with period lattice $\Gamma_{1}:=2{\bf Z}+\tau{\bf Z}$ 
and ${\rm cn}(\sqrt{e_{1}-e_{3}}u)$ (resp. ${\rm dn}(\sqrt{e_{1}-e_{3}}u)$) 
is an even function with period lattice $\Gamma_{2}:=2{\bf Z}+(1+\tau){\bf Z}$ (resp. $\Gamma_{3}:={\bf Z}+2\tau{\bf Z}$).
We regard ${\rm sn}(\sqrt{e_{1}-e_{3}}u)$, ${\rm cn}(\sqrt{e_{1}-e_{3}}u)$,
${\rm dn}(\sqrt{e_{1}-e_{3}}u)$ as periodic meromorphic functions on ${\bf C}$ with period lattice 
$\Gamma_{1}\cap\Gamma_{2}\cap\Gamma_{3}=2{\bf Z}+2\tau{\bf Z}$. 
Namely, they are meromorphic functions on $E_{\tau}$.
\par
By \cite[p.204]{HurwitzCourant29}, \cite[p.69]{Mumford83}, the theta constants are defined as
\begin{equation}
\label{eqn:infinite:product:theta:constant:genus:1}
\begin{array}{ll}
\theta_{0}(\tau)=\theta_{01}(\tau)
&:=
\prod_{n=1}^{\infty}(1-e^{2\pi in\tau})(1-e^{\pi i(2n-1)\tau})^{2},
\\
\theta_{2}(\tau)=\theta_{10}(\tau)
&:=
2e^{\pi i\tau/4}
\prod_{n=1}^{\infty}(1-e^{2\pi in\tau})(1+e^{2\pi in\tau})^{2},
\\
\theta_{3}(\tau)=\theta_{00}(\tau)
&:=
\prod_{n=1}^{\infty}(1-e^{2\pi in\tau})(1+e^{\pi i(2n-1)\tau})^{2}.
\end{array}
\end{equation}
Then $\sqrt{e_{2}-e_{3}}=\pi\,\theta_{2}(\tau)^{2}$ and $\sqrt{e_{1}-e_{3}}=\pi\,\theta_{3}(\tau)^{2}$ by \cite[p.202 Eq.(6)]{HurwitzCourant29}.

\subsection
{A $(2,2,2)$-model of a Kummer surface of product type}
\label{subsec:another:model:Kummer:surface}
\par
For $\tau_{1},\tau_{2}\in{\frak H}$ and $z_{1},z_{2}\in{\bf C}$, we define functions $x_{1},x_{2},y_{0},y_{1},y_{2}$ on $E_{\tau_{1}}\times E_{\tau_{2}}$ by
\begin{equation}
\label{eqn:uniformization}
\begin{array}{lll}
x_{1}
&:=
{\rm cn}(\sqrt{e_{1}^{1}-e_{3}^{1}}z_{1})
&=
\sqrt{\frac{\wp(z_{1};1,\tau_{1})-e^{1}_{1}}{\wp(z_{1};1,\tau_{1})-e^{1}_{3}}},
\\
x_{2}
&:=
{\rm dn}(\sqrt{e^{1}_{1}-e^{1}_{3}}z_{1})
&=
\sqrt{\frac{\wp(z_{1};1,\tau_{1})-e^{1}_{2}}{\wp(z_{1};1,\tau_{1})-e^{1}_{3}}},
\\
y_{0}
&:=
\frac{{\rm sn}(\sqrt{e^{1}_{1}-e^{1}_{3}}z_{1})}{{\rm sn}(\sqrt{e^{2}_{1}-e^{2}_{3}}z_{2})}
&=
\sqrt{\frac{e^{1}_{1}-e^{1}_{3}}{e^{2}_{1}-e^{2}_{3}}}
\sqrt{\frac{\wp(z_{2};1,\tau_{2})-e^{2}_{3}}{\wp(z_{1};1,\tau_{1})-e^{1}_{3}}},
\\
y_{1}
&:=
\frac{{\rm sn}(\sqrt{e^{1}_{1}-e^{1}_{3}}z_{1})}{{\rm sn}(\sqrt{e^{2}_{1}-e^{2}_{3}}z_{2})}
{\rm cn}(\sqrt{e^{2}_{1}-e^{2}_{3}}z_{2})
&=
\sqrt{\frac{e^{1}_{1}-e^{1}_{3}}{e^{2}_{1}-e^{2}_{3}}}
\sqrt{\frac{\wp(z_{2};1,\tau_{2})-e^{2}_{1}}{\wp(z_{1};1,\tau_{1})-e^{1}_{3}}},
\\
y_{2}
&:=
\frac{{\rm sn}(\sqrt{e^{1}_{1}-e^{1}_{3}}z_{1})}{{\rm sn}(\sqrt{e^{2}_{1}-e^{2}_{3}}z_{2})}
{\rm dn}(\sqrt{e^{2}_{1}-e^{2}_{3}}z_{2})
&=
\sqrt{\frac{e^{1}_{1}-e^{1}_{3}}{e^{2}_{1}-e^{2}_{3}}}
\sqrt{\frac{\wp(z_{2};1,\tau_{2})-e^{2}_{2}}{\wp(z_{1};1,\tau_{1})-e^{1}_{3}}},
\end{array}
\end{equation}
where 
$e^{i}_{1}:=\wp(\frac{1}{2};1,\tau_{i})$,
$e^{i}_{2}:=\wp(\frac{1+\tau_{i}}{2};1,\tau_{i})$,
$e^{i}_{3}:=\wp(\frac{\tau_{i}}{2};1,\tau_{i})$.
Consider the map
$$
f\colon
E_{\tau_{1}}\times E_{\tau_{2}}
\ni([z_{1}],[z_{2}])
\to
(1:x_{1}:x_{2}:y_{0}:y_{1}:y_{2})
\in{\bf P}^{5}.
$$
Then $f(-[z_{1}],-[z_{2}])=f([z_{1}],[z_{2}])$ and $f$ is well-defined if and only if $\wp(z_{1};1,\tau_{1})\not=\infty$ or $\wp(z_{2};1,\tau_{2})\not=\infty$. 
Since $\wp(z_{1};1,\tau_{1})=\wp(z_{2};1,\tau_{2})=\infty$ if and only if $([z_{1}],[z_{2}])$ is a point of order $2$ of $E_{\tau_{1}}\times E_{\tau_{2}}$
and since the indeterminacy of $f$ is resolved by blowing-up the points of order $2$,
$f$ is a birational morphism from $K(E_{\tau_{1}}\times E_{\tau_{2}})$ to
\begin{equation}
\label{eqn:Kummer:product:type:(2,2,2)}
X_{(\lambda_{1},\lambda_{2})}
:=
\left\{
(x_{0}:x_{1}:x_{2}:y_{0}:y_{1}:y_{2})\in{\bf P}^{5};\,
\begin{array}{rl}
(1-\lambda_{1})x_{0}^{2}+\lambda_{1}x_{1}^{2}-x_{2}^{2}
&=
0,
\\
\lambda_{2}x_{0}^{2}-\lambda_{2}x_{1}^{2}-y_{0}^{2}+y_{2}^{2}
&=
0,
\\
x_{0}^{2}-x_{1}^{2}-y_{0}^{2}+y_{1}^{2}
&=
0
\end{array}
\right\},
\end{equation}
where $\lambda_{i}=(e^{i}_{2}-e^{i}_{3})/(e^{i}_{1}-e^{i}_{3})=\theta_{2}(\tau_{i})^{4}/\theta_{3}(\tau_{i})^{4}$ 
is the cross ratio of the branch points of $\wp_{E_{\tau_{i}}}$.
Then $K_{(\tau_{1},\tau_{2})}$ is isomorphic to the minimal resolution of $X_{(\lambda_{1},\lambda_{2})}$. 
In the notation of Sect.\,\ref{sec:expression:Phi:special}, we deduce from \eqref{eqn:Kummer:product:type:(2,2,2)} that
$X_{(\lambda_{1},\lambda_{2})}=X_{M(\lambda_{1},\lambda_{2})}$, where
\begin{equation}
\label{eqn:matrix:Kummer:product:type:(2,2,2)}
M(\lambda_{1},\lambda_{2})
:=
\begin{pmatrix}
\lambda_{1}-1&-\lambda_{1}&1&0&0&0
\\
\lambda_{2}&-\lambda_{2}&0&-1&0&1
\\
1&-1&0&-1&1&0
\end{pmatrix}
\in M_{3,6}({\bf C}).
\end{equation}
Since $\Delta_{123}(M(\lambda_{1},\lambda_{2}))=\Delta_{456}(M(\lambda_{1},\lambda_{2}))=0$,
we get $M(\lambda_{1},\lambda_{2})\not\in M^{o}_{3,6}({\bf C})$. Indeed $X_{(\lambda_{1},\lambda_{2})}$ has nodes. 
For $\langle J\rangle\not=\binom{123}{456}$, the minimal resolution of
$$
Y_{(\lambda_{1},\lambda_{2}),\langle J\rangle}
:=
X_{(\lambda_{1},\lambda_{2})}/\iota_{\langle J\rangle}
$$
is an Enriques surface and the value $\|\Phi(Y_{(\lambda_{1},\lambda_{2})),\langle J\rangle}\|$ is well defined.
\par
The involutions $\iota_{\langle J\rangle}$ are {\em involutions of even type} on $K(E_{\tau_{1}}\times E_{\tau_{2}})$ in the sense of \cite{Mukai11}. 
Let $\epsilon$ be the involution on $K(E_{\tau_{1}}\times E_{\tau_{2}})$ induced by the involution 
$-1_{E_{\tau_{1}}}\times{\rm id}_{E_{\tau_{2}}}$ on $E_{\tau_{1}}\times E_{\tau_{2}}$.
For $a=(a_{1},a_{2})\in(\frac{1}{2}{\bf Z}/{\bf Z})^{\oplus2}\setminus\{0\}$ and $b=(b_{1},b_{2})\in(\frac{1}{2}{\bf Z}/{\bf Z})^{\oplus2}\setminus\{0\}$,
let $\sigma_{(a,b)}$ be the involution on $K(E_{\tau_{1}}\times E_{\tau_{2}})$ induced by the translation 
$(z_{1},z_{2})\mapsto(z_{1}+a_{1}+a_{2}\tau_{1},z_{2}+b_{1}+b_{2}\tau_{2})$ on $E_{\tau_{1}}\times E_{\tau_{2}}$. 
Using the transformation rules for ${\rm sn}(\sqrt{e_{1}-e_{3}}u)$, ${\rm cn}(\sqrt{e_{1}-e_{3}}u)$, ${\rm dn}(\sqrt{e_{1}-e_{3}}u)$
under the translations $u\mapsto u+1$, $u\mapsto u+\tau$, $u\mapsto u+1+\tau$ (cf. \cite[p.215 Tabelle I]{HurwitzCourant29}), we have
$$
\epsilon\circ\sigma_{(a,b)}=\iota_{\langle J\rangle},
$$
where the correspondence between $\binom{a}{b}=\binom{a_{1}\,a_{2}}{b_{1}\,b_{2}}$ and $\langle J\rangle\not=\binom{123}{456}$ is given as follows:
$$
\begin{array}{cccccccccccc}
\binom{a}{b}
&
\vline
&
\binom{1\,0}{1\,0}
&
\binom{1\,0}{0\,1}
&
\binom{1\,0}{1\,1}
&
\binom{0\,1}{1\,0}
&
\binom{0\,1}{0\,1}
&
\binom{0\,1}{1\,1}
&
\binom{1\,1}{1\,0}
&
\binom{1\,1}{0\,1}
&
\binom{1\,1}{1\,1}
\\
\hline
\langle J\rangle
&
\vline
&
\binom{135}{246}
&
\binom{134}{256}
&
\binom{136}{245}
&
\binom{146}{235}
&
\binom{156}{234}
&
\binom{145}{236}
&
\binom{125}{346}
&
\binom{124}{356}
&
\binom{126}{345}
\end{array}
$$
Here we used the notation 
$\binom{a_{1}\,a_{2}}{b_{1}\,b_{2}}:=\binom{(a_{1}+a_{2}\tau_{1})/2}{(b_{1}+b_{2}\tau_{2})/2}\in E_{\tau_{1}}\times E_{\tau_{2}}$.

\subsection
{Periods of Kummer surfaces of product type}
\label{subsec:period:Kummer:product:type}
\par
For $(\tau_{1},\tau_{2})\in{\frak H}\times{\frak H}$, set
$$
{\Bbb K}:=\phi\left({\rm pr}_{1}^{*}H^{1}(E_{\tau_{1}},{\bf Z})\wedge{\rm pr}_{2}^{*}H^{1}(E_{\tau_{2}},{\bf Z})\right)
\subset 
H^{2}(E_{\tau_{1}}\times E_{\tau_{2}},{\bf Z}),
$$
where ${\rm pr}_{i}\colon E_{\tau_{1}}\times E_{\tau_{2}}\to E_{\tau_{i}}$ is the projection. 
If $(\tau_{1},\tau_{2})\in{\frak H}\times{\frak H}$ is generic enough, then ${\Bbb K}$ is the transcendental lattice of $K(E_{\tau_{1}}\times E_{\tau_{2}})$.
\par
Let $\{\alpha_{i},\beta_{i}\}$ be the canonical basis of $H_{1}(E_{\tau_{i}},{\bf Z})$ such that
\begin{equation}
\label{eqn:definition:symplectic:basis}
\int_{\alpha_{i}}dz_{i}=2\tau_{i},
\qquad
\int_{\beta_{i}}dz_{i}=2,
\end{equation}
where $z_{i}$ is the coordinate of $E_{\tau_{i}}$.
We define the cocycles $\alpha_{i}^{\lor}$, $\beta_{i}^{\lor}$, ${\bf a}_{i}^{\lor}$, ${\bf b}_{i}^{\lor}$, $\Gamma_{ij}^{\lor}$ in the same way as in
Sects~\ref{subsect:Periods:ppas}, \ref{subsect:Periods:Jacobian:Kummer}. 
Then ${\Bbb K}$ is equipped with the basis $\{\Gamma_{13}^{\lor},\Gamma_{24}^{\lor},\Gamma_{14}^{\lor},\Gamma_{23}^{\lor}\}$
and is regarded as ${\Bbb U}(-2)\oplus{\Bbb U}(-2)$.
Since $[dz_{i}]=2(\alpha_{i}^{\lor}-\tau_{i}\beta_{i}^{\lor})$ by \eqref{eqn:definition:symplectic:basis}, we get
\begin{equation}
\label{eqn:period:Enriques:Kummer:quotient}
\phi({\bf C}[dz_{1}\wedge dz_{2}])
=
{\bf C}(\Gamma^{\lor}_{12}+\tau_{1}\tau_{2}\,\Gamma^{\lor}_{34}+\tau_{1}\Gamma^{\lor}_{23}-\tau_{2}\Gamma^{\lor}_{14}).
\end{equation}

\subsection
{The restriction of $\Phi$ to ${\frak H}\times{\frak H}$}
\label{subsec:restriction:Phi:H:H:even:type}
\par
Let $\pi\colon{\mathcal X}\to{\frak H}\times{\frak H}$ be the family of Kummer surfaces such that 
$\pi^{-1}(\tau_{1},\tau_{2})=\widetilde{X}_{(\lambda_{1},\lambda_{2})}$, where $\widetilde{X}_{(\lambda_{1},\lambda_{2})}$ is the minimal resolution of
$X_{(\lambda_{1},\lambda_{2})}$. 
For a partition $\langle J\rangle$, the fiberwise involution $\iota_{\langle J\rangle}$ on $X_{(\lambda_{1},\lambda_{2})}$ induces an involution on ${\mathcal X}$, 
which is again denoted by $\iota_{\langle J\rangle}$.
Since ${\frak H}\times{\frak H}$ is contractible, the family $\pi\colon{\mathcal X}\to{\frak H}\times{\frak H}$ is topologically trivial
and hence there exists a marking $\alpha_{\langle J\rangle}\colon R^{2}\pi_{*}{\bf Z}\cong{\Bbb L}_{K3}$ such that 
the condition \eqref{eqn:alpha:+:-} is satisfied fiberwise for $\iota_{\langle J\rangle}$.
Choosing $\alpha_{\langle J\rangle}$ suitably, we may also assume by \cite[Prop.\,4.5]{Sterk91} that 
$\alpha_{\langle J\rangle}(\Gamma^{\lor}_{34})={\frak e}_{\ell}\in{\Bbb U}(\ell)$ when the isotropic vector $\alpha_{\langle J\rangle}(\Gamma^{\lor}_{34})$ 
has level $\ell$ in $\LAM$. 
We fix such a marking $\alpha_{\langle J\rangle}$ of $\pi\colon{\mathcal X}\to{\frak H}\times{\frak H}$ and set
$$
\Omega_{\alpha_{\langle J\rangle}({\Bbb K})}:=\{[\eta] \in\Omega_{\LAM}^{+};\,\eta\in\alpha_{\langle J\rangle}({\Bbb K})\otimes{\bf C}\}.
$$
Let $\varpi_{\langle J\rangle}\colon{\frak H}\times{\frak H}\to\Omega_{\alpha_{\langle J\rangle}({\Bbb K})}^{+}$ be the period map 
for the marked family of Enriques surfaces $(\pi\colon{\mathcal X}/\iota_{\langle J\rangle}\to{\frak H}\times{\frak H},\alpha_{\langle J\rangle})$.
By \eqref{eqn:period:Enriques:Kummer:quotient}, we get 
$$
\varpi_{\langle J\rangle}(\tau_{1},\tau_{2})
=
[\alpha_{\langle J\rangle}(\Gamma^{\lor}_{12})+\tau_{1}\tau_{2}\alpha_{\langle J\rangle}(\Gamma^{\lor}_{34})
+\tau_{1}\alpha_{\langle J\rangle}(\Gamma^{\lor}_{23})-\tau_{2}\alpha_{\langle J\rangle}(\Gamma^{\lor}_{14})].
$$
Set
${\frak f}':=\alpha_{\langle J\rangle}(\Gamma^{\lor}_{12})$, ${\frak e}_{\ell}:=\alpha_{\langle J\rangle}(\Gamma^{\lor}_{34})$,
${\frak a}:=\alpha_{\langle J\rangle}(\Gamma^{\lor}_{23})$, ${\frak b}:=\alpha_{\langle J\rangle}(\Gamma^{\lor}_{14})$
as in the proof of Lemma~\ref{lemma:expression:z}.
Define $z_{\langle J\rangle}\in{\Bbb M}_{\ell}\otimes{\bf R}+i\,{\mathcal C}_{{\Bbb M}_{\ell}}^{+}$ by
$$
z_{\langle J\rangle}(\tau_{1},\tau_{2})
:=
z_{\langle J\rangle}({\rm diag}(\tau_{1},\tau_{2}))
=
\frac{{\frak f}'}{2}
-
\langle{\frak f}_{\ell}/\ell,{\frak f}'+\tau_{1}{\frak a}-\tau_{2}{\frak b}\rangle\frac{{\frak e}_{\ell}}{2}
-\frac{\frak f}{\ell}
+
\frac{1}{2}(\tau_{1}{\frak a}-\tau_{2}{\frak b}).
$$
Then $z_{\langle J\rangle}(\tau_{1},\tau_{2})=(A+B\tau_{11}+D\tau_{22})/2$ by Lemma~\ref{eqn:definition:period:z},
where the constants $A,B,D\in{\Bbb M}_{\ell}$ were given in the proof of Lemma~\ref{eqn:definition:period:z}.
Since $z_{\langle J\rangle}(\tau_{1},\tau_{2})$ satisfies the relation \eqref{eqn:definition:period:z}, it follows from the definition of $\Phi_{\ell}(z)$ that
$$
(\varpi_{\langle J\rangle}^{*}\Phi_{\ell})(\tau_{1},\tau_{2})=\Phi_{\ell}(z_{\langle J\rangle}(\tau_{1},\tau_{2})).
$$

\begin{lemma}
\label{lemma:limit:Phi:level:1:cusp}
Let $\ell\in\{1,2\}$ be the level of $\alpha_{\langle J\rangle}(\Gamma^{\lor}_{34})$ in $\LAM$. Then
$$
\lim_{(\tau_{1},\tau_{2})\to(+i\infty,+i\infty)}(\varpi_{\langle J\rangle}^{*}\Phi_{\ell})(\tau_{1},\tau_{2})=2-\ell.
$$
\end{lemma}

\begin{pf}
Since $\Im z_{\langle J\rangle}/\Im\tau_{1}\in{\mathcal C}_{{\Bbb M}_{\ell}}^{+}$, 
we get $B=\lim_{\Im\tau_{1}\to+\infty}2\Im z_{\langle J\rangle}/\Im\tau_{1}\in\overline{\mathcal C}^{+}_{{\Bbb M}_{\ell}}$. 
Hence $\langle\lambda,B\rangle\geq0$ for all $\lambda\in\overline{\mathcal C}_{{\Bbb M}_{\ell}}^{+}$, 
where the equality holds if and only if $\lambda\in{\bf R}_{\geq0}B$. 
Similarly, we get $\langle\lambda,D\rangle\geq0$ for all $\lambda\in\overline{\mathcal C}_{{\Bbb M}_{\ell}}^{+}$,
where the equality holds if and only if $\lambda\in{\bf R}_{\geq0}D$. 
Since the isotropic vectors $B$, $D$ of the Lorentzian lattice ${\Bbb M}_{\ell}$ are not parallel and thus
$[\overline{\mathcal C}_{{\Bbb M}_{\ell}}^{+}\setminus{\bf R}_{\geq0}B]
\cup
[\overline{\mathcal C}_{{\Bbb M}_{\ell}}^{+}\setminus{\bf R}_{\geq0}D]
=
\overline{\mathcal C}_{{\Bbb M}_{\ell}}^{+}\setminus\{0\}$, we get
\begin{equation}
\label{eqn:limit:exponential:6.3}
\lim_{(\tau_{1},\tau_{2})\to(+i\infty,+i\infty)}
|e^{\pi i\langle\lambda,z_{\langle J\rangle}(\tau_{1},\tau_{2})\rangle}|
=
\lim_{(\tau_{1},\tau_{2})\to(+i\infty,+i\infty)}
e^{-\frac{\pi}{2}(\langle B,\lambda\rangle\Im\tau_{1}+\langle D,\lambda\rangle\Im\tau_{2})}
=
0
\end{equation}
for all $\lambda\in\overline{\mathcal C}_{{\Bbb M}_{\ell}}^{+}\setminus\{0\}$.
Since \eqref{eqn:Phi:expansion:level1}, \eqref{eqn:Phi:expansion:level2} converge absolutely 
for $\Im z\in{\Bbb M}_{\ell}\gg0$, we get the result by substituting \eqref{eqn:limit:exponential:6.3} 
into the explicit expressions \eqref{eqn:Phi:expansion:level1}, \eqref{eqn:Phi:expansion:level2}.
\end{pf}

Recall that the holomorphic $2$-form $\omega_{M(\lambda_{1},\lambda_{2})}$
on (the regular part of) $X_{(\lambda_{1},\lambda_{2})}$ was defined in Sect.\,\ref{sec:expression:Phi:special}.
For simplicity, write $\omega_{(\lambda_{1},\lambda_{2})}$ for $\omega_{M(\lambda_{1},\lambda_{2})}$.

\begin{lemma}
\label{lemma:norm:2-form:Kummer:product:type:(2,2,2)}
The following equality holds
\begin{equation}
\label{eqn:2-form:Kummer:product:type:(2,2,2)}
f^{*}\omega_{(\lambda_{1},\lambda_{2})}
=
2^{-3}\pi^{2}\theta_{3}(\tau_{1})^{2}\theta_{3}(\tau_{2})^{2}\,dz_{1}\wedge dz_{2}.
\end{equation}
\end{lemma}

\begin{pf}
In the coordinates of ${\bf P}^{5}$ given by $(x_{1},y_{0},x_{0},x_{2},y_{2},y_{1})$, 
$M(\lambda_{1},\lambda_{2})$ is of the form $(K,I_{3})$, $K\in M_{3}({\bf C})$
and the uniformization \eqref{eqn:uniformization} satisfies $x_{0}=1$ as required in Sect.\,\ref{sec:expression:Phi:special}.
Since $f^{*}dx_{1}=-\sqrt{e_{1}^{1}-e_{3}^{1}}\,{\rm sn}(z_{1}\sqrt{e_{1}^{1}-e_{3}^{1}})\,{\rm dn}(z_{1}\sqrt{e_{1}^{1}-e_{3}^{1}})\,dz_{1}$ and
$$
f^{*}dy_{0}
=
-\sqrt{e_{1}^{2}-e_{3}^{2}}\,
\frac{{\rm sn}(z_{1}\sqrt{e_{1}^{1}-e_{3}^{1}})}{{\rm sn}(z_{2}\sqrt{e^{2}_{1}-e^{2}_{3}})^{2}}
{\rm cn}(z_{2}\sqrt{e^{2}_{1}-e^{2}_{3}}){\rm dn}(z_{2}\sqrt{e^{2}_{1}-e^{2}_{3}})\,dz_{2}
\mod dz_{1}
$$
by the definitions of $x_{1},x_{2},y_{0},y_{1},y_{2}$, we deduce from \eqref{eqn:formula:2-form:X} the desired equality:
$$
f^{*}\omega_{(\lambda_{1},\lambda_{2})}
=
f^{*}(\frac{dx_{1}\wedge dy_{0}}{2^{3}x_{2}y_{1}y_{2}})
=
\frac{1}{2^{3}}\prod_{i=1,2}\sqrt{e^{i}_{1}-e^{i}_{3}}dz_{1}\wedge dz_{2}
=
\frac{\pi^{2}}{2^{3}}\theta_{3}(\tau_{1})^{2}\theta_{3}(\tau_{2})^{2}dz_{1}\wedge dz_{2},
$$
where the last equality follows from \cite[p.202 Eq.(6)]{HurwitzCourant29}.
\end{pf}

\begin{lemma}
\label{lemma:minor:Kummer:product:type:(2,2,2)}
Set $\Delta_{\langle J\rangle}:=\Delta_{\langle J\rangle}(M(\lambda_{1},\lambda_{2}))$.
Then
$$
\Delta_{\langle J\rangle}^{2}
=
\begin{cases}
\begin{array}{l}
(\lambda_{2}-1)^{2}
\\
\lambda_{2}^{2}
\\
1
\\
\lambda_{1}^{2}(\lambda_{2}-1)^{2}
\\
\lambda_{1}^{2}\lambda_{2}^{2}
\\
\lambda_{1}^{2}
\\
(\lambda_{1}-1)^{2}
\\
(\lambda_{1}-1)^{2}\lambda_{2}^{2}
\\
(\lambda_{1}-1)^{2}(\lambda_{2}-1)^{2}
\end{array}
\end{cases}
=
\begin{cases}
\begin{array}{l}
\frac{\theta_{0}(\tau_{2})^{8}}{\theta_{3}(\tau_{2})^{8}}
\\
\frac{\theta_{2}(\tau_{2})^{8}}{\theta_{3}(\tau_{2})^{8}}
\\
1
\\
\frac{\theta_{2}(\tau_{1})^{8}\theta_{0}(\tau_{2})^{8}}{\theta_{3}(\tau_{1})^{8}\theta_{3}(\tau_{2})^{8}}
\\
\frac{\theta_{2}(\tau_{1})^{8}\theta_{2}(\tau_{2})^{8}}
{\theta_{3}(\tau_{1})^{8}\theta_{3}(\tau_{2})^{8}}
\\
\frac{\theta_{2}(\tau_{1})^{8}}{\theta_{3}(\tau_{1})^{8}}
\\
\frac{\theta_{0}(\tau_{1})^{8}}{\theta_{3}(\tau_{1})^{8}}
\\
\frac
{\theta_{0}(\tau_{1})^{8}\theta_{2}(\tau_{2})^{8}}
{\theta_{3}(\tau_{1})^{8}\theta_{3}(\tau_{2})^{8}}
\\
\frac
{\theta_{0}(\tau_{1})^{8}\theta_{0}(\tau_{2})^{8}}
{\theta_{3}(\tau_{1})^{8}\theta_{3}(\tau_{2})^{8}}
\end{array}
\end{cases}
\hbox{ if }
\langle J\rangle
=
\begin{cases}
\begin{array}{l}
\binom{124}{356}
\\
\binom{125}{346}
\\
\binom{126}{345}
\\
\binom{134}{256}
\\
\binom{135}{246}
\\
\binom{136}{245}
\\
\binom{145}{236}
\\
\binom{146}{235}
\\
\binom{156}{234}.
\end{array}
\end{cases}
$$
\end{lemma}

\begin{pf}
By \eqref{eqn:matrix:Kummer:product:type:(2,2,2)}, 
the first equality is elementary. 
The second equality follows from the first one because $\lambda_{i}=\theta_{2}(\tau_{i})^{4}/\theta_{3}(\tau_{i})^{4}$ and
$1-\lambda_{i}=-\theta_{0}(\tau_{i})^{4}/\theta_{3}(\tau_{i})^{4}$.
\end{pf}

\begin{lemma}
\label{lemma:period:integral:Kummer:product:type:(2,2,2)}
Let $g\colon K_{(\tau_{1},\tau_{2})}\to X_{(\lambda_{1},\lambda_{2})}$ be the minimal resolution.
If $\langle J\rangle\not=\binom{123}{456}$, then
$$
\Delta_{\langle J\rangle}^{2}\cdot\left(\int_{\Gamma_{34}}\frac{2}{\pi^{2}}g^{*}\omega_{(\lambda_{1},\lambda_{2})}\right)^{4}
=
\theta_{\epsilon}(\tau_{1})^{8}\theta_{\delta}(\tau_{2})^{8},
$$
where the correspondence between $(\epsilon,\delta)$ and $\langle J\rangle$ is given as follows:
$$
\begin{array}{cccccccccccc}
(\epsilon,\delta)
&
\vline
&
(2,2)
&
(2,0)
&
(2,3)
&
(0,2)
&
(0,0)
&
(0,3)
&
(3,2)
&
(3,0)
&
(3,3)
\\
\hline
\langle J\rangle
&
\vline
&
\binom{135}{246}
&
\binom{134}{256}
&
\binom{136}{245}
&
\binom{146}{235}
&
\binom{156}{234}
&
\binom{145}{236}
&
\binom{125}{346}
&
\binom{124}{356}
&
\binom{126}{345}
\end{array}
$$
\end{lemma}

\begin{pf}
The cycles $\beta_{i}\in H_{1}(E_{\tau_{i}},{\bf Z})$ is represented by the circle
$({\bf R}+\epsilon\sqrt{-1})/2{\bf Z}$ with $0<\epsilon\ll1$, which does not pass through any points of order $2$ of $E_{\tau_{i}}$.
Regard $f^{*}\omega_{(\lambda_{1},\lambda_{2})}$ as the $2$-form on $A_{(\tau_{1},\tau_{2})}$ by the formula
\eqref{eqn:2-form:Kummer:product:type:(2,2,2)}.
Since $\phi(f^{*}\omega_{(\lambda_{1},\lambda_{2})})=2g^{*}\omega_{(\lambda_{1},\lambda_{2})}$ 
and since ${\bf b}_{1}^{\lor}\wedge{\bf b}_{2}^{\lor}$ is the Poincar\'e dual of $\beta_{1}\times\beta_{2}\subset A_{(\tau_{1},\tau_{2})}$, 
we get 
$$
\begin{aligned}
\int_{\Gamma_{34}}g^{*}\omega_{(\lambda_{1},\lambda_{2})}
&=
\frac{1}{2}
\left\langle
\phi(f^{*}\omega_{(\lambda_{1},\lambda_{2})}),\Gamma_{34}^{\lor}
\right\rangle
=
\left\langle
f^{*}\omega_{(\lambda_{1},\lambda_{2})},{\bf b}_{1}^{\lor}\wedge{\bf b}_{2}^{\lor}
\right\rangle
=
\int_{\beta_{1}\times\beta_{2}}f^{*}\omega_{(\lambda_{1},\lambda_{2})}
\\
&=
\int_{\beta_{1}\times\beta_{2}}\frac{\pi^{2}}{2^{3}}\theta_{3}(\tau_{1})^{2}\theta_{3}(\tau_{2})^{2}\,dz_{1}\wedge dz_{2}
=
\frac{\pi^{2}}{2}\theta_{3}(\tau_{1})^{2}\theta_{3}(\tau_{2})^{2}
\end{aligned}
$$
by Lemma~\ref{lemma:norm:2-form:Kummer:product:type:(2,2,2)} and \eqref{eqn:definition:symplectic:basis}.
This, together with Lemma~\ref{lemma:minor:Kummer:product:type:(2,2,2)}, implies the result.
\end{pf}

\begin{theorem}
\label{theorem:constant:Thomae:Kummer:product:type}
One has ${\frak c}_{1}={\frak c}_{2}={\frak C}=1$ in \eqref{eqn:constant:Phi:generic} and Proposition~\ref{proposition:Thomae:Phi:with:constant}.
\end{theorem}

\begin{pf}
Since ${\frak C}=|{\frak c}_{\ell}|$, it suffices to prove ${\frak c}_{1}={\frak c}_{2}=1$.
Let $\langle J\rangle$ be a partition such that $\alpha_{\langle J\rangle}(\Gamma_{34}^{\lor})={\frak e}_{\ell}$. Then we get
\begin{equation}
\label{eqn:constant:Phi:theta}
\sqrt{{\frak c}_{\ell}}
=
\frac{\varpi_{\langle J\rangle}^{*}\Phi_{\ell}}{\Delta_{\langle J\rangle}^{2}\cdot\left(\int_{\Gamma_{34}}2\pi^{-2}\omega_{(\lambda_{1},\lambda_{2})}\right)^{4}}
=
\frac{(\varpi_{\langle J\rangle}^{*}\Phi_{\ell})(\tau_{1},\tau_{2})}{\theta_{\epsilon}(\tau_{1})^{8}\theta_{\delta}(\tau_{2})^{8}},
\end{equation}
where the first equality follows from Proposition~\ref{proposition:Thomae:Phi:with:constant} and the second from
Lemma~\ref{lemma:period:integral:Kummer:product:type:(2,2,2)}.
By \eqref{eqn:infinite:product:theta:constant:genus:1}, we get
$\lim_{(\tau_{1},\tau_{2})\to(+i\infty,+i\infty)}\theta_{\epsilon}(\tau_{1})^{8}\theta_{\delta}(\tau_{2})^{8}=0$ or $1$ and both cases occur.
Since ${\frak c}_{\ell}$ is a {\em non-zero} constant on ${\frak H}\times{\frak H}$ by Proposition~\ref{proposition:Thomae:Phi:with:constant},
it follows from \eqref{eqn:constant:Phi:theta} and Lemma~\ref{lemma:limit:Phi:level:1:cusp}
that for any $\ell\in\{1,2\}$, there is a choice $\langle J\rangle$ such that
$\alpha_{\langle J\rangle}(\Gamma^{\lor}_{34})$ has level $\ell$ in $\LAM$. Hence $\alpha_{\langle J\rangle}(\Gamma_{34}^{\lor})={\frak e}_{\ell}$.
\par
Let $\ell=1$.
By \eqref{eqn:constant:Phi:theta} and Lemma~\ref{lemma:limit:Phi:level:1:cusp}, 
we get $\lim_{(\tau_{1},\tau_{2})\to(+i\infty,+i\infty)}\sqrt{{\frak c}_{1}}=\pm1$,
since ${\frak c}_{1}$ is a non-zero constant by Proposition~\ref{proposition:Thomae:Phi:with:constant}.
This proves that ${\frak c}_{1}=1$.
\par
Assume $\ell=2$. Set $q_{m}:=\exp(\pi i\tau_{m})$.
Since $\Phi_{2}(z)$ has integral Fourier coefficients by \eqref{eqn:Phi:expansion:level2} and thus
$(\varpi_{\langle J\rangle}^{*}\Phi_{\ell})(\tau_{1},\tau_{2})\in{\bf Z}\{q_{1},q_{2}\}$,
we get $\sqrt{{\frak c}_{2}}\in q_{1}^{-2}q_{2}^{-2}\,{\bf Z}\{q_{1},q_{2}\}$ by \eqref{eqn:infinite:product:theta:constant:genus:1}, \eqref{eqn:constant:Phi:theta}.
Hence $\sqrt{{\frak c}_{2}}\in{\bf Z}$.
Since $|{\frak c}_{2}|=|{\frak c}_{1}|$ and since ${\frak c}_{1}=1$, we get ${\frak c}_{2}=1$.
\end{pf}

\begin{corollary}
\label{cor:value:Phi:Kummer:product:type:(2,2,2):1}
If $\alpha_{\langle J\rangle}(\Gamma_{34}^{\lor})$ is a primitive isotropic vector of $\LAM$ of level $\ell$, then
$$
\Phi_{\ell}\left(\varpi_{\langle J\rangle}(\tau_{1},\tau_{2})\right)=\pm\theta_{\epsilon}(\tau_{1})^{8}\theta_{\delta}(\tau_{2})^{8}
$$
under the correspondence between $\langle J\rangle$ and $(\epsilon,\delta)$ in Lemma~\ref{lemma:period:integral:Kummer:product:type:(2,2,2)}.
In particular, 
$$
\left\|
\Phi(Y_{\langle J\rangle}(\lambda_{1},\lambda_{2}))
\right\|
=
\|\theta_{\epsilon}(\tau_{1})\theta_{\delta}(\tau_{2})\|^{8}.
$$
\end{corollary}

\begin{pf}
Since ${\frak c}_{1}={\frak c}_{2}=1$ by Theorem~\ref{theorem:constant:Thomae:Kummer:product:type}, 
the result follows from \eqref{eqn:constant:Phi:theta}.
\end{pf}


\end{document}